\documentclass[a4paper,10pt]{article}
\usepackage[a4paper, margin=1.2in]{geometry}
\usepackage{mathrsfs,mathtools,amsfonts,latexsym,theorem,amssymb,amsmath,newlfont}
\usepackage{graphicx}

\usepackage{pgf}
\usepackage{fancyhdr}
\usepackage{times}
\usepackage{dsfont} 
\usepackage{accents} 
\usepackage{tikz} 
\usetikzlibrary{calc} 
\usetikzlibrary{decorations.pathreplacing,decorations.markings} 

\tikzset{
  on each segment/.style={
    decorate,
    decoration={
      show path construction,
      moveto code={},
      lineto code={
        \path [#1]
        (\tikzinputsegmentfirst) -- (\tikzinputsegmentlast);
      },
      curveto code={
        \path [#1] (\tikzinputsegmentfirst)
        .. controls
        (\tikzinputsegmentsupporta) and (\tikzinputsegmentsupportb)
        ..
        (\tikzinputsegmentlast);
      },
      closepath code={
        \path [#1]
        (\tikzinputsegmentfirst) -- (\tikzinputsegmentlast);
      },
    },
  },
  mid arrow/.style={postaction={decorate,decoration={
        markings,
        mark=at position .6 with {\arrow[#1]{stealth}}
      }}},
}

\tikzset{
    right angle quadrant/.code={
        \pgfmathsetmacro\quadranta{{1,1,-1,-1}[#1-1]}     
        \pgfmathsetmacro\quadrantb{{1,-1,-1,1}[#1-1]}},
    right angle quadrant=1, 
    right angle length/.code={\def\rightanglelength{#1}},   
    right angle length=2ex, 
    right angle symbol/.style n args={3}{
        insert path={
            let \p0 = ($(#1)!(#3)!(#2)$) in     
                let \p1 = ($(\p0)!\quadranta*\rightanglelength!(#3)$), 
                \p2 = ($(\p0)!\quadrantb*\rightanglelength!(#2)$) in 
                let \p3 = ($(\p1)+(\p2)-(\p0)$) in  
            (\p1) -- (\p3) -- (\p2)
        }
    }
}

\theoremstyle{plain}

\newtheorem{thm}{Theorem}[section]
\newtheorem{defn}[thm]{Definition}
\newtheorem{lem}[thm]{Lemma}
\newtheorem{prop}[thm]{Proposition}
\newtheorem{cor}[thm]{Corollary}
\newtheorem{remark}[thm]{Remark}

\newtheorem{conj}[thm]{Conjecture}

\numberwithin{equation}{section}

\def\beginpf{\noindent {\bf Proof:} \quad}
\def\endpf{\rightline{$\square$}}

\def\H{{\mathcal H}}
\def\E{{\mathcal E}}

\def\NN^*{{\mathbb N}^*}
\def\NN{{\mathbb N}}

\def\RR{{\mathbb R}}

\def\ZZ{{\mathbb Z}}
\def\M{{\cal M}}

\def\F{{\cal F}}

\def\G{{\cal G}}

\renewcommand{\d}{\partial}

\def\<{\langle}
\def\>{\rangle}

\def\restriction#1#2{\mathchoice
              {\setbox1\hbox{${\displaystyle #1}_{\scriptstyle #2}$}
              \restrictionaux{#1}{#2}}
              {\setbox1\hbox{${\textstyle #1}_{\scriptstyle #2}$}
              \restrictionaux{#1}{#2}}
              {\setbox1\hbox{${\scriptstyle #1}_{\scriptscriptstyle #2}$}
              \restrictionaux{#1}{#2}}
              {\setbox1\hbox{${\scriptscriptstyle #1}_{\scriptscriptstyle #2}$}
              \restrictionaux{#1}{#2}}}
\def\restrictionaux#1#2{{#1\,\smash{\vrule height .8\ht1 depth .85\dp1}}_{\,#2}}

\allowdisplaybreaks[3]        

\title{Sufficient Conditions for the Controllability of Wave Equations with a Transmission Condition at the Interface}
\author{ Ludovick \textsc{Gagnon}\footnote{ This work was supported by the European Research Council, 
ERC-2012-ADG, project number 320845:  Semi Classical Analysis of Partial Differential
Equations. Address : Universit\'e de Nice Sophia-Antipolis, CNRS UMR 7351, Laboratoire J.-A. Dieudonn\'e, Parc Valrose, 06108, Nice, France. E-mail address : gagnonl@unice.fr}}

\begin{document}

\maketitle

\begin{abstract}

We consider waves travelling in two different mediums each endowed with a different constant speed of propagation. At the interface between the two mediums, the refraction of the rays of the optic geometry is described by the Snell-Descartes law. We provide sufficient conditions on the geometry of the mediums and on the speed of propagation for the boundary controllability.

\end{abstract}

\tableofcontents

\section{Introduction}

\subsection{Statement of the problem}\label{SecState}

We consider the observability of wave equations with a transmission condition at the interface. Let $\Omega \subset \RR^2$ be an open, bounded and strictly convex domain. Let $\Omega_2 \subset \Omega$ be an open, bounded, strictly convex domain such that there exists $\Omega' \subset \Omega$ be open and simply connected such that $\overline{\Omega}_2 \subset \Omega'$ and $\overline{\Omega'}\subset \Omega$. We define $\Omega_1:=\Omega \setminus \overline{\Omega}_2$. In this setting, the boundary of $\Omega_2$ is $\d \Omega_2$ and the boundary of $\Omega_1$ is $\d \Omega_1=\d \Omega \cup \d \Omega_2$. We assume the boundary $\d \Omega_2$ and $\d \Omega$ to be of class $C^k, k\geq 3$ and with no contact of order $k-1$ with its tangents. The outward unit normal of $x\in \d \Omega_2$ is denoted $n_2(x)$ and the outward unit normal of $x\in \d \Omega$ is denoted $n(x)$. With these notations, we identify the outward unit normal of $\Omega_1$ by $-n_2(x)$ if $x\in \d \Omega_1 \cap \d \Omega_2$ and $n(x)$ if $x\in \d \Omega_1 \cap \d \Omega$. To ease the notations, we identify $\d \Omega_1 \cap \d \Omega$ with $\d \Omega$ and $\d \Omega_1 \cap \d \Omega_2$ with $\d \Omega_2$.   

\begin{figure}[!h]
\begin{center}
\begin{tikzpicture}
\draw (0,0) circle (1cm);
\node[text width=.25cm] at (0,-.55) {$u_2$};
\node[text width=.25cm] at (0,.25) {$\Omega_2$};
\node[text width=.25cm] at (0,-1.52) {$u_1$};
\draw (0,0) circle (2cm);
\node[text width=.25cm] at (-1.6,0) {$\d \Omega_2$};
\node[text width=.25cm] at (0,1.5) {$\Omega_1$};
\node[text width=.25cm] at (-2.52,0) {$\d \Omega$};
\draw[->](1,0) -- (1.6,0) node[midway,above] {$n_2$};
\draw[->](2,0) -- (2.6,0) node[midway,above] {$n$};
\end{tikzpicture}
\end{center}
\caption{Representation of the spatial domain for \eqref{wave}-\eqref{Dirhom} with ${\Omega_2=\{x\in \RR^2 \, | \, \|x\| < 1\}}$, \mbox{$\Omega_1=\{x\in \RR^2 \, | \, 1<\|x\|< 2\}$}.}\label{repressol}
\end{figure}
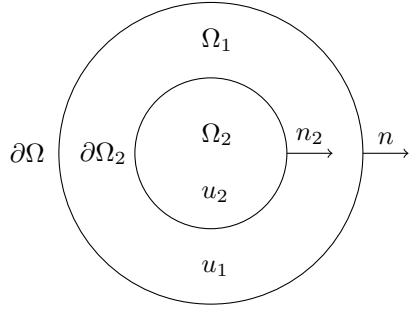

For $T>0$, we consider the set of wave equations
\begin{equation}\label{wave}
\begin{cases}
(\d_t^2-c_i^2 \Delta )u^i(t,x)=0, & (t,x)\in (0,T)\times \Omega_i, \\
u^i(0,x)=u_0^i(x), u_t^i(0,x)=u_1^i(x), & x\in \Omega_i
\end{cases}
\end{equation}
where $c_i >0, i=1,2$ and $(\mathds{1}_{\Omega_1} u_0^1+\mathds{1}_{\Omega_2} u_0^2,\mathds{1}_{\Omega_1} u_1^1+\mathds{1}_{\Omega_2} u_1^2)\in H^1_0(\Omega) \times L^2(\Omega)$. At the interface $\d \Omega_2$, a transmission condition is imposed, 
\begin{equation}\label{transmission} 
\begin{cases}
u^1(t,x)=u^{2}(t,x), & (t,x)\in (0,T)\times \d \Omega_2,\\
c_1^2 \d_{n_2} u^1(t,x)=c_{2}^2  \d_{n_2} u^{2}(t,x), & (t,x)\in (0,T)\times \d \Omega_2.
\end{cases}
\end{equation} 
At the exterior boundary $\d \Omega$, an homogeneous Dirichlet boundary condition is imposed,
\begin{equation}\label{Dirhom}
u^1(t,x)=0, \qquad (t,x)\in (0,T)\times \d \Omega.
\end{equation}
Figure \ref{repressol} illustrates the spatial domain on which the solutions to \eqref{wave}-\eqref{Dirhom} are defined.

The energy of solutions of \eqref{wave}-\eqref{Dirhom} is preserved along time,
\[
E(u^1(t,.),u^2(t,.)):=\sum_{i=1}^2 \int_{\Omega_i} \left| \d_t u^i(t,x) \right|^2+ c_i^2 \left| \nabla u^i(t,x) \right|^2 \, \textrm{dx}=E(u^1(0,.),u^2(0,.)).
\] 

The definition of the observability of solutions to \eqref{wave}-\eqref{Dirhom} is stated as follow.

\begin{defn}\label{defnobs}
Let $\Gamma \subset \d \Omega$. We say that \eqref{wave}-\eqref{Dirhom} is observable in time $T>0$ with respect to $\Gamma$ if and only if there exists $c_T>0$ such that for any initial data $(\mathds{1}_{\Omega_1} u_0^1+\mathds{1}_{\Omega_2} u_0^2,\mathds{1}_{\Omega_1} u_1^1+\mathds{1}_{\Omega_2} u_1^2)\in H^1_0(\Omega) \times L^2(\Omega)$, the solutions to \eqref{wave}-\eqref{Dirhom} satisfy
\begin{equation}\label{ineqobs}
E(u^1(0,.),u^2(0,.))\leq c_T \int_0^T \int_{\Gamma} \left| \d_{n}u^1(t,x)\right|^2 \, \textrm{d}\sigma \textrm{dt}.
\end{equation}
\end{defn}
The observability of \eqref{wave}-\eqref{Dirhom} finds applications in optical fibers and guided waves \cite{Saleh} and also in seismic prospection of Earth inner layers \cite{Hoop}.

We recall that the observability of \eqref{wave}-\eqref{Dirhom} is equivalent (\cite{Brezis}) to the controllability of 
\begin{equation*}
\begin{cases}
(\d_t^2-c_i^2 \Delta )y^i(t,x)=0, & (t,x)\in (0,T)\times \Omega_i, \\
y^1(t,x)=y^{2}(t,x), & (t,x)\in (0,T)\times \d \Omega_2,\\
c_1^2 \d_{n_2} y^1(t,x)=c_{2}^2  \d_{n_2} y^{2}(t,x), & (t,x)\in (0,T)\times \d \Omega_2, \\ 
y^1(t,x)=0, & (t,x)\in (0,T)\times \d \Omega\setminus \Gamma, \\
y^1(t,x)=v(x,t), & (t,x)\in (0,T)\times \Gamma, \\
y^i(0,x)=y_0^i(x), y_t^i(0,x)=y_1^i(x), & x\in \Omega_i 
\end{cases}
\end{equation*}
with $(\mathds{1}_{\Omega_1} y_0^1+\mathds{1}_{\Omega_2} y_0^2,\mathds{1}_{\Omega_1} y_1^1+\mathds{1}_{\Omega_2} y_1^2)\in L^2(\Omega) \times H^{-1}(\Omega)$ and $v\in L^2((0,T)\times \Gamma)$.


\subsection{The problematic}\label{SecStrategy}

Let us present, formally, the main difficulty in the analysis of the observability of \eqref{wave}-\eqref{Dirhom}. In order to do so, we begin by recalling the main ideas behind the proof of the observability for the classical wave equation with the microlocal analysis approach. For $\Omega$ defined as above and $T>0$, consider
\begin{equation}\label{ClassicWave}
\begin{cases}
(\d_{tt}-c^2 \Delta)u(x,t) =0 , & (x,t) \in \Omega \times (0,T), \\
u(x,t)=0, & (x,t) \in \d \Omega \times (0,T), \\
u(x,0)=u_0(x), u_t(x,0)=u_1(x), & x \in \Omega,
\end{cases}
\end{equation}
where $c>0$, $(u_0,u_1)\in H^1_0(\Omega) \times L^2(\Omega)$. We recall that the energy 
\[
E(u(t,.)):=\int_{\Omega} \left| \d_t u(t,x) \right|^2+ c^2 \left| \nabla u(t,x) \right|^2 \, \textrm{dx},
\] 
is preserved along time. The observability inequality for \eqref{ClassicWave} is defined as for \eqref{wave}-\eqref{Dirhom},
\begin{equation}\label{obsineclassical}
E(u)(0)\leq c_T \int_0^T \int_{\Gamma} \left| \d_{n}u(t,x)\right|^2 \, \textrm{d}\sigma \textrm{dt}.
\end{equation}
The proof the observability inequality \eqref{obsineclassical} is done by contradiction and its analysis is reduced to the observability of the rays of the optic geometry of $\Omega$. Rays of \eqref{ClassicWave} travel at constant speed in straight line and are reflected on $\d \Omega$ with the rule "the angle of reflection is equal to the angle of incidence". The observability is obtained in this framework if the \textit{Geometric Control Condition} \cite{BLR} (GCC in short) is satisfied, that is, if every ray of the optic geometry in $\Omega$ intersects (in a non-diffractive way) $\Gamma$ in time $T>0$. 

\begin{figure}[!ht]
\begin{center}
	\includegraphics[height=5cm]{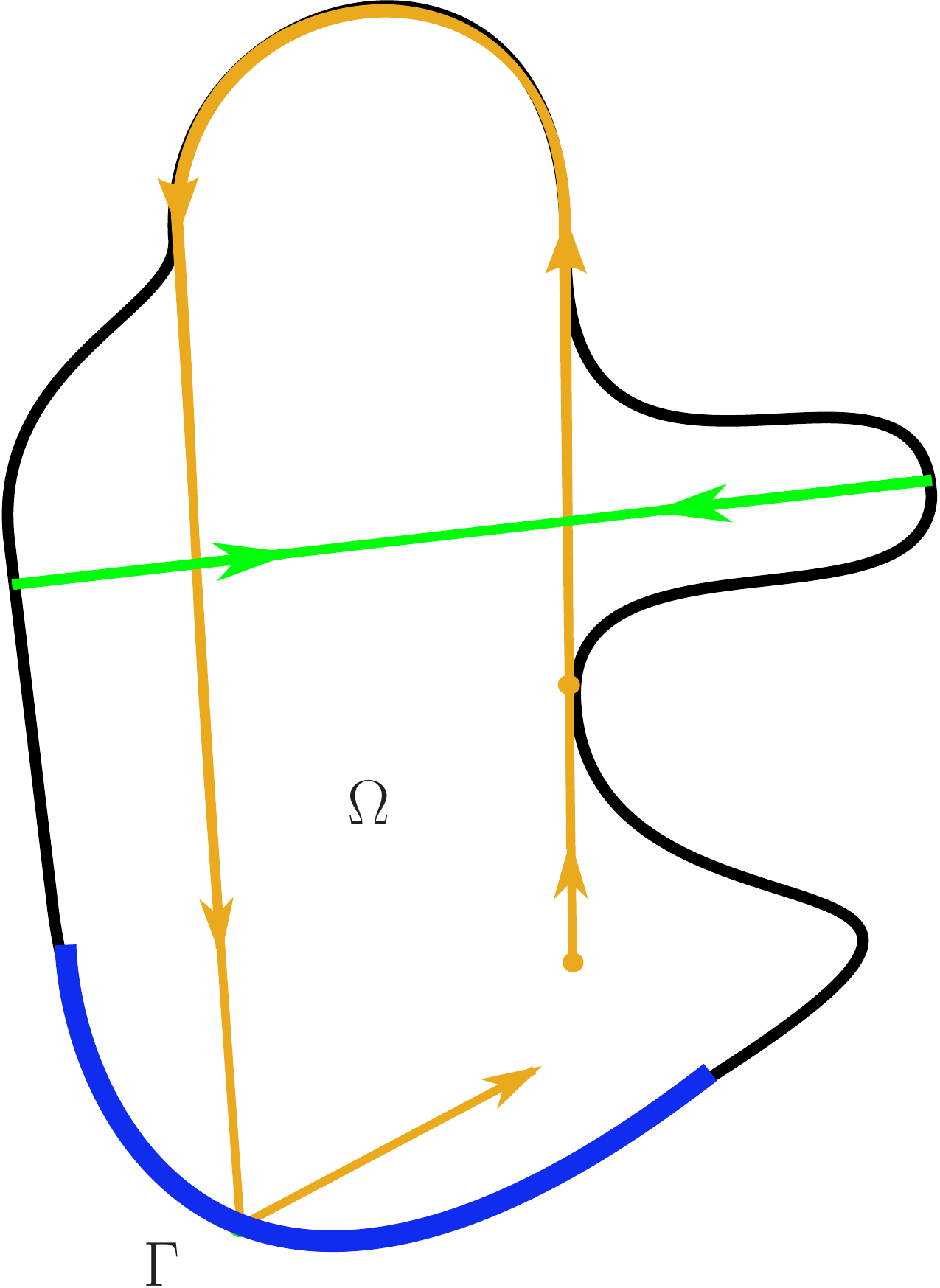}
\end{center}
\caption{\footnotesize{Representation of the propagation of rays for \eqref{ClassicWave}. The blue region is the observability region $\Gamma \subset \d \Omega$. The ray in yellow is a ray observed by the observability region. The ray in green is a trapped ray for any time $T>0$ for $(\Omega,\Gamma)$.}}
\end{figure}

\begin{figure}[!ht]
\begin{center}
	\includegraphics[height=5cm]{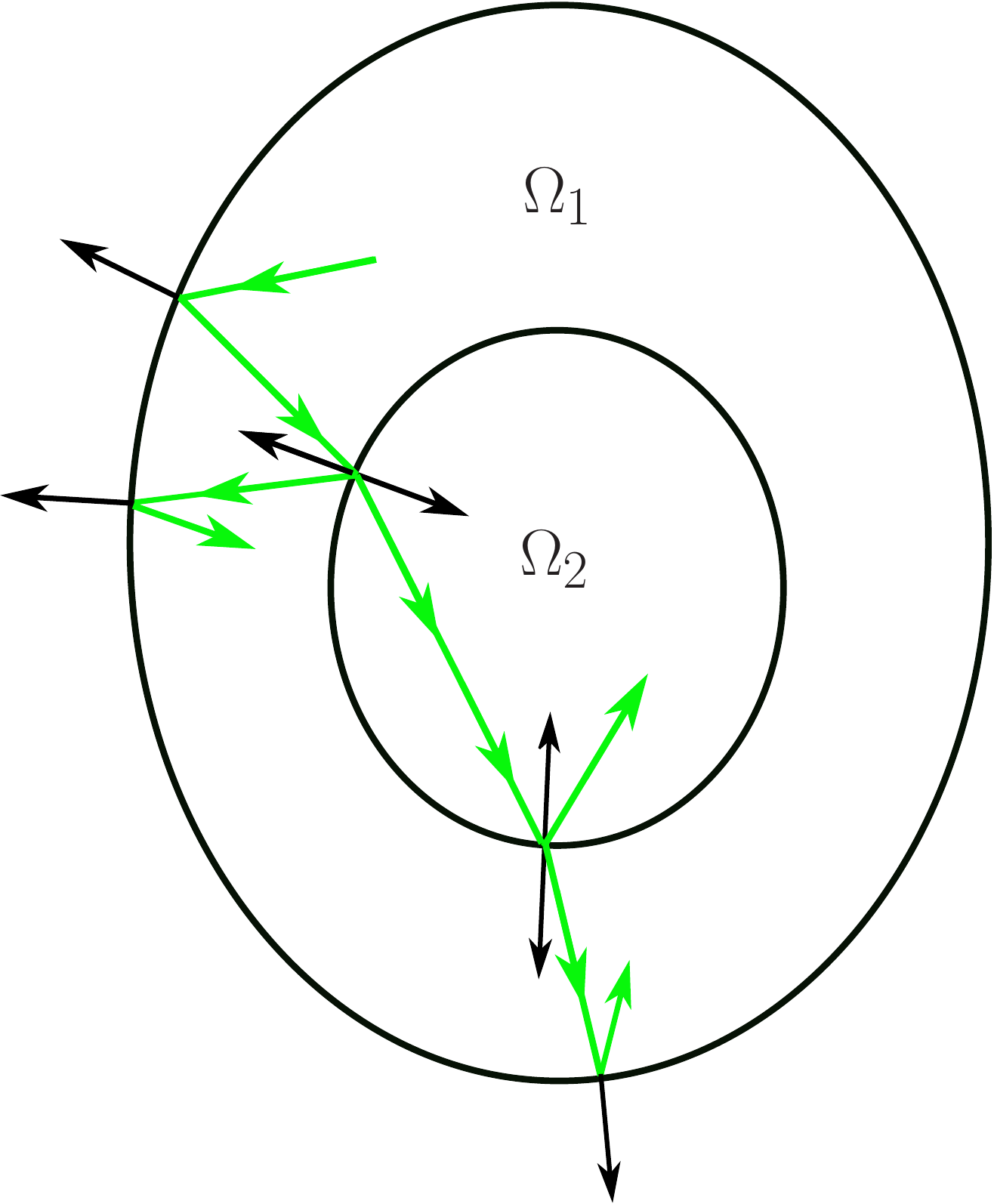}
\end{center}
\caption{\footnotesize{Example of the propagation of rays for \eqref{wave}-\eqref{Dirhom}}}\label{ExPropag}
\end{figure}

It will be shown in Section \ref{SecDefBichar} that rays of \eqref{wave}-\eqref{Dirhom} also travel in straight line at constant speed $c_i$ in $\Omega_i, i=1,2$ and rays in $\Omega_1$ are reflected on $\d \Omega$ according to law of reflection of \eqref{ClassicWave}. Moreover, if a ray, say $\gamma_1^-$, intersects $\d \Omega_2$ with an angle $\theta_1$ (with respect to the outward normal $n_2$), then a ray $\gamma_1^+$ will be reflected in $\Omega_1$ with an angle $\theta_1$. If the Snell-Descartes 
\begin{equation}\label{SnellDescartes}
\dfrac{\sin \theta_1}{c_1} = \dfrac{\sin \theta_2}{c_2},
\end{equation} 
is not vacuous for the values of $c_1,c_2$ and $\theta_1$, then a ray $\gamma_2^+$ will be transmitted in $\Omega_2$ with an angle $\theta_2$ (with respect to the inward normal $-n_2$). If \eqref{SnellDescartes} is vacuous, then $\Omega_2$ acts as an obstacle for the ray $\gamma_1^-$ and only a reflection occurs at the boundary $\d \Omega_2$ for $\gamma_1^-$. The same description holds for rays coming from $\Omega_2$. We stress here that \eqref{SnellDescartes} is never vacuous under the hypothesis $c_2>c_1$ for incoming rays on $\d \Omega_2$ of $\Omega_2$, meaning that every rays from $\Omega_2$ are reflected and transmitted on $\d \Omega_2$. Figure \ref{ExPropag} illustrate an example of propagation of a ray for \eqref{wave}-\eqref{Dirhom}.

\begin{figure}[!ht]
\begin{center}
	\includegraphics[height=5cm]{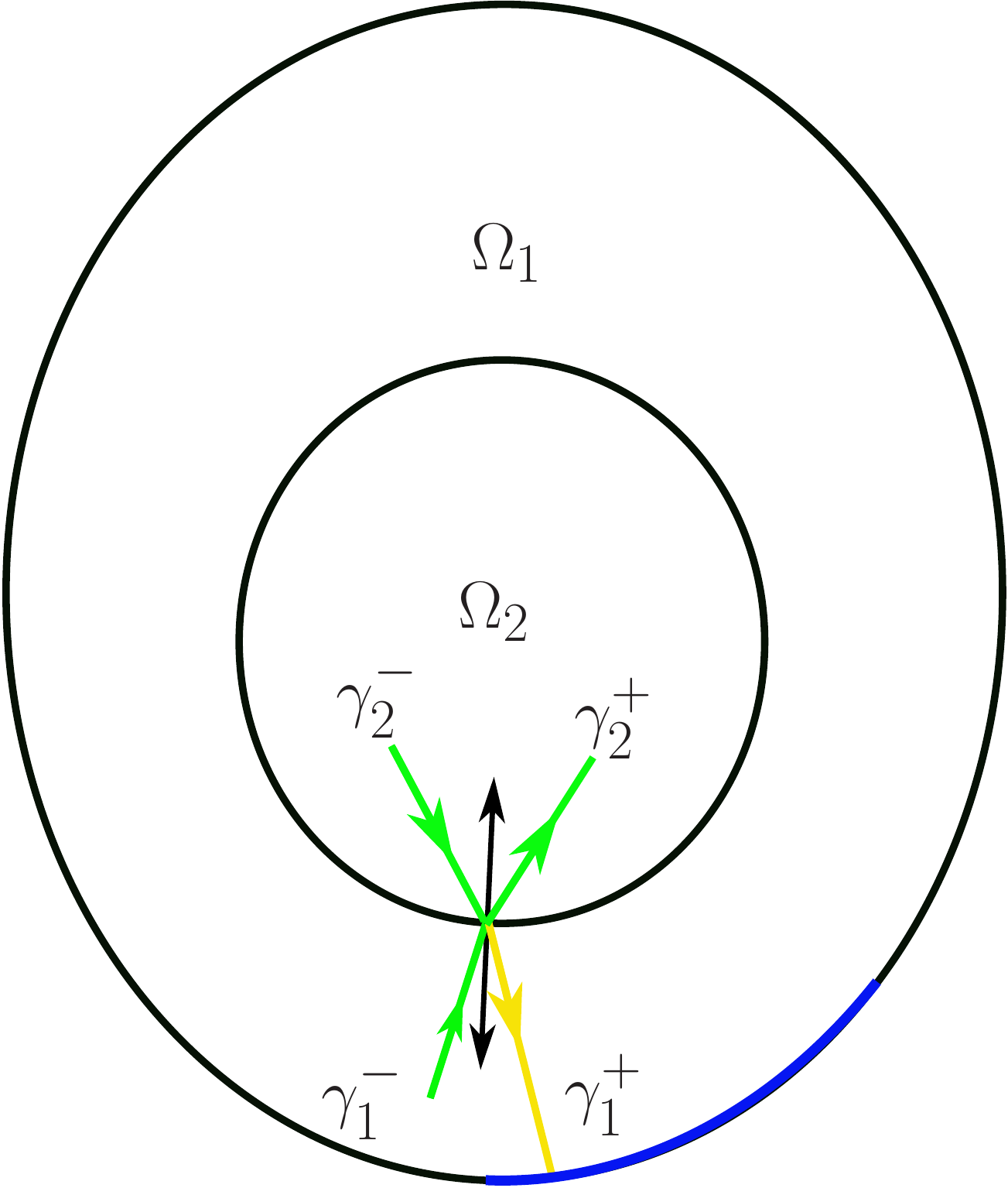}
\end{center}
\caption{\footnotesize{Example of a ray obtained as a superposition of two rays at $\d \Omega_2$ for \eqref{wave}-\eqref{Dirhom}}}\label{Exsuperpo}
\end{figure}

In this framework, unlike the wave equation \eqref{ClassicWave}, the observation of a single ray at $\Gamma$ is not sufficient to conclude on the observability of the starting point of the ray. Indeed, denote $\gamma_1^+$ a ray intersecting $\Gamma$ at time $t$ and suppose this ray have intersected $\d \Omega_2$ at a time $t'<t$. Denote $\theta_1$ the angle between $\gamma_1^+$ and $n_2$ at the point of intersection on $\d \Omega_2$. This ray may be the outgoing ray from a reflection of a ray $\gamma_1^-$ coming from $\Omega_1$ at angle $\theta_1$ or the transmitted ray of a ray $\gamma_2^-$ coming from $\Omega_2$ at an angle $\theta_2$ such that $\theta_1$ and $\theta_2$ satisfy \eqref{SnellDescartes}. From the superposition principle, $\gamma_1^+$ may be obtained simultaneously as a reflection from $\Omega_1$ and a transmission from $\Omega_2$ at time $t'$ (see figure \ref{Exsuperpo}).

\begin{figure}[!ht]
\begin{center}
	a)\includegraphics[height=5cm]{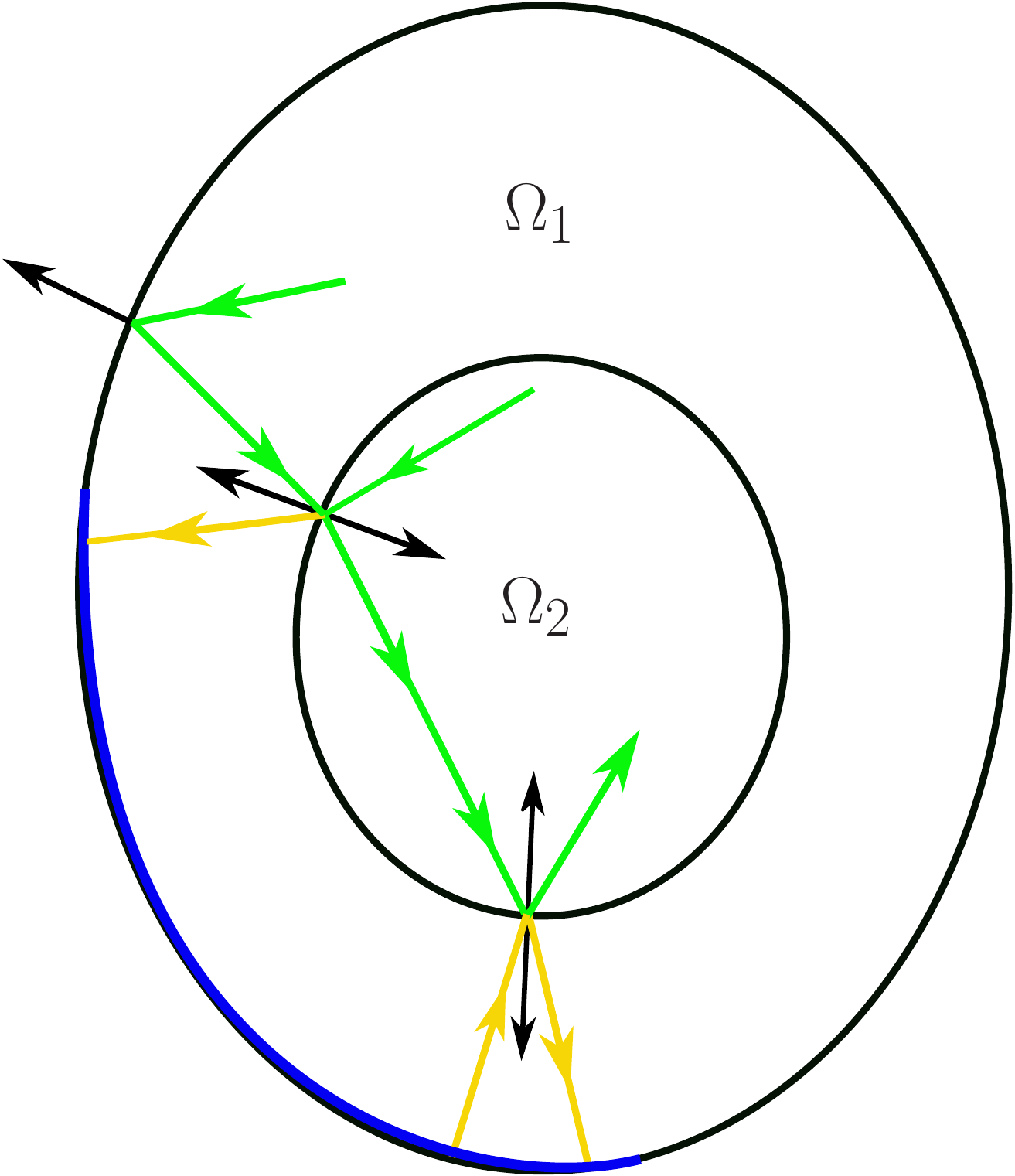} \quad 
		b)\includegraphics[height=5cm]{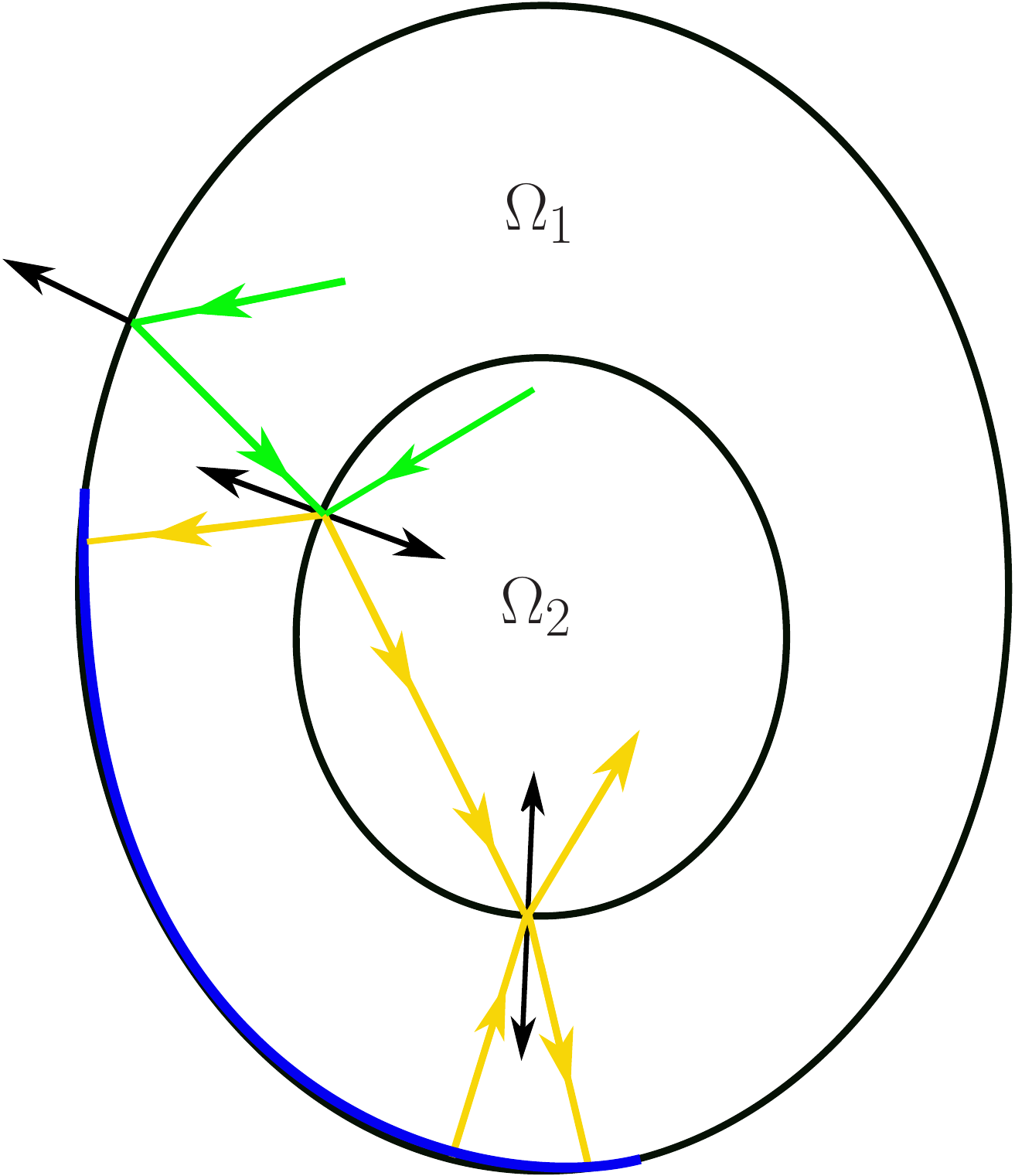} \\ 
		\vspace{.4cm}
				c)\includegraphics[height=5cm]{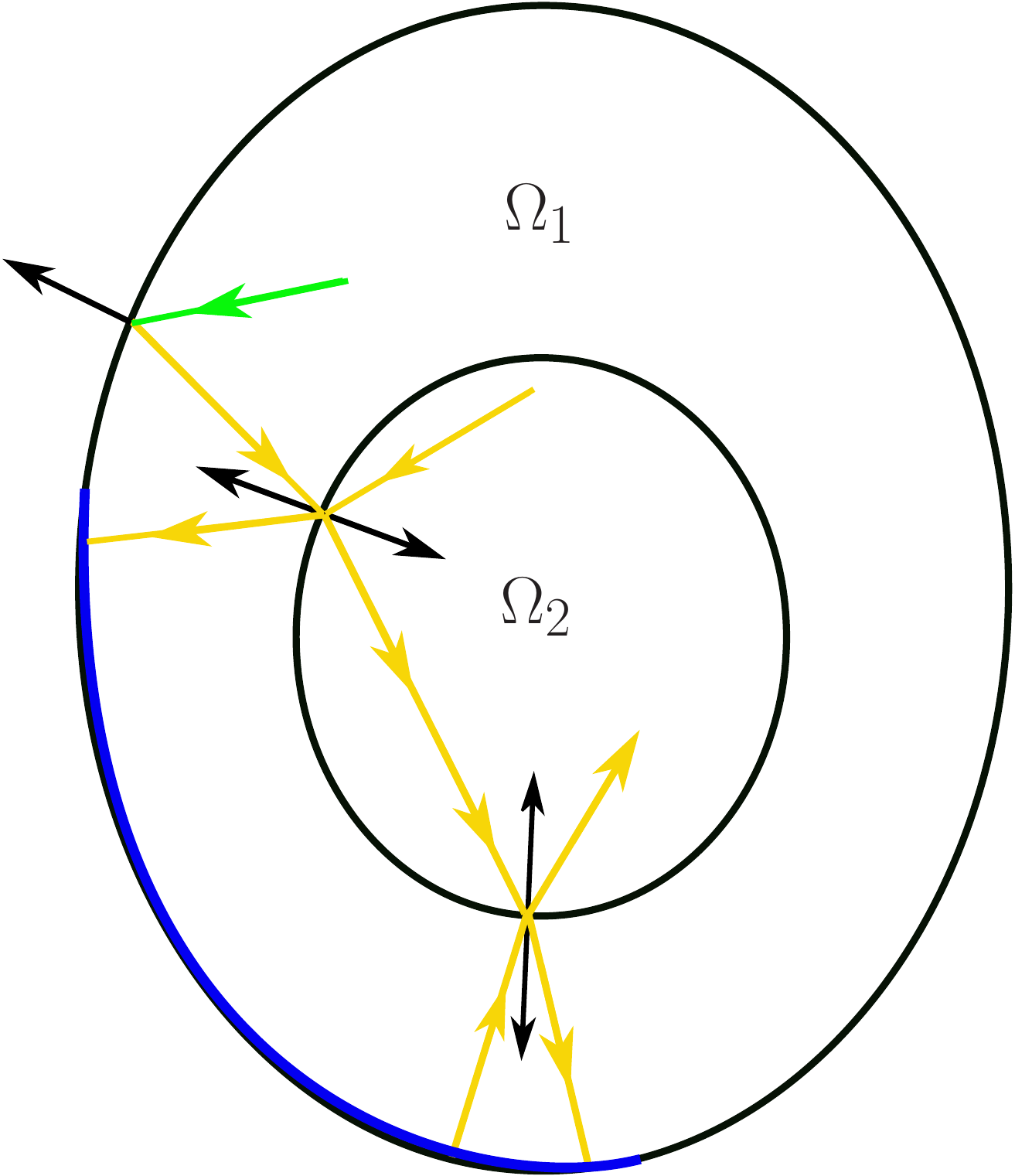} \quad 
				d)\includegraphics[height=5cm]{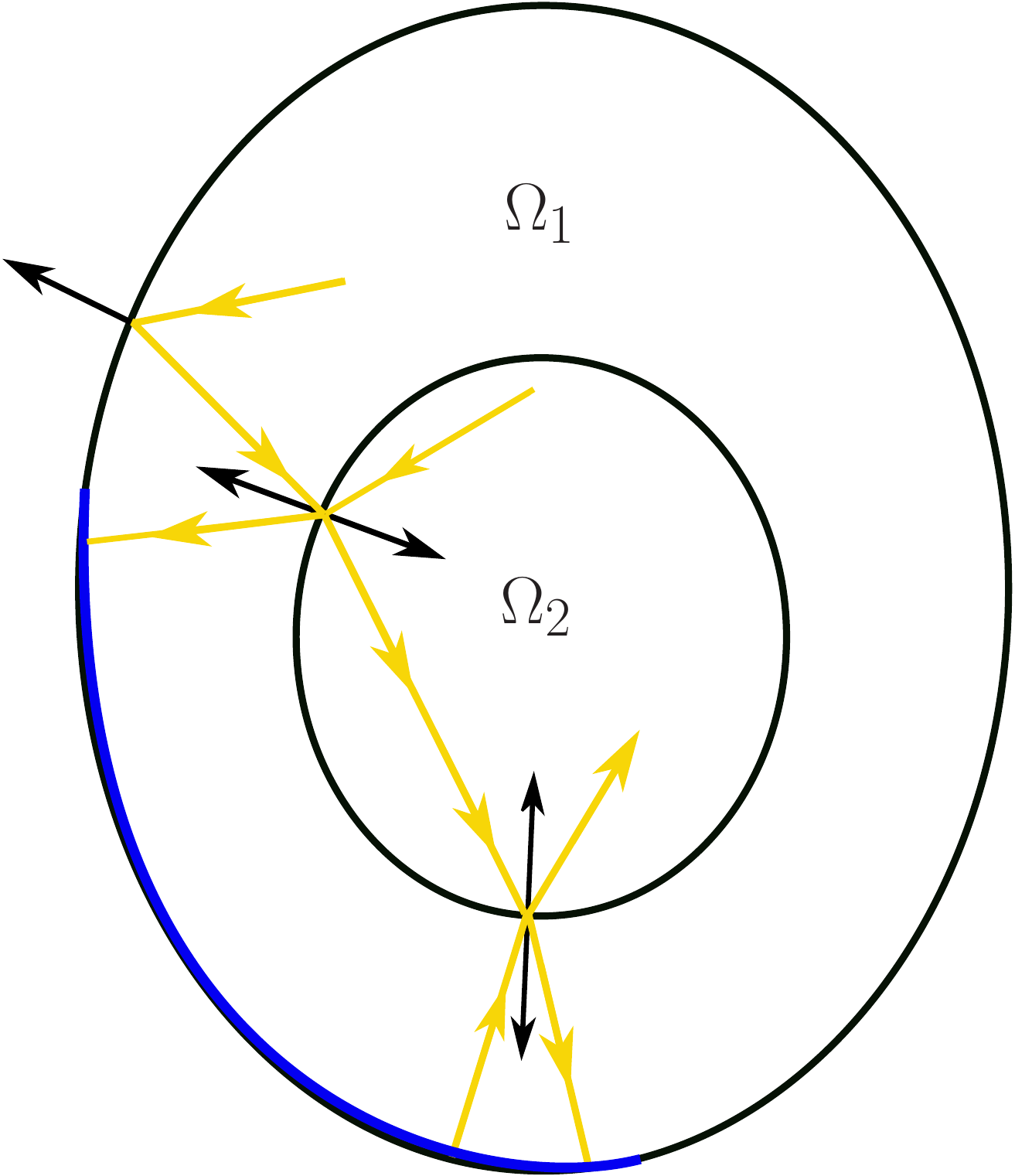} 
\end{center}
\caption{\footnotesize{Recursive use of Corollary \ref{equivsuppmes} and the propagation for the classical wave equation at the exterior boundary. Rays in yellow are observed and rays in green are not observed (yet).}}\label{figurerecurs}
\end{figure}

Interference between $\gamma_1^-$ and $\gamma_2^-$ may prevent one to conclude on the observability of $\gamma_1^-$ or $\gamma_2^-$ from the observability of $\gamma_1^+$ at $\Gamma$ alone. We prove in Corollary \ref{equivsuppmes} that if two outgoing rays from $\d \Omega_2$ are observed, then the two incoming rays are also observed. One then use this result iteratively to obtain the observability of a ray in the general case as depicted in figure \ref{figurerecurs}. Such recursive use of Corollary \ref{equivsuppmes} is stated as Corollary \ref{equivsuppmesgraph}.

If $\Omega$ and $\Omega_2$ are assumed to be strictly convex, $c_2>c_1$ (in order to prevent trapped rays in $\Omega_2$) and $\Gamma = \d \Omega$, then one obtains the observability of \eqref{wave}-\eqref{Dirhom} from Corollary \ref{equivsuppmesgraph}. This result was already noticed by Lebeau, Le Rousseau, Terpolilli and Trelat \cite{Cnam}. When $\Gamma \subsetneq \d \Omega$, a geometrical argument is needed to ensure that Corollary \ref{equivsuppmesgraph} may be applied succesfully to observe every rays of \eqref{wave}-\eqref{Dirhom}. The main contribution of this paper is to provide a geometrical construction allowing to give sufficient conditions for the observability of \eqref{wave}-\eqref{Dirhom} when $\Gamma \subsetneq \d \Omega$. Before introducing the geometrical argument and the main results, let us first introduce some notations. 

\subsection{Notations}

We begin by recalling the definition of the $\Gamma(x_0)$ observability region for the wave equation (\cite{LionsBook1}). Let $x_0 \in \RR^2 \setminus \overline{\Omega}$. Then,
\begin{equation}\label{defgammax0}
\Gamma(x_0):=\{ x\in \d \Omega \, | \, \<(x-x_0),n(x)\> >0 \},
\end{equation}
where $\<.,.\>$ denotes the $\RR^2$ inner product. We recall that if $\Omega$ is strictly convex, then $\Gamma(x_0)$ is a connected part of $\d \Omega$. 

Consider the parametrisation $\delta(s), s\in [0,1]$ of $\d \Omega$ in the counter-clockwise direction, normalized such that $\delta'(s)$ is the unit tangent vector at $\delta(s), s\in [0,1]$ in the direction of propagation of $\delta$. For simplicity, we fix $\delta(0)=\delta(1)$ such that $\|x_0-\delta(0)\|=\min\{\left \|x_0-\delta(s)\| \, \right| \, s\in [0,1]\}$. The parametrisation $\delta_2(s), s\in [0,1]$ of $\d \Omega_2$ is chosen with the same requirements. 

We define what is referred to the collision map $\F$ in the billiard literature (\cite{Billiards}) for $\Omega_1$,
\begin{align*}
\F :  \left(\d \Omega \cup \d \Omega_2 \right)\times \RR^2 & \longrightarrow \left(\d \Omega \cup \d \Omega_2 \right)\times \RR^2, \\
 (x,\xi)\qquad  & \longmapsto \qquad (x^1,\xi^1),
\end{align*}
where $x^1$ is the point where the ray of $\Omega_1$ starting from $x$ travelling in the $\xi$ direction at constant speed and in straight line intersects $\d \Omega \cup \d \Omega_2$ and $\xi^1 \in \RR^2$ is the direction of the outgoing ray reflected according to the law of reflection of \eqref{ClassicWave}. Iterations of the collision map are denoted $\F^n, n\in \ZZ$. For $(x,\xi)\in (\d \Omega \cup \d \Omega_2) \times \RR^2$ we denote the projections $\Pi_x(x,\xi)=x$ and $\Pi_{\xi}(x,\xi)=\xi$. Notice that not all directions $\xi\in \RR^2$ are admissible directions for the domain of definition of $\F$ since half those directions correspond to incoming ray but it is customary to identify these directions to their unique outgoing direction \cite{Billiards}. We define in the same way $\F_2$ the collision map of $\Omega_2$ 
\begin{align}
\F_2 : \left(\d \Omega \cup \d \Omega_2 \right)\times \RR^2 & \longrightarrow  \left(\d \Omega \cup \d \Omega_2 \right) \times \RR^2, \\
 (x,\xi)\quad  & \longmapsto \, \, \,  (x^1,\xi^1),
\end{align}
At the boundary $\d \Omega_2$, we identify $\F(x,\xi)$ to $\F_2(x,\xi)$ according to \eqref{SnellDescartes} in the following sense. Suppose a ray starts from $x$ in the direction $\xi$ intersects $\d \Omega_2$ at $\Pi_x(\F(x,\xi))$ with an angle 
\[
\theta_1=\arccos\left(\frac{\<n_2(\Pi_x(\F(x,\xi))),\Pi_{\xi}(\F(x,\xi))\>}{\|\Pi_{\xi}(\F(x,\xi))\|}\right).
\] 
If a ray is transmitted with an angle $\theta_2$ according to \eqref{SnellDescartes}, then 
\[
\theta_2=\arccos\left(\frac{\<-n_2(\Pi_x(\F(x,\xi))),\Pi_{\xi}(\F_2(x,\xi))\>}{\|\Pi_{\xi}(\F_2(x,\xi))\|}\right).
\]

\subsection{Strategy of the proof and statement of the main results}

We introduce formally the geometrical construction used in this paper to provide sufficient conditions for the observability \eqref{wave}-\eqref{Dirhom}. The geometrical construction is obtained iteratively.  \\

\begin{figure}[!ht]
\begin{center}
	\includegraphics[height=5cm]{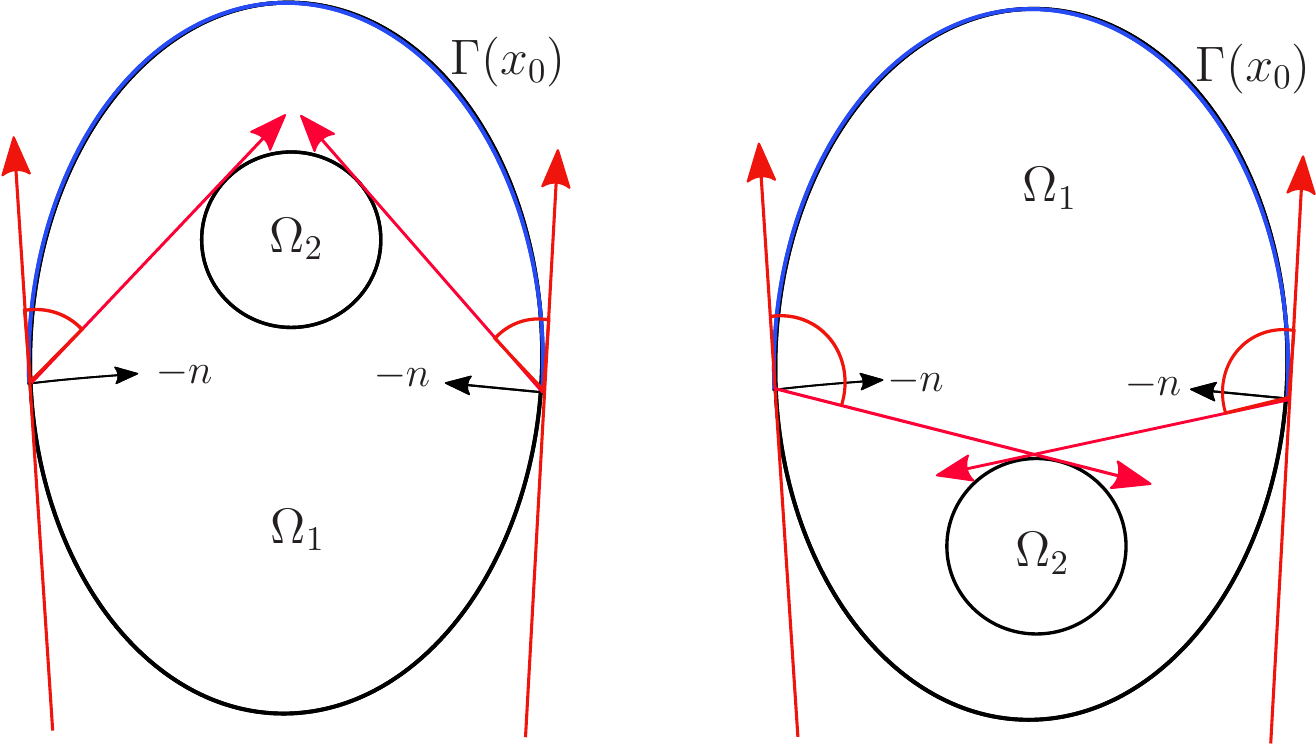} 
\end{center}
\caption{\footnotesize{Left : case where $\Omega_2$ is included in $\Gamma(x_0)$. Right : case where $\Omega_2$ is not included in $\Gamma(x_0)$}}\label{Omega2inclus}
\end{figure}

\begin{figure}[!ht]
\begin{center}
	\includegraphics[height=5cm]{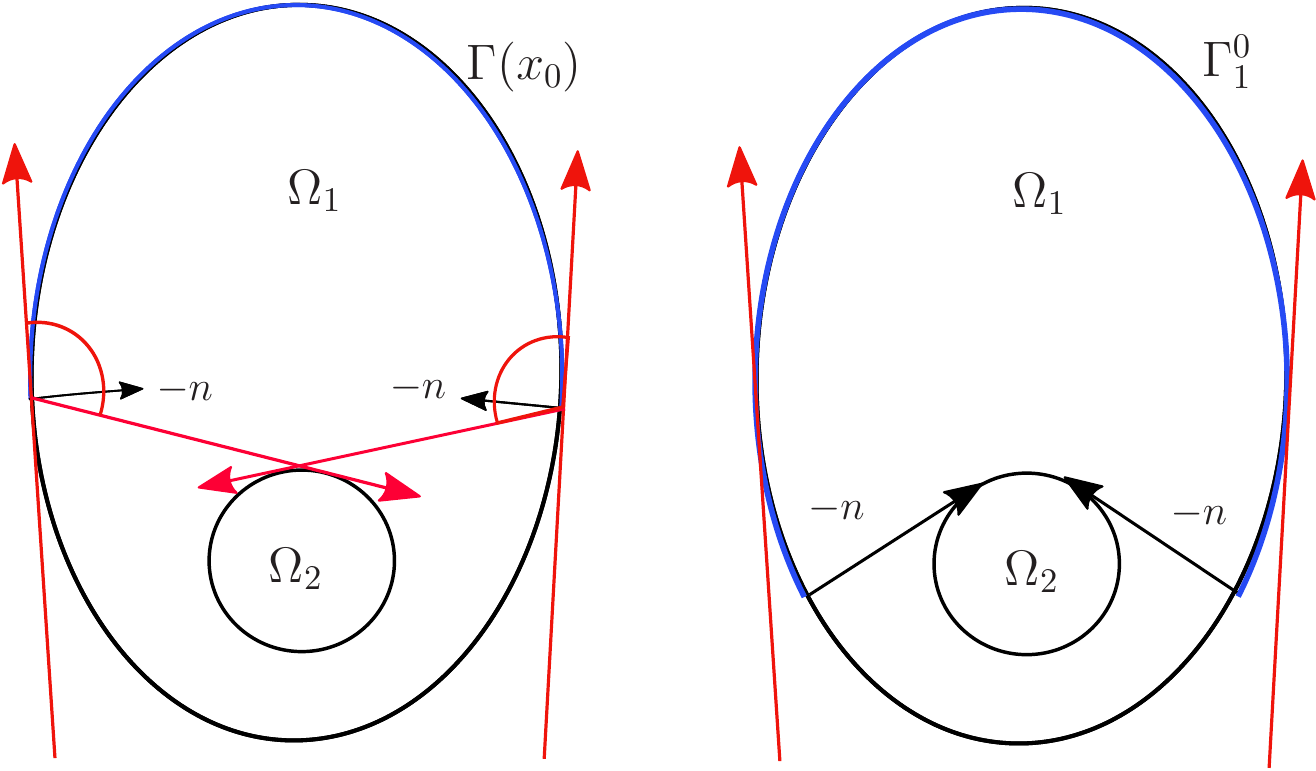} 
\end{center}
\caption{\footnotesize{Left : case where $\Omega_2$ is not included in $\Gamma(x_0)$. Right : extension of $\Gamma(x_0)$ to $\Gamma_1^0$}}\label{PropagG}
\end{figure}

\textbf{Step 1 : Localisation of $\Omega_2$ with respect to $\Gamma(x_0)$} \\
If $\Omega_2$ is not included in $\Gamma(x_0)$ (see figure \ref{Omega2inclus}, on the right), then one extend the boundary of $\Gamma(x_0)$ until the ray starting from the boundary of $\Gamma(x_0)$ and travelling in the $-n$ direction intersects (tangentially) $\overline{\Omega}_2$ (see figure \ref{PropagG}). We denote $\Gamma_1^0$ this extension of $\Gamma(x_0)$. If $\Omega_2$ is included in $\Gamma(x_0)$ (see figure \ref{Omega2inclus}, on the left), one defines $\Gamma_1^0:=\Gamma(x_0)$. 

\begin{figure}[!ht]
\begin{center}
	\includegraphics[height=5cm]{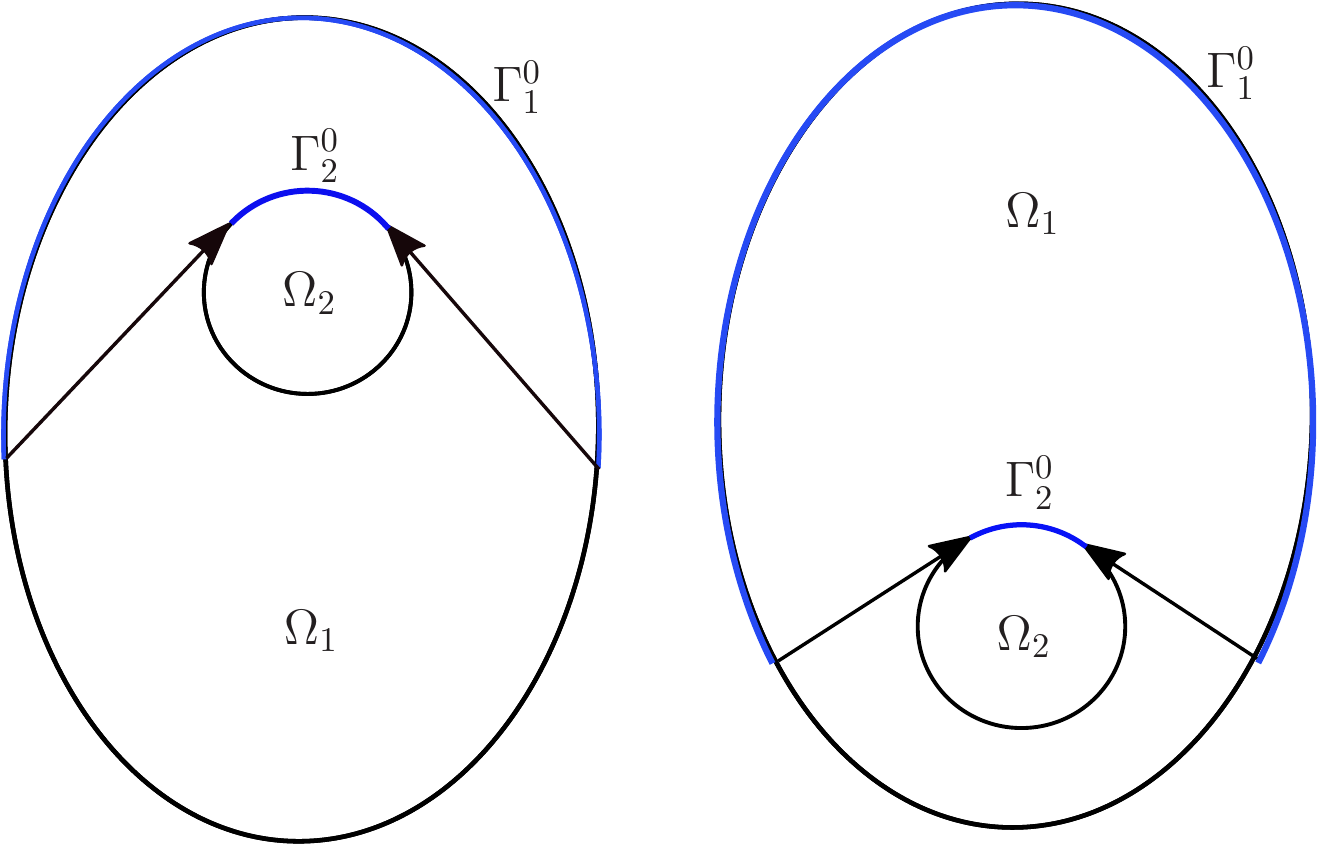} 
\end{center}
\caption{\footnotesize{Creation of $\Gamma_2$. Left : case where $\Omega_2$ is included in $\Gamma(x_0)$. Right : case where $\Omega_2$ is not included in $\Gamma(x_0)$}}\label{CreationG2}
\end{figure}

\textbf{Step 2 : Creation of $\Gamma_2^0$}\\
We create the boundary of $\Gamma_2^0 \subset \d \Omega_2$ with the tangential contact obtained during the localisation of $\Omega_2$ (see figure \ref{CreationG2}, on the left) or with the tangential contact of rays starting from the boundary of $\Gamma_1^0$ and travelling in the $-n$ direction if $\Gamma(x_0)$ was extended (see figure \ref{CreationG2}, on the right).\newline

\textbf{Step 3 : Extension of $\Gamma_2^0$ to $\Gamma_2^1$}\\

\begin{figure}[!ht]
\begin{center}
	\includegraphics[height=5cm]{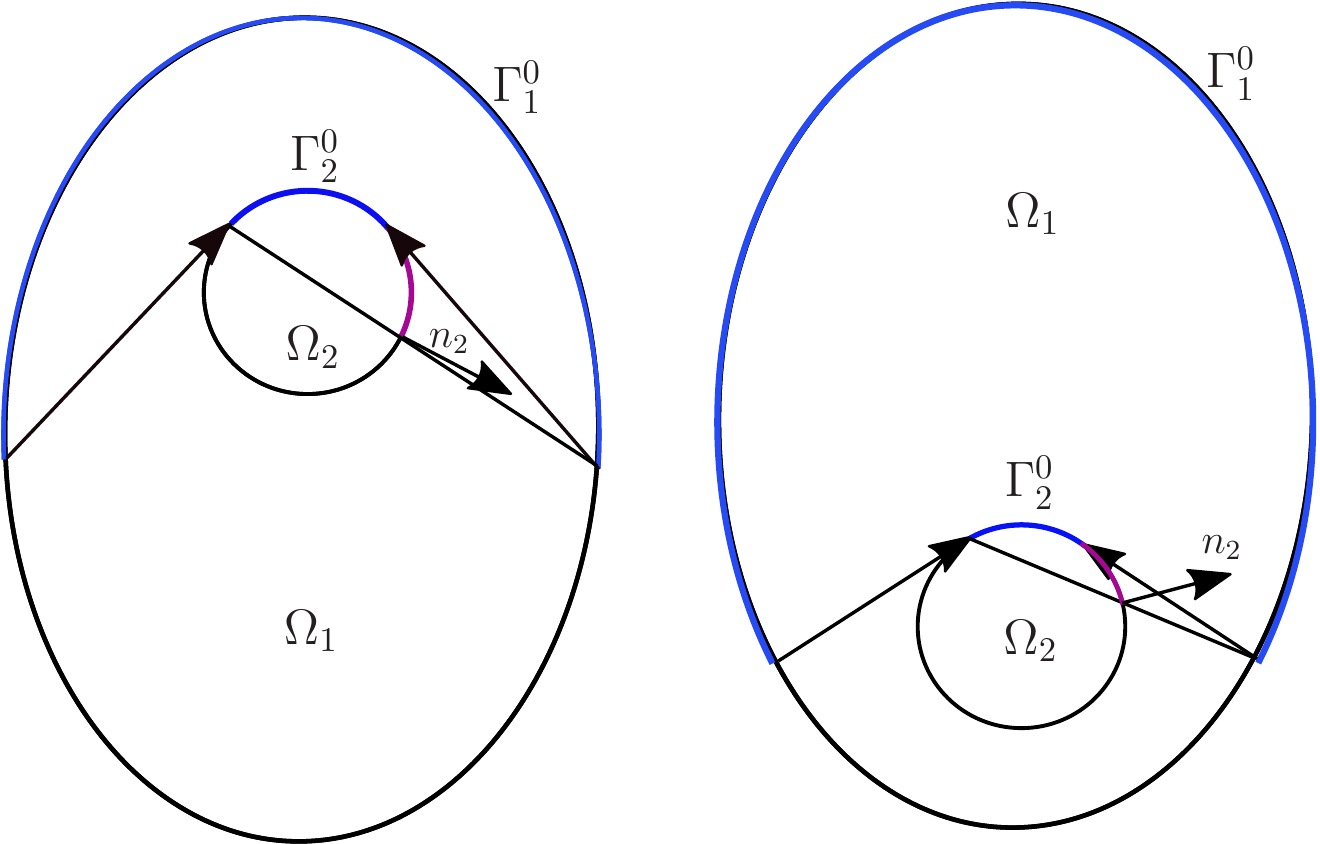} 
\end{center}
\caption{\footnotesize{Extension of one of the boundary of $\Gamma_2^0$ in two different cases.}}\label{ExtensionG2}
\end{figure}

We extend $\Gamma_2^0$ to $\Gamma_2^1$ in the following way. Draw the line between one point of $\d \Gamma_2^0$ and the opposite point of the boundary $\d \Gamma_1^0$ (see figure \ref{ExtensionG2}). If this line has no intersection with $\d \Omega_2$, then there is no extension of $\Gamma_2^0$ from this boundary. Otherwise, $\Gamma_2^1 \supset \Gamma_2^0$ is bounded by these intersections and is such that for rays starting from a point in $\Gamma_2^1$ and travelling in the $n_2$ direction intersect $\Gamma_1^0$ (see figure \ref{ExtensionG2}).

\begin{figure}[!ht]
\begin{center}
	\includegraphics[height=5cm]{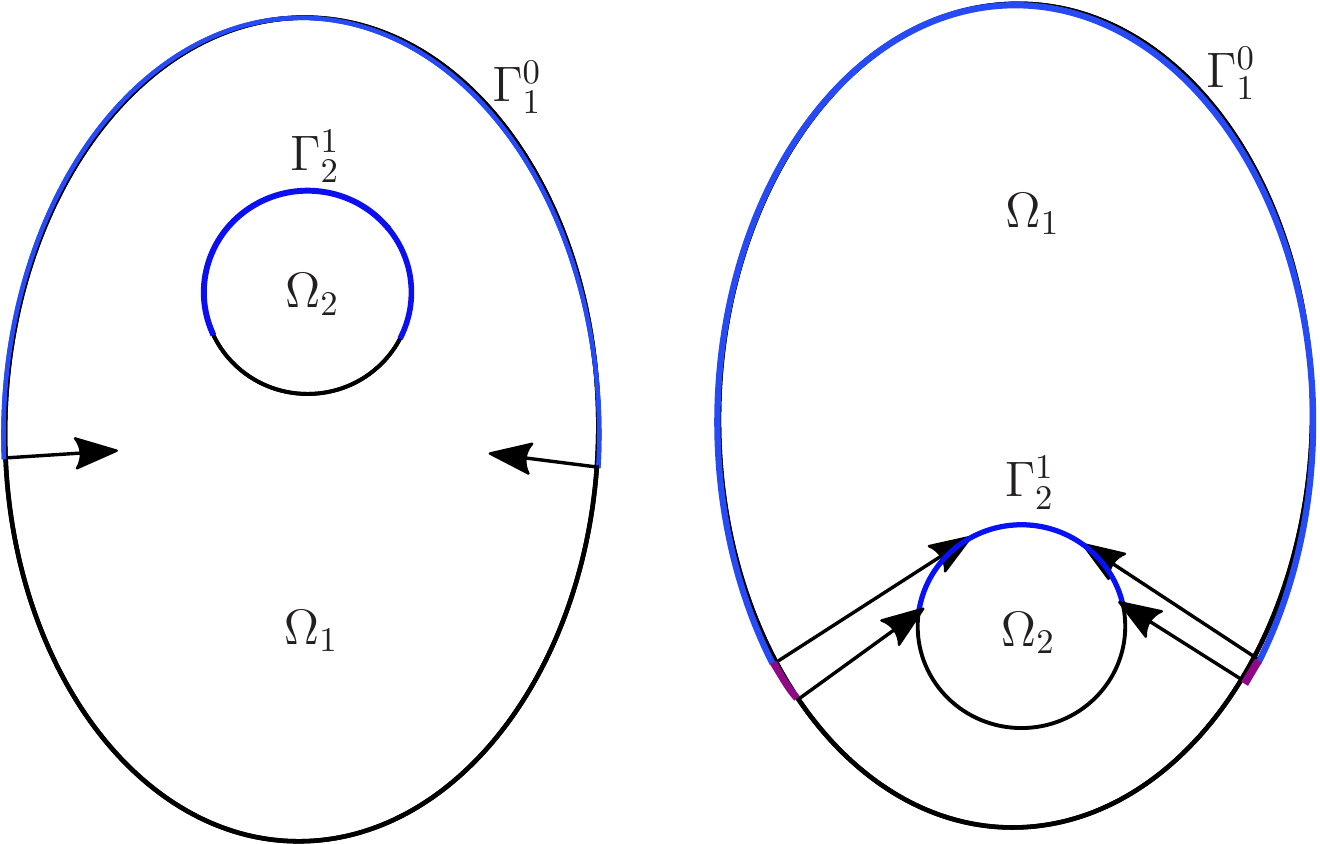} 
\end{center}
\caption{\footnotesize{Extension of one of the boundary of $\Gamma_1^0$. Left : $\Gamma_1^0$ is not extended. Right : $\Gamma_1^0$ is extended.}}\label{ExtensionG1}
\end{figure}

\textbf{Step 4 : Extension of $\Gamma_1^0$ to $\Gamma_1^1$}\\
We extend the boundary of $\Gamma_1^0$ until rays starting from the boundary of $\Gamma_1^0$ in the $-n$ direction intersect the boundary of $\Gamma_2^1$ (see figure \ref{ExtensionG1}) 

\textbf{Step 5 : Iteration of step 3 and 4}

We iterate step 3 and 4 to extend $\Gamma_i^{n-1}$ to $\Gamma_i^{n}, i=1,2, n\in \NN$ until $\Gamma_2^{n-1}=\Gamma_2^n$. We prove in Lemma \ref{extensiongamma12} it happens in a finite number of iterations of step 3 and 4. More precisely, it happens either because there are no longer intersection of the lines defined in step 3 with the boundary of $\Omega_2$, or because the rays starting from the extended region of $\Gamma_2^n$ and propagating in the $n_2$ direction do not intersect $\Gamma_1^n$. Once the iteration process is over, we denote $\Gamma_i:=\Gamma_i^n$. 

In each step of the iteration process, we prove that the extended region is observable, meaning that every rays starting from this region are observed. Therefore, once the extension is proved to be an observable region, one concludes on the observation of a ray when it intersects this region. We are able to prove that the iteration process allows to obtain an observable region on $\Omega_2$ that satisfies GCC.
\begin{lem}
Let $\Gamma_2$ be defined by the iteration process. Then there exists $x_0^2\in \RR^2 \setminus \overline{\Omega}_2$ such that $\Gamma_2=\Gamma(x_0^2)$.  
\end{lem}

\begin{figure}[!ht]
\begin{center}
	\includegraphics[height=5cm]{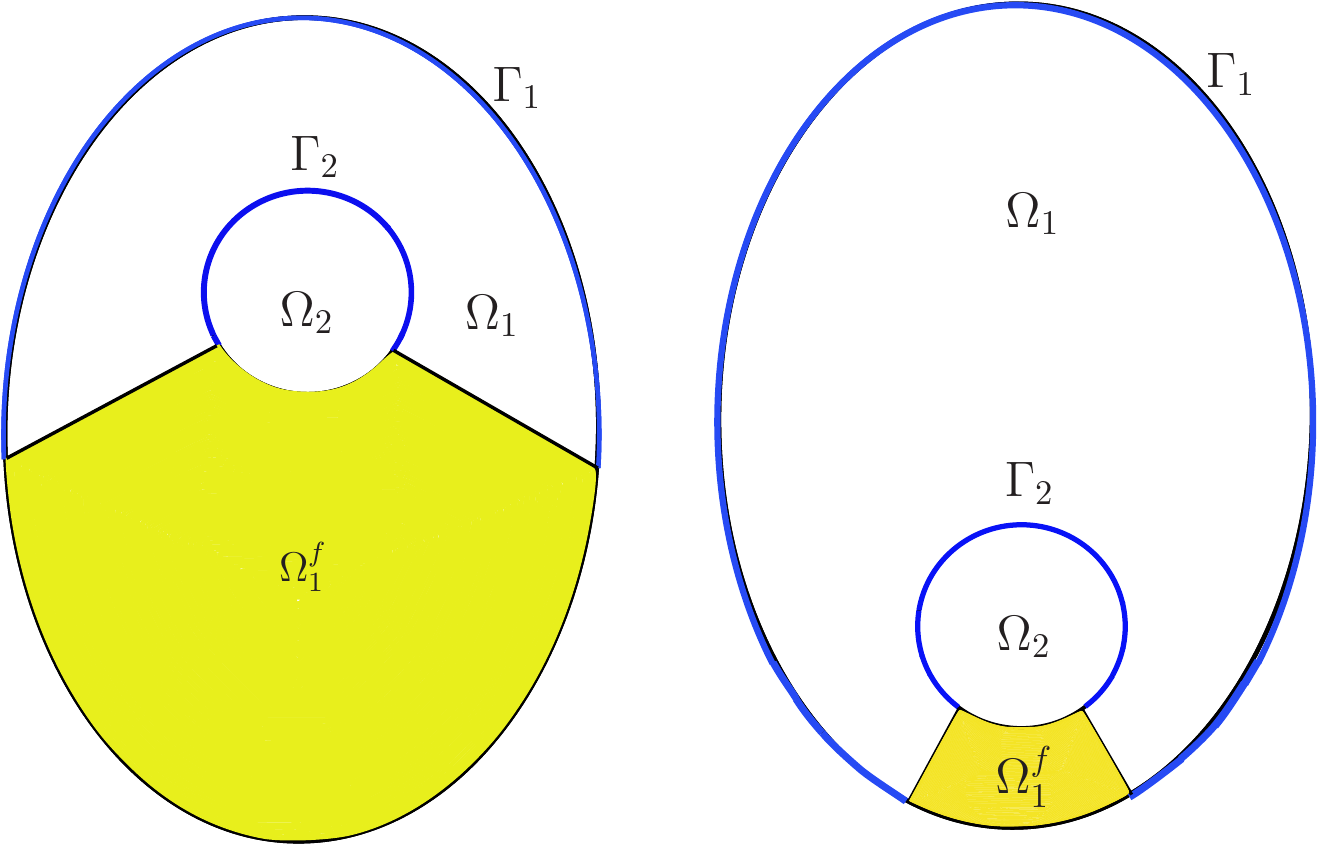} 
\end{center}
\caption{\footnotesize{Definition of $\Omega_1^f$.}}\label{Defof}
\end{figure}

We define $\Omega_1^f\subset \Omega_1$ the remaining part of $\Omega$ bounded by $\d \Omega \setminus \Gamma_1$, $\d \Omega_2 \setminus \Gamma_2$ and the lines $l(\delta(s_i^1),\delta_2(s_i^2))$ passing through $ \delta(s_i^1)\in \d \Gamma_1, \delta_2(s_i^2) \in \d \Gamma_2, i=1,2$ (see figure \ref{Defof}). At this point, rays propagating in $\Omega_1^f$ were not shown to be observable. Their observability relies on the type of geometry $\Omega_1^f$ defines. 

We introduce the notion of uniformly escaping geometry. For convenience, we identify $n(\delta(s))^\perp$ with $\delta'(s)$, that is the tangential vector in the direction of propagation.
\begin{defn}[Uniformly escaping geometry]
We say that $\Omega_1^f$ is a uniformly escaping geometry if the application 
\begin{align}\label{directionmap}
\M:  (\d \Omega_2 \setminus \Gamma_2) \times \RR^2 & \longrightarrow \qquad \qquad \RR \\
(x,\xi) \quad & \longmapsto  \<\xi,n\left(\Pi_x\left(\F\left(x,\xi \right)\right)\right)^\perp\>
\end{align}
is nondecreasing for $s \mapsto \M(\delta_2(s),n_2(\delta_2(s))), \delta_2(s)\in \d \Omega_2 \setminus \Gamma_2$.
\end{defn}

The name escaping geometry refers to the work of Miller in \cite{Escape} on escape functions where GCC is reinterpreted in terms of escaping rays. In \cite{Escape}, a ray is said to have escaped for \eqref{ClassicWave} through the observability region $\Gamma$ if the ray intersects $\Gamma$ in a non-diffractive way. Such a ray is then not reflected on $\Gamma$ and is assumed to have escaped $\Omega$. The geometric control condition is then equivalent to ask every rays to have exited $\Omega$ through $\Gamma$ in time $T>0$. The notion of uniformly escaping geometry is an adaptation of this definition in the framework of \eqref{wave}-\eqref{Dirhom}. 

We consider the boundary $l(\delta(s_i^1),\delta_2(s_i^2)), \delta(s_i^1)\in \d \Gamma_1,  \delta_2(s_i^2) \in \d \Gamma_2, i=1,2$ to be escaping for $\Omega_1^f$. Rays that go through this boundary are not reflected back to $\Omega_1^f$. These rays will propagate to $\Omega_1 \setminus \Omega_1^f$ and intersect eventually $\Gamma_1$ or $\Gamma_2$, regions that we will prove to be observable. We consider that no rays can escape through $\d \Omega \setminus \Gamma_1$ and that only the rays intersecting $\d \Omega_2 \setminus \Gamma_2$ in the $-n_2$ direction can escape through this boundary. Indeed, from the Snell-Descartes laws \eqref{SnellDescartes}, such rays are transmitted to $\Omega_2$ no matter what $c_1$ and $c_2$ are. 

\begin{figure}[!ht]
\begin{center}
	\includegraphics[height=2cm]{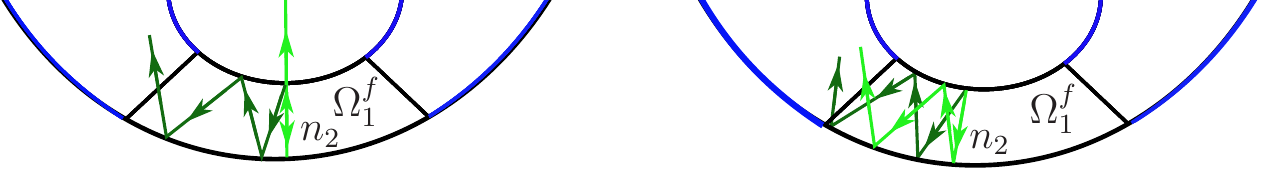} 
\end{center}
\caption{\footnotesize{Example of a uniformly escaping geometry. Left : (light green) ray from $x$ in the $n_2$ direction such that $\M(x,n_2(x))=0$ - (dark green) half-ray from the same point propagating in the negative tangential direction (with respect to $\delta'$). Right : (light green) ray from $x$ in the $n_2$ direction such that $\M(x,n_2(x))<0$ - (dark green) half-ray from the same point propagating in the negative tangential direction (with respect to $\delta'$). }}\label{Uegfig}
\end{figure}

In our context, the uniformly escaping geometry ensures that every rays propagating in $\Omega_1^f$ will be observed, again no matter what $c_1$ and $c_2$ are. Indeed, by definition, a ray starting in the $\M(x,n_2(x))=0, x\in \d \Omega_2 \setminus \Gamma_2$ region and in the direction $n_2(x)$ satisfy $\F^2(x,n_2(x))=(x,n_2(x))$ by definition. Since this is an escaping direction for $\d \Omega_2 \setminus \Gamma_2$, this ray is assumed to have escaped $\Omega_1^f$ (see the light green ray in figure \ref{Uegfig} on the left). A ray starting in the $\M(x,n_2(x))<0$ region in the $n_2(x)$ direction will eventually escape through $l(\delta(s_2^1),\delta_2(s_2^2)), \delta(s_2^1)\in \d \Gamma_1,  \delta_2(s_2^2) \in \d \Gamma_2$ thanks to the nondecreasing assumption on $\M$ (see the light green ray in figure \ref{Uegfig} on the right). The description is symmetric for rays in the $\M(x,n_2(x))>0$ region in the $n_2(x)$ direction. The complete picture can be deduced by this analysis. Indeed, consider a ray in the $M(x,n_2(x))<0$ region propagating in the opposite direction of the parametrisation of $\delta$. It is always possible to follow such a ray since $\xi \neq n_2(x)$ implies that $x$ is a point where a reflection occurs. One can then choose to follow the half-ray propagating in the opposite direction of the parametrisation $\delta$. Since the ray propagating in the $n_2(x)$ direction have escaped in the opposite direction of propagation, then the half-ray propagating in the same direction also escape through the same boundary (see the dark green ray in figure \ref{Uegfig} on the right). The case $M(x,n_2(x))<0$ is symmetric and one can follow either half-ray in the region $M(x,n_2(x))=0$ as both half-ray escapes through one of the two boundary $l(\delta(s_i^1),\delta_2(s_i^2)), \delta(s_i^1)\in \d \Gamma_1,  \delta_2(s_i^2) \in \d \Gamma_2, i=1,2$ (see the dark green ray in figure \ref{Uegfig} on the left for the propagation of one of the half-ray). 

We are finally now able to state the main results of this paper. 

\begin{thm}\label{main}
Let $x_0 \in \RR^2 \setminus \overline{\Omega}$ and $\Gamma(x_0)$ be defined as in \eqref{defgammax0}. If $\Omega_1^f$ is a uniformly escaping geometry, then for every $c_2>c_1$, there exists a time $T>0$ such that \eqref{wave}-\eqref{Dirhom} is observable. 
\end{thm}

Notice that the weaker condition of assuming $\M$ is negative, null and positive on respective parts of $\d \Omega \cap \Omega_1^f$ and $\d \Omega_2 \cap \Omega_1^f$ is not sufficient since it does not prevent change of curvature of the boundary and therefore may create focusing dynamics (see \cite{Billiards}) that wouldn't leave $\Omega_1^f$. If $\Omega_1^f$ does not satisfy the uniformly escaping condition, one has to impose much more restrictive conditions to ensure observability. 

\begin{thm}\label{thmtrap}
Let $x_0 \in \RR^2 \setminus \overline{\Omega}$ and $\Gamma(x_0)$ be defined as in \eqref{defgammax0}. If $\Omega_1^f$ is not a uniformly escaping geometry but if $\Omega_1^f$ and $c_2>c_1$ are such that for every $x\in \d \Omega \setminus \Gamma_1$ and $\xi \in \RR^2$ such that $\Pi_x(\F(x,\xi))\in \d \Omega_2 \setminus \Gamma_2$ and $\Pi_x(\F^2(x,\xi))\in \d \Omega \setminus \Gamma_1$, the transmitted rays in $\Omega_2$ are such that $\Pi_x(\F_2^2(x,\xi^{\pm})) \in \Gamma_2$, then \eqref{wave}-\eqref{Dirhom} is observable in some time $T>0$. 
\end{thm}

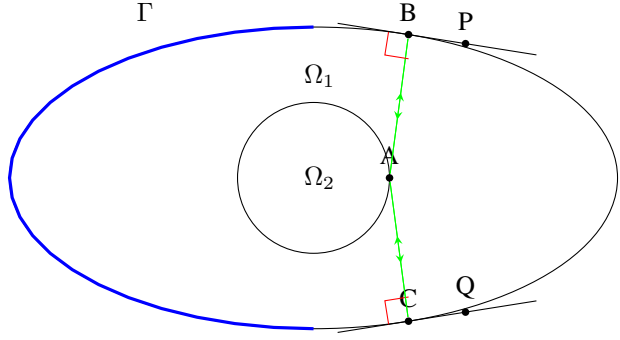
\begin{figure}[!h]
\begin{center}
\begin{tikzpicture}[dot/.style={circle,inner sep=1pt,fill,label={#1},name=#1},
  extended line/.style={shorten >=-#1,shorten <=-#1},
  extended line/.default=1cm]
 \draw (0,0) ellipse (4cm and 2cm);
\node [dot=A] at (1,0) {};
\node [dot=B] at (1.25,1.9) {};
\node [dot=C] at (1.25,-1.9){};
\node [dot=P] at (2,1.78) {};
\node [dot=Q] at (2,-1.78) {};
\draw (0,0) circle (1cm);
\node[text width=.25cm] at (0,0) {$\Omega_2$};
\node[text width=.25cm] at (0,1.35) {$\Omega_1$};
\node[text width=.25cm] at (-2.2,2.2) {$\Gamma$};
\draw[-,green,postaction={mid arrow=green}](A) -- (B) ;
\begin{scope}[shift={(B)},rotate=180]
\draw[-,green,postaction={mid arrow=green}](B) -- (A) ;
\end{scope}
\draw [red,right angle symbol={P}{B}{A}];
\draw[extended line] (B) -- (P) ;
\draw[extended line] (C) -- (Q) ;
\draw [red,right angle symbol={Q}{C}{A}];
\draw[-,green,postaction={mid arrow=green}](A) -- (C)  ;
\begin{scope}[shift={(C)},rotate=180]
\draw[-,green,postaction={mid arrow=green}] (C) -- (A)  ;
\end{scope}
\draw [blue,very thick,domain=90:270] plot ({4*cos(\x)}, {2*sin(\x)});
\end{tikzpicture}
\end{center}
\caption{Representation of the spatial domain for \eqref{wave}-\eqref{Dirhom} with trapped rays not encountering $\Gamma \subset \d \Omega_1$ satisfying GCC. Here ${\Omega_2=\{x\in \RR^2 \, | \, |x| < 1\}}$, $\Omega_1=\{(x_1,x_2)\in \RR^2 \, | |x_1/4|^2+|x_2/2|^2<1 \textrm{ and } x\notin \Omega_2 \}$ and $c_1=1$, $c_2=\sqrt{2}$.}\label{figuretrapped}
\end{figure}

Rays propagating in $\Omega_1^f$ if $\Omega_1^f$ is not assumed to be a uniformly escaping geometry are difficult to analyse. One can exhibit examples where $\Omega_2$ is an obstacle for certain direction of rays incoming from $\Omega_1$ and create trapped rays for \eqref{wave}-\eqref{Dirhom} even if $c_2>c_1$ (see figure \ref{figuretrapped}). Therefore, one has to ensure (by changing the coefficient $c_2>c_1$ or by taking a larger observability region $\Gamma$) that such trapped rays does not exist for $\Omega_1^f$. 

But even the hypothesis that there are no trapped rays in $\Omega_1^f$ is not sufficient to conclude on the observability of \eqref{wave}-\eqref{Dirhom} as the transmitted rays in $\Omega_2$ could intersect $\d \Omega_2 \setminus \Gamma_2$ and be transmitted back to $\Omega_1^f$. One has to keep track of these rays and to make sure that they are eventually observed by $\Gamma$ (in the sense of Corollary \ref{equivsuppmesgraph}). More generally rays in $\Omega_1^f$ consist of rays that are bouncing back and forth between $\d \Omega_2 \cup \overline{\Omega_1^f}$ and $\d \Omega \cup \overline{\Omega_1^f}$ that are not of finite cycle and that are never escaping $\Omega_1^f$. For these rays, one has to make sure that the transmitted to $\Omega_2$ provides the observation. But recall that such rays may very well intersect $\d \Omega_2 \setminus \Gamma_2$ and create transmitted rays back in $\Omega_1^f$. Therefore a fine analysis of the collision map is required in this case in order to state sharp sufficient conditions for non-uniformly escaping geometry $\Omega_1^f$. This analysis is beyond the scope of this paper and sufficient conditions in this case are presented in Theorem \ref{thmtrap}. 

The $\Gamma(x_0)$ condition may seem unexpected for Theorem \ref{main} and Theorem \ref{thmtrap} in the context of a proof relying on microlocal arguments. It was already proved in the paper of Bardos, Lebeau and Rauch (\cite{BLR}) that there exists strictly convex $\Omega$ and non-connected $\Gamma$ that satisfies GCC. Since $\Gamma(x_0)$ is connected for strictly convex $\Omega$, then $\nexists x_0 \in \RR^2$ such that $\Gamma=\Gamma(x_0)$ for such $\Gamma$. It was proved later on by Miller (\cite{Escape}) that $\Gamma(x_0)$ always imply GCC using escape functions. 

We prove here that GCC is more general than the $\Gamma(x_0)$ condition even when $\Omega$ is strictly convex and $\Gamma$ is connected.

\begin{lem}\label{GCCmoregen} 
There exist $\Omega \subset \RR^2$ open, strictly convex domains and $\Gamma \subset \d \Omega$ connected that satisfy GCC such that $\nexists x_0 \in \RR^2$ such that $\Gamma=\Gamma(x_0)$.
\end{lem}

The proof relies on the specific dynamic of the billiards on the ellipse (billiards have the same rules than the optic geometry; see \cite{Billiards} for a good introduction on billiards). More precisely, we use the 1-dimensional folliation of its phase space \cite{Billiards}. The (strong) Birkhoff's conjecture states that this 1-dimensional folliation, or integrability, is particuliar to the ellipse.
\begin{conj}[Birkhoff]
The only integrable billiards in two dimension are the ellipses.
\end{conj}

\begin{figure}[!ht]
\begin{center}
	\includegraphics[height=5cm]{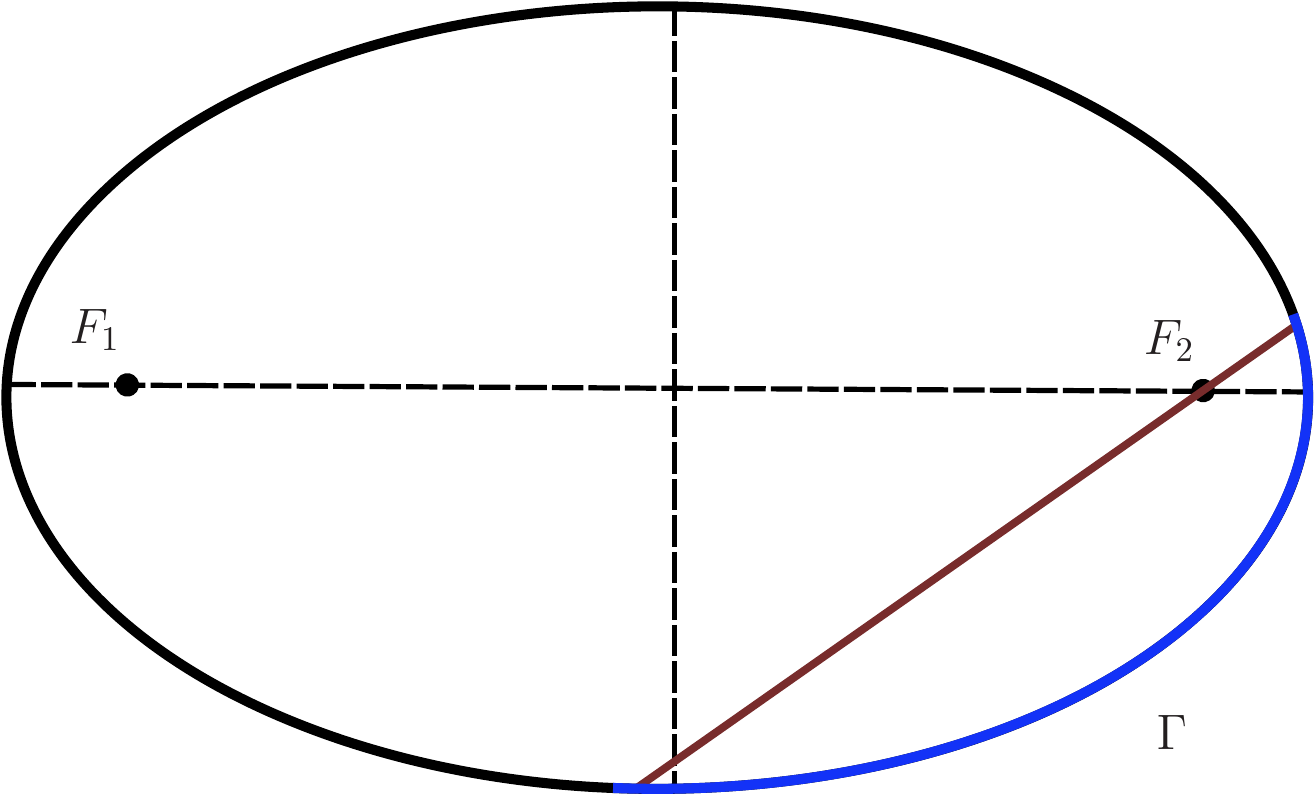} 
\end{center}
\caption{\footnotesize{Example of an observable region $\Gamma \subset \d \Omega$ for \eqref{ClassicWave}.}}\label{ellips}
\end{figure}

Figure \ref{ellips} illustrates an observable region $\Gamma \subset \d \Omega$ given by Lemma \ref{GCCmoregen}. 

From Lemma \ref{GCCmoregen} we deduce the rather surprinsing corollary. 

\begin{cor}
There exists $\Omega \subset \RR^2$ open, strictly convex domains and $\Gamma \subset \d \Omega$ such that $\Gamma$ and $\d \Omega \setminus \Gamma$ satisfy GCC. 
\end{cor} 
 
Therefore, hypothesis $\Gamma=\Gamma(x_0)$ is not necessary condition, even in the case where $\Omega_1^f$ is a uniformly escaping region. The $\Gamma(x_0)$ condition however provides an understanding on the propagation of rays starting from $x\in \Gamma \setminus \Gamma(x_0)$ and moving in the $-n(x)$ direction. Indeed, the propagation of any such rays for any strictly convex $\Omega$ can be approximated by replacing $\d \Omega \setminus \Gamma(x_0)$ by the cone defined by $l(x_0,x_i)$, $x_i\in \d \Gamma(x_0), i=1,2$. The rays propagating in the cone provide a lower bound on the displacement toward $\Gamma(x_0)$ for any rays propagating in the $-n$ direction starting from $\d \Omega \setminus \Gamma(x_0)$. This property is fundamental in the present analysis and will be referred as the conical assumption. 

\subsection{State of the art}

Only few results are available in the literature regarding the observability problem of \eqref{wave}-\eqref{Dirhom}. It was first raised by Lions in \cite{LionsBook1} where the necessary condition $c_2>c_1$ was stated. Later on, Miller obtained in \cite{Miller} the propagation of the defect measures of \eqref{wave}-\eqref{Dirhom}. The controllability was then stated in a conference by Lebeau, Le Rousseau, Terpolilli and Trelat when $c_2>c_1$ and $\Gamma=\d \Omega$ \cite{Cnam}. We also mention the work of Baudouin, Mercado and Osses \cite{Baudouin} on global Carleman estimates for \eqref{wave}-\eqref{Dirhom}.

The exponential decays in presence of memory and damping of \eqref{wave}-\eqref{Dirhom} was addressed by Calvacanti, Coelho and Domingos Cavalcanti in \cite{Cavalcanti}. One can also found in the same paper a good review of the literature on other transmission problems. 

We cite the work of Burq and Lebeau on system and its application to exponential decay of the energy of solutions to the Lam\'e system \cite{LameBG} (see also \cite{LZ}). The exact controllability of the Lam\'e system was obtained by Belhassem and Raymond \cite{LameBJ}).

\subsection{Structure of the paper}

The paper is organized as follow. In Section \ref{SecDefBichar}, we describe the generalized bicharacteristics for \eqref{wave}-\eqref{Dirhom} following closely the presentation in \cite[Section 2]{LameBJ}. The propagation of the defect measure is then presented in Section \ref{SecDefMeasure}. It relies on the previous work of Burq and Lebeau on the Lam\'e system done in \cite{LameBG} and is presented as in \cite{LameBJ}. The microlocal analysis at the interface is very similar and we therefore only sketch the proofs. Section \ref{Secmain} is devoted to the proof of Theorem \ref{main} and \ref{thmtrap}. Lemma \ref{GCCmoregen} is proved in Section \ref{SecGCCmoregen}. 

\textbf{Acknowledgements}

The author is grateful to Nicolas Burq, Belhassen Dehman and Gilles Lebeau for their help and precious advices.

\section{The generalized bicharacteristics and the propagation theorem}

\subsection{Definition of the generalized bicharacteristics}\label{SecDefBichar}

Let $M_i=\RR_t \times \Omega_i, i=1,2$. We consider $T^*(M_i)=\{ (t,x,\tau_i,\xi_i) \in \RR^6 \, | \, (t,x)\in \RR_t \times \Omega_i \}$. The principal symbol of $P_i=\d_t^2-c_i^2 \Delta$ is given by 
\[
p_i(t,x;\tau_i,\xi_i)=\tau_i^2 - c_i^2 |\xi_i|^2.
\]

\textbf{Propagation of the bicharacteristics inside $\Omega_i$}\newline

The propagation of the bicharacteristics $\gamma_i$ in $\textrm{Char }P_i:=\{(t,x,\tau_i,\xi_i)\in T^*(M_i) \, | \, p_i(t,x;\tau_i,\xi_i)=0\}$ is obtained as the solution of 
\begin{equation}\label{bichar}
\dfrac{d}{d\sigma}(t(\sigma),x(\sigma),\tau_i(\sigma),\xi_i(\sigma))=H_{p_i}(t(\sigma),x(\sigma),\tau_i(\sigma),\xi_i(\sigma)), \quad \sigma \in \RR,
\end{equation}
where $H_{p_i}$ is the Hamiltonian associated to $p_i$. The dual coordinates $(\tau_i,\xi_i)$ solution to \eqref{bichar} are constant with respect to $\sigma$ and the rays - the projection over $\RR_t \times \Omega_i$ of the bicharacteristics - are solution to
\[
\begin{cases}
\dot{t}(\sigma)=-2 \tau_i, & \sigma\in \RR, \\
\dot{x}(\sigma)=2 c_i^2 \xi_i, & \sigma\in \RR.
\end{cases}
\]

\textbf{Description of the bicharacteristics near $\d \Omega$}\newline

For $x\in \d \Omega$, consider the local geodesic coordinates $(x_n,x')\in \RR^2$ such that $x=(x_n,x') \Leftrightarrow (x_n,x')=(0,0)$ and, locally, $\Omega=\{x_n>0\}, \d \Omega = \{ x^n=0\}$ and $x'>0$ is the counter-clockwise direction of propagation of $\delta$. Likewise, we write the dual coordinate $\xi_1=(\xi_1',\xi_1^n)$ where $\xi_1'>0$ is the tangential direction in the counter-clockwise direction and $\xi_1^n>0$ is in the direction of the inward normal $-n(x)$.

The Laplacian takes locally the form 
\[
\Delta=\d_{x_n}^2 + R(x_n,x',D_{x'}). 
\]
where $R(x_n,x',D_{x'}) $ is a second order tangential elliptic operator of real principal symbol $r(x_n,x',\xi'_1)$. Since $\Omega$ is assumed to be strictly convex, every bicharacteristics defined on $\textrm{Char } P_1$ intersect $\d \Omega$ in a \textit{non-diffractive} way. In other words, if we denote $\sigma_0$ the moment when the ray intersects $\d \Omega$, then the bicharacteristic defined on $T^*(\RR_t \times \RR^2)$ and in the direction $\xi_1$ escapes $\Omega$ for $\sigma>\sigma_0$ small. 

We decompose $T^*(\d \Omega)=\H \cup \G^+ \cup \E$ where
\begin{align*}
\H:=&\left\{(t,x,\tau_1,\xi_1)\in T^*(\d \Omega) \, | \, r(x_n,x',\xi_1)>|\tau_1|/c_1\right\}, \\
\G^+:=&\left\{(t,x,\tau_1,\xi_1)\in T^*(\d \Omega) \, | \, r(x_n,x',\xi_1)=|\tau_1|/c_1\right\},\\
\E:=&\left\{(t,x,\tau_1,\xi_1)\in T^*(\d \Omega) \, | \, r(x_n,x',\xi_1)<|\tau_1|/c_1\right\}.
\end{align*}
These sets are the hyperbolic, gliding and elliptic set, respectively, of $T^*(\d \Omega \times (0,T))$. From the strict convexity of $\Omega$, the glancing set $\G^-$ of $\d \Omega$ is reduced to $\G^-=\emptyset$. \newline

\textit{Reflection : $\rho_1 \in \H$.}\newline

If a bicharacteristic intersects $\d \Omega$ transversally, that is $x(\sigma_0) \in \d \Omega$ for some $\sigma_0 \in \RR$ and $c_1|\xi_1'| < |\tau_1|$, then 
\[
p_1(t,x,\tau_1,\xi_1',\xi_1^n)(\sigma_0)=\tau_1^2-c_1^2(|\xi_1^n|^2+r(x_n(\sigma_0),x'(\sigma_0),\xi_1')^2)=0,
\]
admits at most two real roots. Therefore at these points  
\begin{equation}\label{bicharrefl}
\lim_{\sigma\rightarrow \sigma_0^+} \xi_1^+(\sigma)=\xi_1'(\sigma_0)-\dfrac{\sqrt{r}}{c_1}n(x(\sigma_0)), \qquad \lim_{\sigma\rightarrow \sigma_0^-} \xi_1^-(\sigma)=\xi_1'(\sigma_0)+\dfrac{\sqrt{r}}{c_1}n(x(\sigma_0)).
\end{equation}
The bicharacteristic components $(x^{\pm}(\sigma),\xi_1^{\pm})$ obey 
\begin{equation}\label{rays}
\begin{cases}
x^+(\sigma)=x+2c_1^2\sigma \xi_1^+,& 0 < \sigma < \epsilon, \\
x^-(\sigma)=x+2c_1^2\sigma \xi_1^-,& -\epsilon< \sigma< 0. \\
\end{cases}
\end{equation}
For $\sigma>0$ (resp. $\sigma<0$), let $\gamma_1^+(\sigma,\rho_1):=(t(\sigma),x^+(\sigma),\tau_1,\xi_1^+)$ (resp. $\gamma_1^-(\sigma,\rho_1):=(t(\sigma),x^-(\sigma),\tau_1,\xi_1^-)$) the outgoing (resp. incoming) bicharacteristic of $\rho_1$. The generalized bicharacteristic path satisfy $\gamma(0,\rho_1)=\rho_1$ and 
\begin{equation}\label{bicharreflec}
\gamma(0,\rho_1)=\begin{cases} \gamma_1^+(\sigma,\rho_1), & 0 < \sigma < \epsilon, \\ 
 \gamma_1^-(\sigma,\rho_1), & -\epsilon < \sigma <0,
 \end{cases}
\end{equation}
for $\epsilon>0$ small.

\begin{figure}[!ht]
\begin{center}
	\includegraphics[height=5cm]{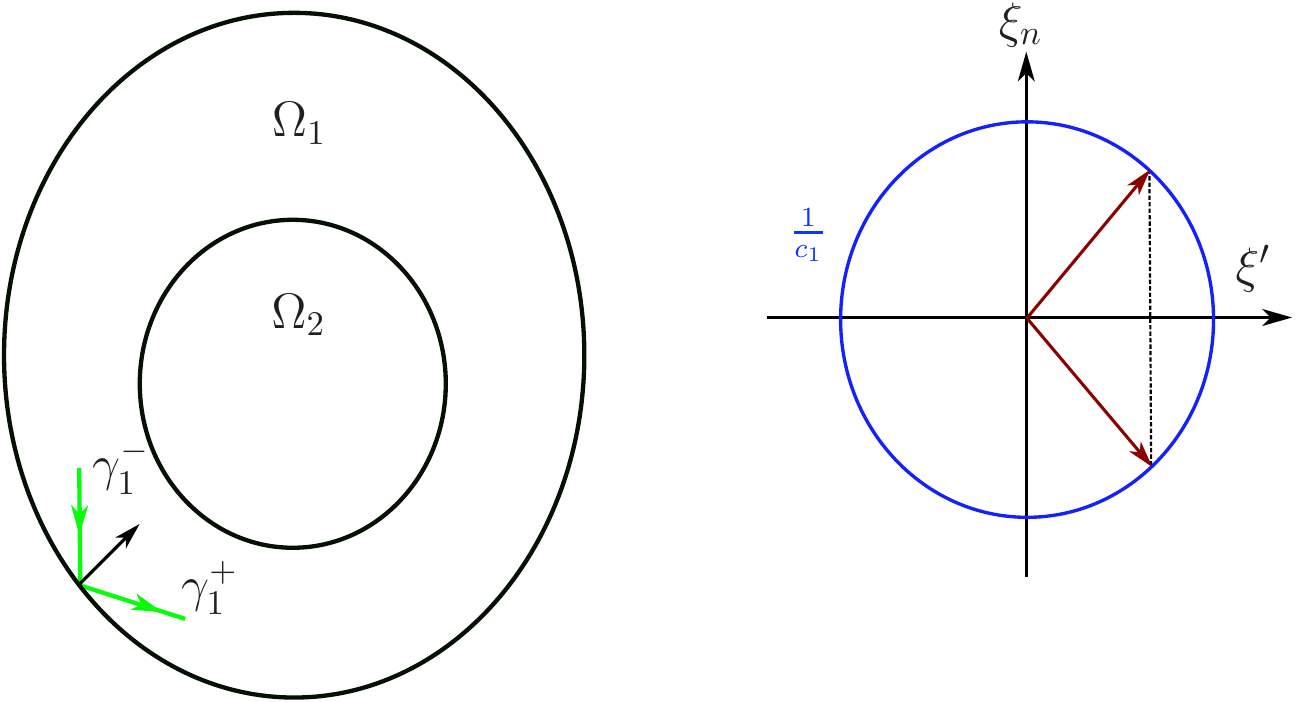} 
\end{center}
\caption{\footnotesize{Rays $(t,x(t))$ of the bicharacteristics near $\rho_1\in \H$}}
\end{figure}

\textit{Gliding rays : $\rho_1 \in \G^+$.}\newline

The gliding ray on $\d \Omega$ is described by the following equation on the bicharcateristic 
\begin{equation}\label{gliding}
\begin{cases}
\widetilde{\gamma_1}(\sigma)=(t(\sigma),x(\sigma),\tau_1(\sigma),\xi_1(\sigma)), \sigma\in [a,b], \textrm{such that }(t(\sigma),x(\sigma),\tau_1(\sigma),\xi_1(\sigma)) \in \G^+ \\
\dot{\tau_1}(\sigma)=0, \dot{t}(\sigma)=-2 \tau_1, \dot{x}(\sigma)=2c_1 \xi'_1(\sigma) \textrm{ and } D\xi_1(\sigma)=0, \\
\textrm{where } D \textrm{ denotes the covariant derivative over } \d \Omega_1.  
\end{cases}
\end{equation}
In other words, $\sigma\in [a,b] \rightarrow x(\sigma)\in \d \Omega$ is a geodesic curve such that $\dot{x}(\sigma)=2c_i \xi'_1(\sigma)$.\newline

\begin{figure}[!ht]
\begin{center}
	\includegraphics[height=5cm]{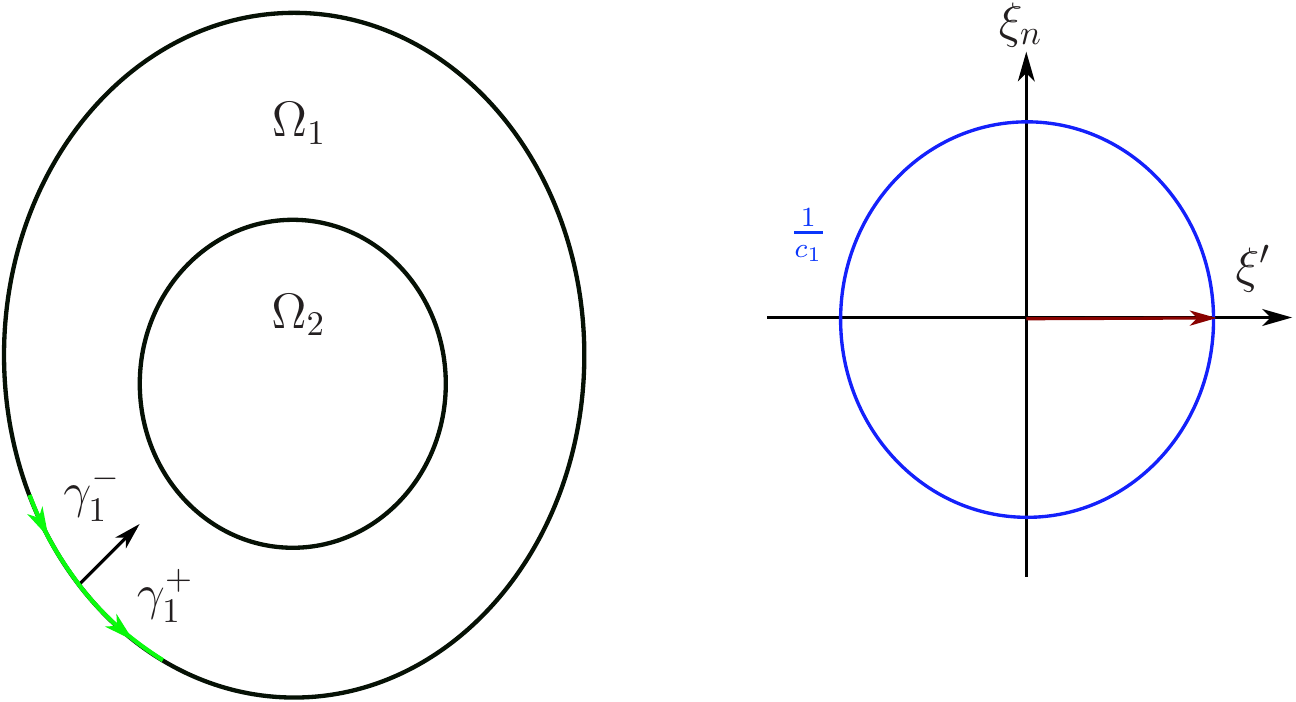} 
\end{center}
\caption{\footnotesize{Rays $(t,x(t))$ of the bicharacteristics near $\rho_1\in \G^+$}}
\end{figure}

\textbf{Description of the bicharacteristics near $\d \Omega_2$.}\newline

We use again the local geodesic coordinates near $\d \Omega_2$. In this case, locally, $\Omega_2=\{x_n >0\}, \d \Omega_2=\{ x_n =0 \}$, $\Omega_1=\{x_n<0\}$ and $x'>0$ is in the direction of propagation of $\delta_2$. For the dual component $\xi_i=(\xi_i',\xi_i^n), i=1,2$, $\xi_1^n>0$ (resp. $\xi_2^n<0$) is in the direction of the inward normal $n_2(x)$ (with respect to $\Omega_1$) (resp. $-n_2(x)$) and $\xi_i'>0$ is in the direction of the propagation of $\delta_2$. The abuse of notation introduced in the introduction shorten the notation $T^*((\d \Omega_1 \cap \d \Omega_2)\times (0,T)) \times T^*(\d \Omega_2 \times (0,T))$ to $T^*(\d \Omega_2 \times (0,T))$. Notice that for $(\rho_1,\rho_2)\in T^*(\d \Omega_2 \times (0,T))$, the transmission boundary condition \eqref{transmission} implies
\begin{equation}\label{transmissionxin}
c_1^2 \xi_1^n = c_2^2 \xi_2^n, 
\end{equation} 
Therefore $T^*(\d \Omega_2 \times (0,T))$ is decomposed as $(\H^1\times \H^2) \cup (\H^1\times \G^{2,+}) \cup (\H^1\times \E^2) \cup (\G^{1,-}\times \E^2) \cup (\E^1\times \E^2)$, where
\begin{align*}
\H^1\times \H^2=&\{\left. (\rho_1,\rho_2)\in T^*(\d \Omega_2 \times (0,T)) \, \right| \, |\xi_i'|<|\tau_i|/c_i, i=1,2  \textrm{ and } c_1^2 \xi_1^n = c_2^2 \xi_2^n \}, \\ 
\H^1\times\G^{2,+}=&\{\left. (\rho_1,\rho_2)\in T^*(\d \Omega_2 \times (0,T)) \, \right| \, |\xi_1'|<|\tau_1|/c_1, |\xi_2'|=|\tau_2|/c_2 \}, \\
\H^1\times\E^2=&\{\left. (\rho_1,\rho_2)\in T^*(\d \Omega_2 \times (0,T)) \, \right| \, |\xi_1'|<|\tau_1|/c_1, |\xi_2'|>|\tau_2|/c_2 \}, \\
 \G^{1,-}\times\E^2=&\{\left. (\rho_1,\rho_2)\in T^*(\d \Omega_2 \times (0,T)) \, \right| \, |\xi_1'|=|\tau_1|/c_1, |\xi_2'|>|\tau_2|/c_2 \}, \\
  \E^1\times\E^2=&\{\left. (\rho_1,\rho_2)\in T^*(\d \Omega_2 \times (0,T)) \, \right| \, |\xi_i'|>|\tau_i|/c_i, i=1,2 \}. 
\end{align*}
The strict convexity of $\Omega_2$ implies $\G^{2,-}=\emptyset$ and $\G^{1,+}=\emptyset$ on $\d \Omega_2$.\newline 

\textit{Reflection - Transmission : $(\rho_1,\rho_2) \in \H^1\times \H^2$.}\newline

If a bicharacteristic, incoming from $\Omega_i$, intersects $\d \Omega_2$ transversally, that is $x(\sigma_0) \in \d \Omega_2$ for some $\sigma_0 \in \RR$ and $c_i|\xi_i'| < |\tau_i|$, then 
\begin{equation}\label{symbol}
p_i(t,x,\tau_i,\xi_i',\xi_i^n)(\sigma_0)=\tau_i^2-c_i^2(|\xi_i^n|^2+r(x_n(\sigma_0),x'(\sigma_0),\xi_i')^2)=0,
\end{equation}
admits at most two real roots. Therefore at these points  
\begin{equation}\label{bicharreflfrontiere2}
\lim_{\sigma\rightarrow \sigma_0^+} \xi_i^+(\sigma)=\xi_i'(\sigma_0)-\dfrac{\sqrt{r}}{c_i}n(x(\sigma_0)), \qquad \lim_{\sigma\rightarrow \sigma_0^-} \xi_i^-(\sigma)=\xi_i'(\sigma_0)+\dfrac{\sqrt{r}}{c_i}n(x(\sigma_0)).
\end{equation}
Moreover, $\xi_i^\pm$, \eqref{transmissionxin} and \eqref{symbol} defines $\xi_j^\pm, j=3-i$. 

\begin{figure}[!ht]
\begin{center}
	\includegraphics[height=5cm]{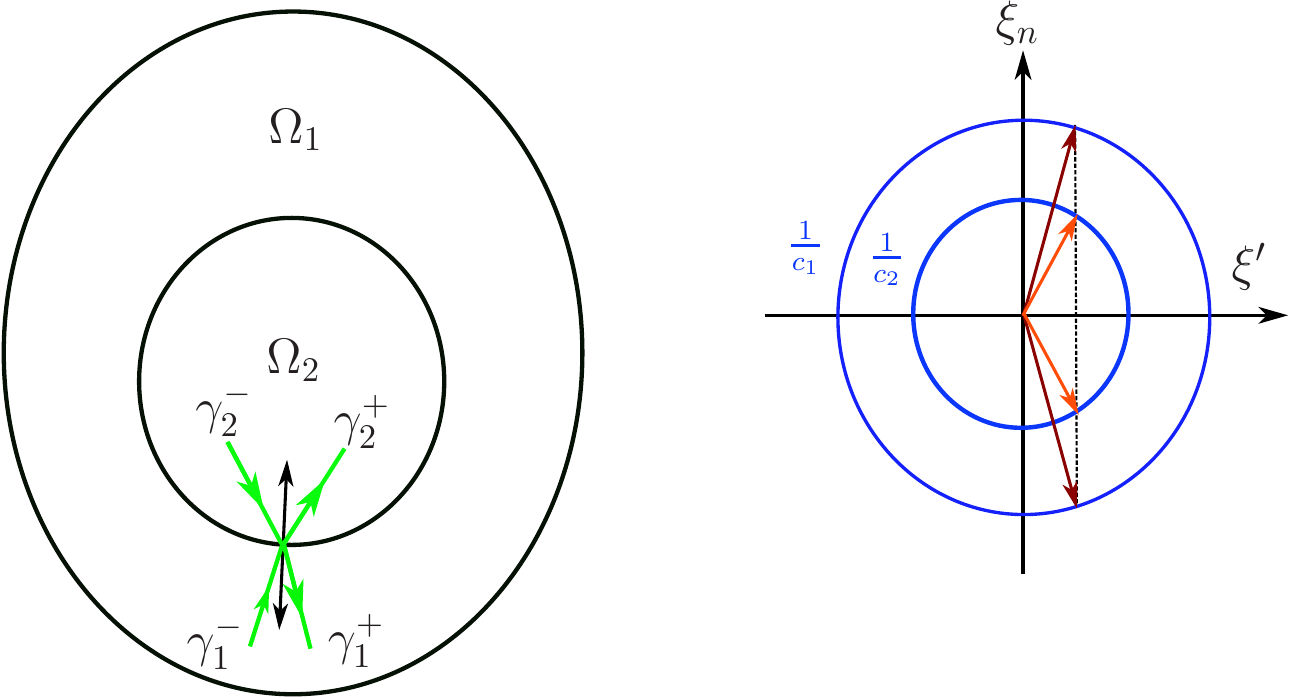} 
\end{center}
\caption{\footnotesize{Rays $(t,x(t))$ of the bicharacteristics near a point $(\rho_1,\rho_2)\in \H^1\times \H^2$}}\label{bicharHH}
\end{figure}

For $\sigma>0$ let $\gamma_i^+(\sigma,\rho_i)=(t(\sigma),x_i^+(\sigma),\tau_i(\sigma),\xi_i^+(\sigma))$ be the reflected bicharacteristic, let $\gamma_j^+(\sigma,\rho_j)=(t(\sigma),x_j^+(\sigma),\tau_j(\sigma),\xi_j^+(\sigma))$ be the transmitted bicharacteristic and let $\gamma_i^-(\sigma,\rho_i)=(t(\sigma),x_i^-(\sigma),\tau_i(\sigma),\xi_i^-(\sigma))$ be the incoming bicharacteristic of $\rho_i$ ($x_{i,j}^{\pm}$ being defined by \eqref{rays}). The generalized bicharacteristic path is such that $\gamma(0,\rho_i)=\rho_i$ and 
\[
\gamma(\sigma,\rho_i)=\begin{cases} \gamma_i^+(\sigma,\rho_i), & 0 < \sigma < \epsilon, \\ 
\gamma_j^+(\sigma,\rho_j), & 0 < \sigma < \epsilon, \\ 
 \gamma_i^-(\sigma,\rho_i), & -\epsilon < \sigma <0.
 \end{cases}
\]
Figure \ref{bicharHH} illustrates the outgoing bicharacteristics from an incoming bicharacteristic $\gamma_1^-$ or $\gamma_2^-$ at a point $(\rho_1,\rho_2)\in \H^1 \times \H^2$. 

Notice that when $(\rho_1,\rho_2) \in \H^1 \times ( \H^2\cup \G^{2,+})$, the assumption $c_1<c_2$ and the Snell-Descartes condition \eqref{SnellDescartes} implies 
\begin{equation}\label{angledescartes}
\theta_1<\theta_2, \quad \forall (\rho_1,\rho_2) \in \H^1\times (\H^2\cup \G^{2,+}).
\end{equation}\newline

\textit{Reflection - Gliding transmission : $(\rho_1,\rho_2) \in \H^1\times \G^{2,+}$.}\newline

\begin{figure}[!ht]
\begin{center}
	\includegraphics[height=5cm]{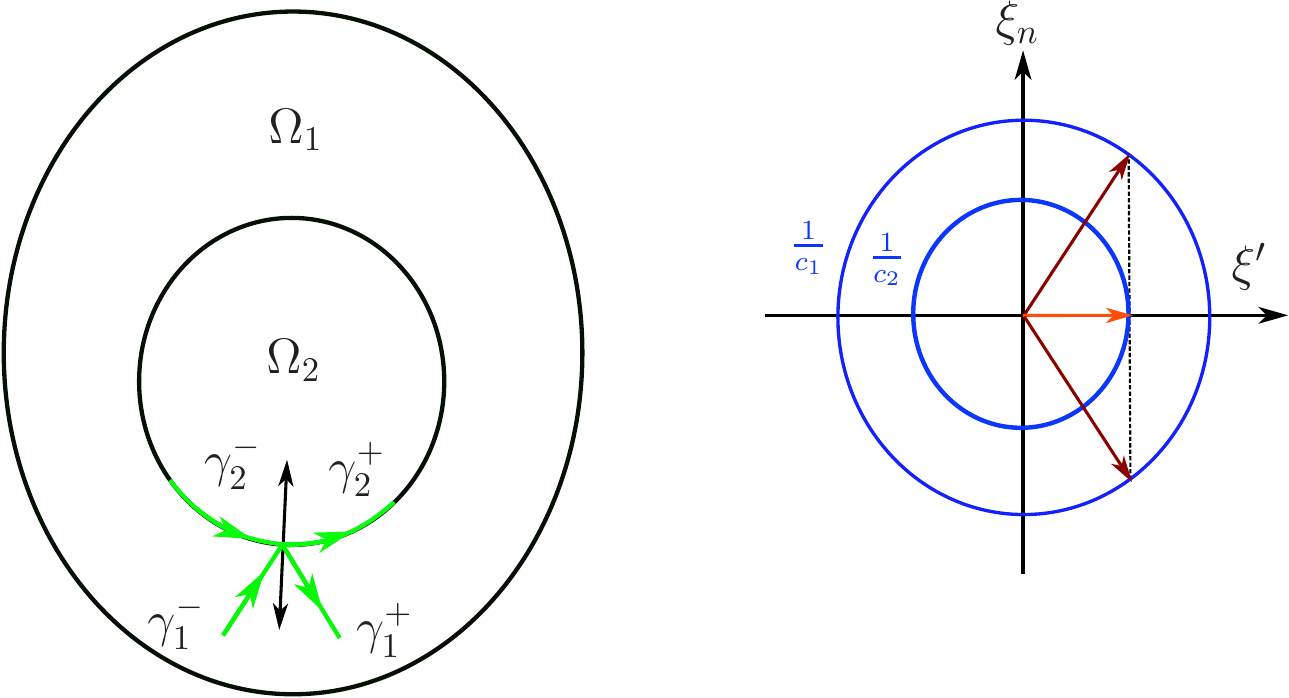} 
\end{center}
\caption{\footnotesize{Rays $(t,x(t))$ of the bicharacteristics near $(\rho_1,\rho_2)\in \H^1\times \G^{2,+}$}}\label{bicharHG}
\end{figure}

If an incoming bicharacteristic from $\Omega_1$ intersects $\d \Omega_2$ traversally at $x(\sigma_0) \in \d \Omega_2$ for some $\sigma_0 \in \RR$ at the critical angle $\sin \theta_1=c_1/c_2$, then $p_1=0$ defined by\eqref{symbol} admits two real roots. The dual variable $\xi_2=(\xi_2^n,\xi_2')$ is such that $\xi_2^n=0$ and $\xi_1'=\xi_2'$. Using the same notations than for the reflected-transmitted case for the half-bicharacteristics $\gamma_1^\pm$, the generalized bicharacteristic path is such that 
\[
\gamma(\sigma,\rho_1)=\begin{cases} \gamma_1^+(\sigma,\rho_1), & 0 < \sigma < \epsilon, \\ 
\widetilde{\gamma}_2(\sigma,\rho_2), & 0 < \sigma < \epsilon, \\ 
 \gamma_1^-(\sigma,\rho_1), & -\epsilon < \sigma <0,
 \end{cases}
\]
where $\widetilde{\gamma}_2$ is defined as the gliding ray for $\Omega_2$ following the definition \eqref{gliding} of gliding rays on $\d \Omega$. Figure \ref{bicharHG} illustrates the outgoing bicharacteristics from an incoming bicharacteristic $\gamma_1^-$ or $\widetilde{\gamma}_2^-$ at a point $(\rho_1,\rho_2)\in \H^1 \times \G^{2,+}$.  \newline

\textit{Reflection in $\Omega_1$: $(\rho_1,\rho_2) \in \H^1\times \E^2$.}\newline

\begin{figure}[!ht]
\begin{center}
	\includegraphics[height=5cm]{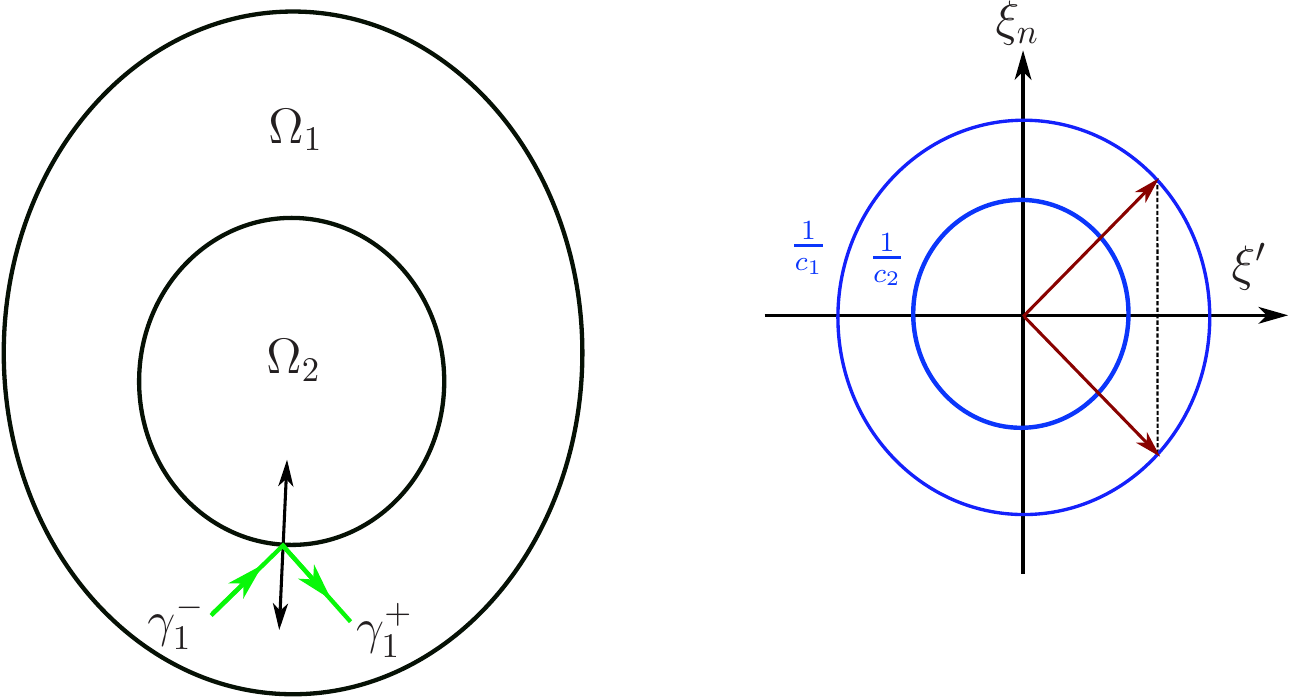} 
\end{center}
\caption{\footnotesize{Rays $(t,x(t))$ of the bicharacteristics near $(\rho_1,\rho_2)\in \H^1\times \E^2$}}\label{bicharHE}
\end{figure}

If the angle of incidence of the incoming bicharacteristic from $\Omega_1$ is strictly greater than the critical angle 
\begin{equation}\label{anglecritique}
\sin \theta_1=c_1/c_2,
\end{equation}
it is reflected in $\Omega_1$ according to \eqref{bicharrefl} and \eqref{bicharreflec} without transmission to $\Omega_2$. In this case the domain $\Omega_2$ acts as an obstacle. Figure \ref{bicharHE} illustrates an example of reflection in $\Omega_1$ without transmission in $\Omega_2$. \newline

\textit{Strictly diffractive rays : $(\rho_1,\rho_2) \in \G^{1,-}\times \E^2$.}\newline

\begin{figure}[!ht]
\begin{center}
	\includegraphics[height=5cm]{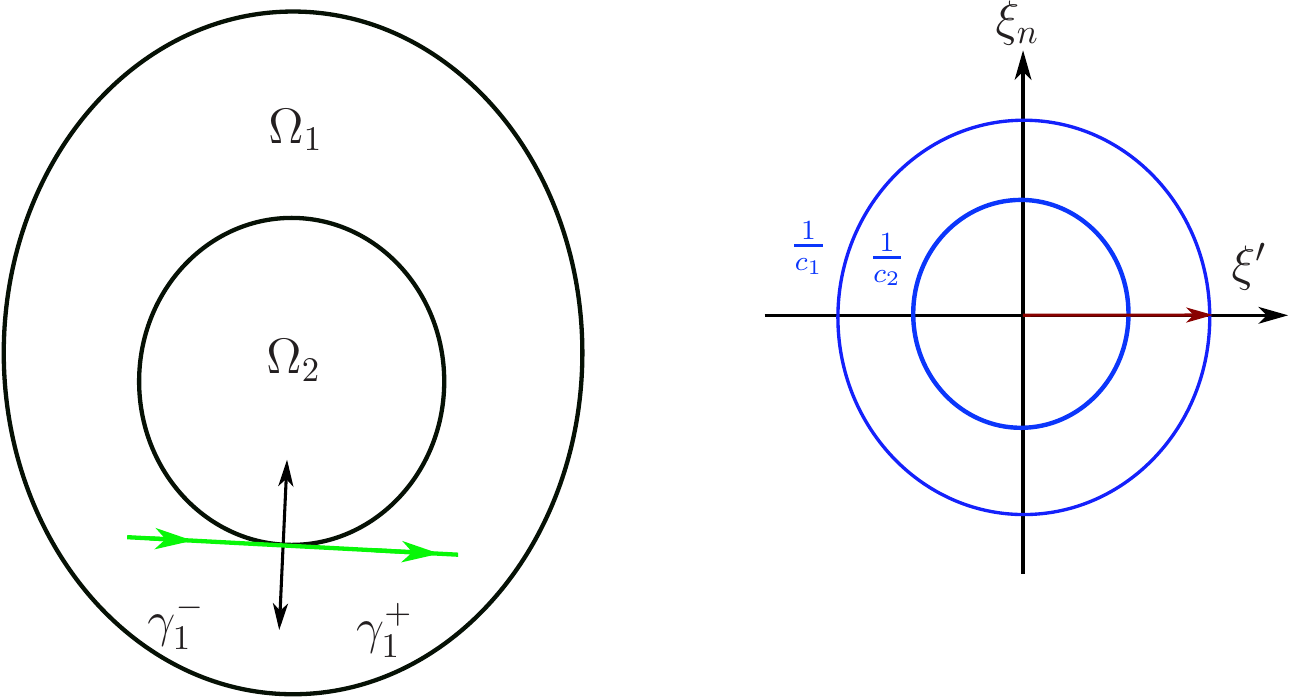} 
\end{center}
\caption{\footnotesize{Rays $(t,x(t))$ of the bicharacteristics near $(\rho_1,\rho_2)\in \G^{1,-}\times \E^2$}}\label{bicharGE}
\end{figure}

In this case, the contact between the incoming ray and $\d \Omega_2$ is \textit{strictly diffractive}, meaning that the incoming bicharacteristic from $\Omega_1$ defined over $T^*(\RR_t \times \RR^2)$ stays in $\Omega_1$ for $|\sigma-\sigma_0|>\epsilon$ for $\epsilon>0$ small and where $\sigma_0$ is the moment when the bicharacteristic interesected $\d \Omega_1$. In this case, there are no transmission in $\Omega_2$ under the hypothesis $c_2>c_1$. Figure \ref{bicharGE} illustrates an example of a glancing ray. \newline

\textbf{Full description of the generalized bicharacteristic flow of \eqref{wave}-\eqref{Dirhom}.}\newline

We gather the description above to describe the generalized bicharacteristic flow of \eqref{wave}-\eqref{Dirhom}. The bicharacteristics propagate in straight line and constant speed $c_i$ in $\Omega_i$. They are connected to the broken bicharacteristics near $\d \Omega$ and $\d \Omega_2$ according to their angle of incidence. As there are no tangential angle of incidence for $\d \Omega$ under the strict convexity assumption of $\Omega$, the gliding ray of $\d \Omega$ is not connected to the bicharacteristic flow in $\Omega_1$. Therefore a gliding ray on $\d \Omega$ exists if and only if it started as a gliding ray on $\d \Omega$. The gliding ray of $\d \Omega_2$ is connected to the bicharacteristic flow in $\Omega$ by the critical angle of incidence \eqref{anglecritique}. A gliding ray on $\d \Omega_2$ emits transmitted rays $\gamma_1^+$ (according to the description $(\rho_1,\rho_2)\in H^1\times \G^{2,+}$) along its trajectory on $\d \Omega_2$. 

\subsection{Microlocal defect measures}\label{SecDefMeasure}

If $(u_i^k)_{k\in \NN}$ are bounded sequences of solutions of \eqref{wave}-\eqref{Dirhom} weakly converging to $0$ in $H^1_{loc}((0,T)\times \Omega_i)$, then one can associated to each sequences a microlocal defect measure $\mu_i$ which are orthogonal in the sense of measure theory (\cite{LameBG}, \cite[Lemma 3.30]{Duyck}). Moreover, the measures $\mu_1$ and $\mu_2$ are each supported in $\textrm{Char} P_1\cup \H \cup \G^+ \cup \H^1 \cup \G^{1,-}$ and $\textrm{Char} P_2\cup \H^2 \cup \G^{2,+}$ respectively. This property is known as the elliptic regularity theorem for the microlocal defect measures \cite{LameBG}. 

We now analyse the property of these measures. In the interior $T^*(\RR_t \times \Omega_i)$, since the solutions $u^1$ and $u^2$ of \eqref{wave}-\eqref{Dirhom} evolve independently in $\Omega_i, i=1,2$, we make use of the classical measure propagation theorem of Gerard \cite{Gerard}. Near the boundary $\d \Omega$, the analysis of \cite{BLR} yields the propagation theorem for the wave equation with homogeneous Dirichlet boundary condition. 
\begin{cor}\label{ondesbordpropag}
For $\rho_1\in \H \times \G^{+}$, we have the following equivalence
\begin{align*}
\left((\gamma_1^-)\cap \textrm{supp}(\mu_1)\right)= \emptyset \Leftrightarrow \left((\gamma_1^+)\cap \textrm{supp}(\mu_1)\right) = \emptyset. 
\end{align*}
\end{cor}
An adaptation of the results of \cite{LameBG} (see also \cite{LameBJ}) allows the characterisation of the propagation of the support of the measure along the generalized bicharacteristics near $\d \Omega_2$.

\begin{prop}\label{propag}
With the above notations, we have 
\begin{enumerate}
\item If $(\rho_1,\rho_2) \in (\H^1 \cup \G^{1,-} \cup \E^1) \times \E^2$, then $\mu_2=0$ near $\rho_2$. Therefore,
\begin{enumerate}
\item if $\rho_1 \in \E^1$, then $\mu_1=0$ near $\rho_1$,
\item otherwise, $\rho_1 \in \H^1 \cup \G^{1,-}$ and the support of $\mu_1$ propagates from $\gamma_1^-$ to $\gamma_1^+$. 
\end{enumerate}
\item otherwise $\rho_2\in \H^2 \cup \G^{2,+}$ and $\rho_1 \in \H^1$. In this case, if $\gamma_i^- \cap \textrm{supp}(\mu_i) = \emptyset$, for $i=1$ or $i=2$, then for $j=3-i$, the support of $\mu_j$ propagates from $\gamma_j^-$ to $\gamma_j^+$ and from $\gamma_j^-$ to $\gamma_i^+$.  
\end{enumerate}
\end{prop}

As a consequence of the conversion of the total mass (\cite[Proof of Theorem 4]{LameBG}), we obtain

\begin{cor}\label{equivsuppmes}
For $\rho_2\in \H^2 \times \G^{2,+}$, we have the following equivalence
\begin{align*}
\left((\gamma_1^-)\cap \textrm{supp}(\mu_1)\right) &\cup \left((\gamma_2^-)\cap \textrm{supp}(\mu_2)\right) = \emptyset \\
&\Updownarrow \\
\left((\gamma_1^+)\cap \textrm{supp}(\mu_1) \right)&\cup\left( (\gamma_2^+)\cap \textrm{supp}(\mu_2)\right) = \emptyset. 
\end{align*}
\end{cor}

\begin{remark}
\begin{enumerate}
\item If one of the two statements in Corollary \ref{equivsuppmes} holds true, then $\mu_i =0, i=1,2$ near $(\rho_1,\rho_2)$. 
\item By time-reversibility, one also obtain the characterisation of two incoming bicharacteristics from $\Omega_1$ or $\Omega_2$. If these two incoming bicharacteristics of, say $\Omega_1$, do not intersect the support of the measure $\mu_1$, then the outgoing bicharacteristics in $\Omega_2$ do not intersect the support of $\mu_2$.  
\end{enumerate}
\end{remark}

We begin by proving Proposition \ref{propag}.

\beginpf 

1. Let $(\rho_1,\rho_2) \in (\H^1 \cup \G^{1,-} \cup \E^1) \times \E^2$. It is a classical result from the elliptic theory at the boundary that $\mu_2=0$ near $\rho_2 \in \E^2$ (see for instance \cite[Appendix A.1]{LameBG}). Therefore 1.a) and 1.b) follow from the classical propagation results for the wave equation at the boundary. 

2. Consider now $(\rho_1,\rho_2)\in \H^1 \times (\H^2 \cup \G^{2,+})$. We explain how to adapt the results of \cite[Appendix A.2]{LameBG}. We recall the notations of Section \ref{SecDefBichar}. The local geodesic coordinates $(x_n,x')$ are such that, locally, we have $\Omega_2=\{x_n >0\}, \d \Omega_2=\{x_n=0\}$ and $\Omega_1=\{x_n <0\}$. Recall that near this point, 
\[
P_i=c_i^2(\d_{x_n}^2+R(x_n,x',D_{x'}))
\]
where the principal symbol of $R$ is $r(x_n,x',\xi')$ and $r(x_n,x',\xi')\geq c|\xi'|^2$. One can then use \cite[Lemma A.1]{LameBG} to factorise the pseudodifferential operators in two different ways
\begin{align*}
P_i=-c_i^2(D_{x_n}-\Lambda^+(x_n,x',D_{x'}))(D_{x_n}-\Lambda^-(x_n,x',D_{x'}))+T(x_n,x',D_{x'}) \\
P_i=-c_i^2(D_{x_n}-\widetilde{\Lambda}^-(x_n,x',D_{x'}))(D_{x_n}-\widetilde{\Lambda}^+(x_n,x',D_{x'}))+\widetilde{T}(x_n,x',D_{x'}) 
\end{align*}
where $\Lambda^\pm$ and $\widetilde{\Lambda}^\pm$ are tangential pseudodifferential operators of order 1 and such that $\sigma_1(\Lambda^\pm)=\pm \sqrt{r}$ and $\sigma_1(\widetilde{\Lambda}^\pm)=\pm \sqrt{r}$ and where $T$ and $\widetilde{T}$ are tangential pseudodifferential operators of order $-\infty$. A microlocalisation near of $\gamma_i^{\pm}$ is done using $q_0(x',\xi_i')$ a symbol of order $0$ and equal to $1$ in a conical neighborhood of $\rho_i$ and of compact support. Let us remark that at a point $(\rho_1,\rho_2)\in \H^1 \times (\H^2 \cup \G^{2,+})$, the tangential components of $\rho_1=(x_n,x',\xi_1^n,\xi_1')$ and $\rho_2=(x,y,\xi_2^n,\xi_2')$ are equal : $(x',\xi_1')=(x',\xi_2')$. Therefore, the same symbol $q_0$ may be used for the microlocalisation near $\rho_1$ and $\rho_2$. The symbol is then propagated by the Hamiltonian 
\[
(c_i\d_{x_n} \mp H_{\sqrt{r}})q_i^{\pm}=0, \quad \restriction{q_i^{\pm}}{x_n=0}=q_0. 
\]
Consider $\varphi \in C_0^\infty( \RR_x)$ to be equal to $1$ near $0$ and of compact support near $0$. If we denote $Q_i^\pm = \textrm{Op}(\varphi q_i^\pm)$ and $Q_0 = \textrm{Op}(\varphi q_0)$, then we obtain the same results as in \cite[Appendice A.2]{LameBG}, that is, if we consider bounded sequences $(u_1^k,u_2^k) \subset H^1(]-1,0[\times Y) \times H^1(]0,1[\times Y)$, where $Y=\{ x' \, | \, |x'| <1 \} $, such that 
\begin{align*}
P_1 u_1& \rightarrow 0 \textrm{ in } L^2(]-1,0[\times Y), \\
P_2 u_2& \rightarrow 0 \textrm{ in } L^2(]0,1[\times Y).
\end{align*}
and if we suppose, without loss of generality, that $\gamma_1^+ \cup \textrm{supp}(\mu_1) = \emptyset$, then (recall that we can deduce that $\restriction{u_i^k}{x=0}$ and $\restriction{D_x u_i^k}{x=0}$ are bounded in $H^1_{\rho_i}$ and $L^2_{\rho_i}$ respectively) 
\[
c_1^2Q_0(\restriction{D_{x_n} u_1^k}{x_n=0}-\Lambda^- \restriction{u_1^k}{x_n=0}) \rightarrow 0 \textrm{ in } L^2(Y). 
\]
Applying the tangential pseudodifferential operator $\Lambda^-$ to the Dirichlet boundary condition yields
\[
\Lambda^-u_1^k=\Lambda^-u_2^k, \textrm{ in } L^2(Y),
\]
and a Lopatinski condition on $u_2^k$ is obtained
\begin{equation}\label{lopau2}
c_1^2Q_0\left(\restriction{\dfrac{c_2^2 D_{x_n} u_2^k}{c_1^2}}{x_n=0}-\Lambda^- \restriction{u_2^k}{x_n=0}\right) \rightarrow 0 \textrm{ in } L^2(Y). 
\end{equation}
In particular, the traces of solution $u_k^2$ are decoupled from the traces of $u_k^1$. Therefore we conclude by the classical propagation result if $\gamma_1^- \cap \textrm{supp}(\mu_1)=\emptyset$. If $\gamma_2^- \cap \textrm{supp}(\mu_2)=\emptyset$, we use the following lemma to conclude

\begin{lem}
Let $(\rho_1,\rho_2) \in \H^1 \times (\H^2 \cup \G^{2,-})$. Suppose $\gamma_1^+ \cup \textrm{supp}(\mu_1) = \emptyset$. Then neighborhoods of $\rho_2$ intersect $\textrm{supp}(\mu_2)$ if and only if it intersects $ \textrm{supp}(\textrm{meas}(\restriction{u_2}{x=0})) \cup  \textrm{supp}(\textrm{meas}(\restriction{D_{n_2}u_2}{x=0}))$
\end{lem}

\beginpf

We only need to prove the only if part and the result comes from \eqref{lopau2}. 

\endpf

This conclude the proof of Proposition \ref{propag}. 

\endpf

We now turn to the proof of Corollary \ref{equivsuppmes}

\beginpf 

Let us consider the case where $\gamma_1^-$ and $\gamma_2^-$ intersects the boundary $\d \Omega_2$ at a $\H^1\times \H^2$ point. The case $\H^1\times \G^{2,+}$ is treated similarly. Consider $\left((\gamma_1^-)\cap \textrm{supp}(\mu_1)\right) \cup \left((\gamma_2^-)\cap \textrm{supp}(\mu_2)\right) = \emptyset$, we obtain from the proof of Proposition \ref{propag}
\begin{equation}\label{reltrace}
c_i^2Q_0(\restriction{D_x u_i^k}{x=0}-\Lambda^- \restriction{u_i^k}{x=0}) \rightarrow 0 \textrm{ in } L^2(Y). 
\end{equation}
Using the Neumann condition yields
\[
(c_1^2-c_2^2)Q_0(\Lambda^- \restriction{u_i^k}{x=0}) \rightarrow 0 \textrm{ in } L^2(Y), 
\] 
which is a tangential pseudodifferential operator of order 1 
and elliptic at $\rho_i$. It follows that 
\[
c_i^2Q_0\restriction{D_x u_i^k}{x=0} \rightarrow 0 \textrm{ in } L^2(Y). 
\]
The ellipticity implies that $u_i^k$ converges to 0 in $H^1$ and the standard propagation results implies Corollary \ref{equivsuppmes}. 

\endpf

The following corollary follows from Proposition \ref{propag} and Corollary \ref{equivsuppmes}. We provide a sufficient condition for the propagation of the measure of the generalized bicharacteristic flow of \eqref{wave}-\eqref{Dirhom} as illustrated in figure \ref{figurerecurs}. To this end, let us introduce some notations. 
%

Let $(\rho_1,\rho_2) \in T^*(\d \Omega_2 \times (0,T))$. We denote $(\rho_1^1,\rho_2^1) \in T^*(\d \Omega_2 \times (0,T))$ the coordinates of the intersection of $\gamma_1^+$ with $\d \Omega_2$ (after reflection on $\d \Omega$). Denote $(\rho_1^2,\rho_2^2) \in  T^*(\d \Omega_2 \times (0,T))$ the coordinates of the intersection of $\gamma_2^+$ with $\d \Omega_2$. Denote $(\rho_1^{-a},\rho_2^{-a}) \in  T^*(\d \Omega_2 \times (0,T))$, for $a=1,2$ the intersection of $\gamma_a^-$ with $\d \Omega_2$. Finally, for a finite sequence $a_i=\pm 1, \pm 2, 1\leq i \leq n$, denote $(\rho_1^{a_1\ldots a_n},\rho_2^{a_1\ldots a_n})$ the iteration of this procedure with the requirement that a sequence possesses consecutive $2$ and $-2$ (or the converse) if and only if the ray is a trapped ray of $\Omega_2$. Otherwise consecutives $1$ and $-1$ (or the converse) or $2$ and $-2$ (or the converse) consists to following a ray say $\gamma_1^+$ forward in time and then backward in time after the intersection with $\d \Omega_2$ which is not described by the bicharcteristic flow of \eqref{wave}-\eqref{Dirhom}.

\begin{cor}\label{equivsuppmesgraph}

For $(\rho_1,\rho_2) \in T^*(\d \Omega_2 \times (0,T))$, if there exists a finite sequence $\{a_i\}_{1\leq i \leq n},a_i=\pm 1, \pm 2$ such that for every $1\leq i \leq n-1$ there exists a half-ray $\gamma_{|b|}^{\textrm{sign}(b)}(\rho_1^{a_1\ldots a_i},\rho_2^{a_1\ldots a_i})\cup \textrm{supp}(\mu_{|b|})= \emptyset$ where $b=\pm 1$ or $b=\pm 2$ is such that $b\neq -a_i$ or $b\neq a_{i+1}$ and such that there exist two half-rays $\gamma_{|b_j|}^{\textrm{sign}(b_j)}(\rho_1^{a_1\ldots a_n},\rho_2^{a_1\ldots a_n})\cup \textrm{supp}(\mu_{|b_j|})= \emptyset$ where $j=1,2$, $b_j = \pm 1$ or $ b_j = \pm 2$ are such that $b_j \neq -a_n, j=1,2$, then $\mu_1,\mu_2=0$ near $(\rho_1,\rho_2)$. 
\end{cor}

\beginpf

Proposition \ref{propag} allows to use a recursive argument on the sequence of half-rays defined by the sequence $\{a_i\}_{1\leq i \leq n}$. We prove the result for $n=1$. Since two-half rays do not intersect the support of their respective measure near $(\rho_1^{a_1}, \rho_2^{a_1})$, we conclude by Corollary \ref{equivsuppmes} that $\mu_1,\mu_2 =0$ near $(\rho_1^{a_1}, \rho_2^{a_1})$. Corollary \ref{equivsuppmes} implies $\gamma_{|a_1|}^{\textrm{sign}(-a_1)}(\rho_1^{a_1}, \rho_2^{a_1}) \cup \textrm{supp}(\mu_{|a_1|}) = \emptyset$. Moreover $\gamma_{|a_1|}^{\textrm{sign}(-a_1)}(\rho_1^{a_1}, \rho_2^{a_1})$ intersects $\d \Omega_2$ at $(\rho_1, \rho_2)$ by definition. By hypothesis, there exists another half-ray that do not intersect the support of its respective measure. Therefore we conclude this case by Corollary \ref{equivsuppmes}. The general case follows easily by iterating this process.

\endpf

\section{Proof of the observability of waves equations with a transmission condition}\label{Secmain}

We provide in this section the proof of Theorem \ref{main} and Theorem \ref{thmtrap} based on the geometrical construction introduced in Section \ref{SecStrategy}. We begin by proving Theorem \ref{main} which is divided in several steps. The proof of Theorem \ref{thmtrap} is then presented, which only differs from the proof of Theorem \ref{main} at the last step. 

\beginpf

\textbf{Step 1 : Definition of $\Gamma_1^0$}\newline

Let $0\leq s_1< s_2 < 1$ be such that $\{\delta( s_1), \delta( s_2) \} \in \d \Gamma(x_0)$. We denote
\begin{equation}\label{tangentOmega2}
L(\delta(s_i)):=\max\left\{ \left. \dfrac{\<(x-\delta(s_i)),(-1)^{i+1}\delta'(s_i)\>}{\|x-\delta(s_i)\|} \, \right| \, \forall x\in \d \Omega_2 \right\}, \quad i=1,2,
\end{equation}
\begin{figure}[!ht]
\begin{center}
	\includegraphics[height=5cm]{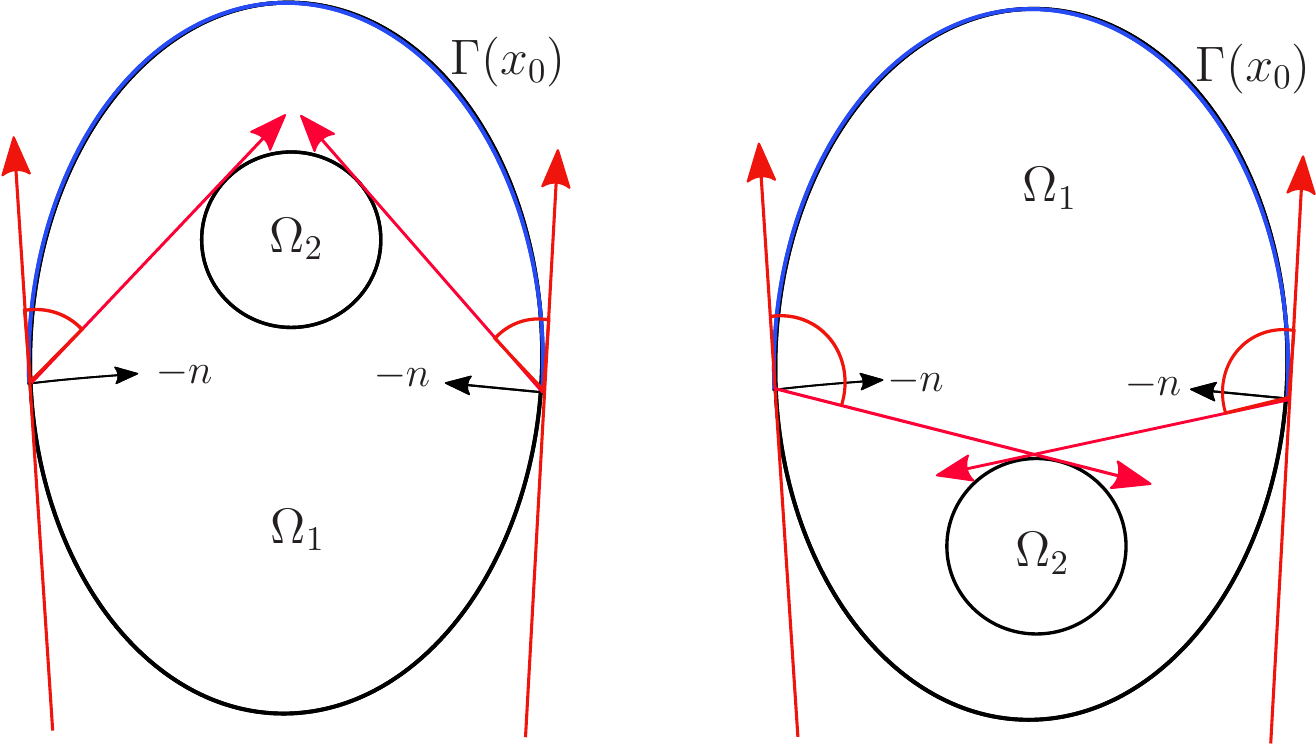}
\end{center}
\caption{\footnotesize{Localisation of $\Omega_2$. Case $L(\delta(s_i))<0, i=1,2$ on the left and $L(\delta(s_i))>0, i=1,2$ on the right.}}\label{Omega2inclus2}
\end{figure}

\begin{figure}[!ht]
\begin{center}
	\includegraphics[height=5cm]{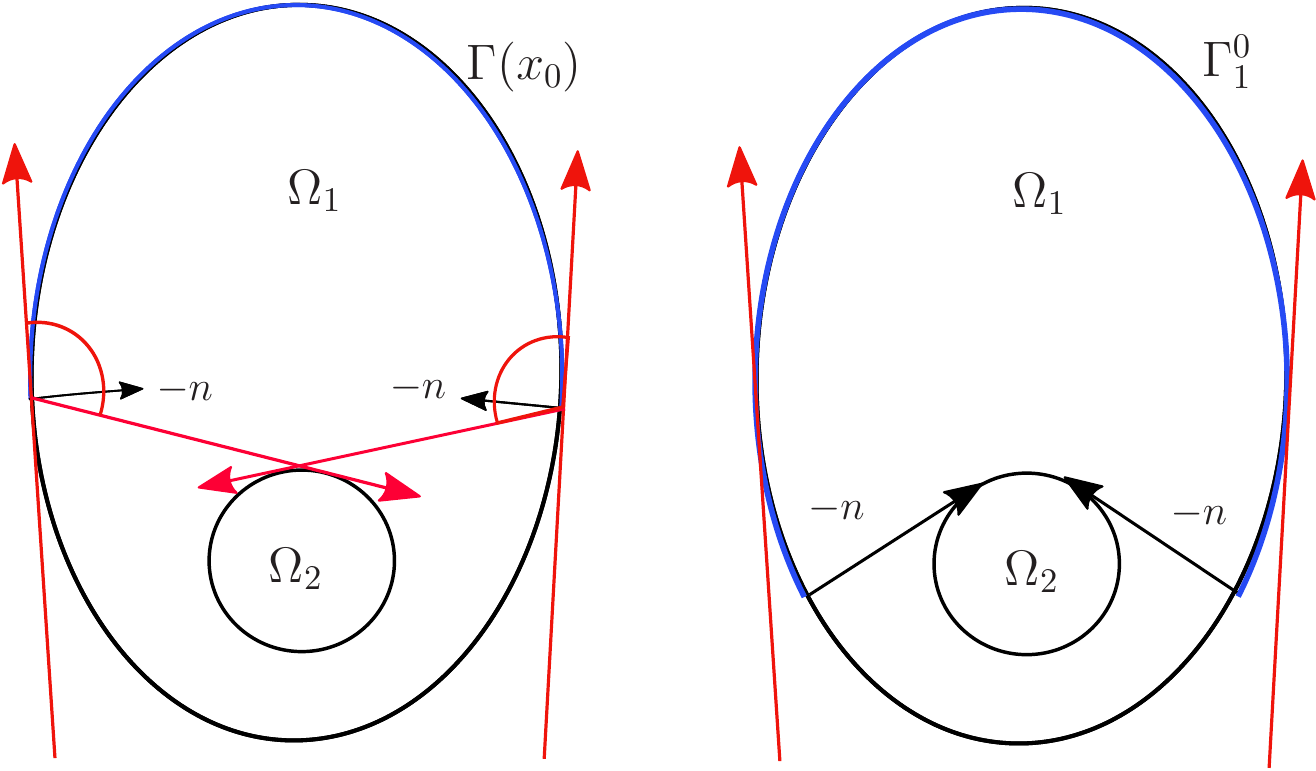} 
\end{center}
\caption{\footnotesize{Left : case where $\Omega_2$ is not included in $\Gamma(x_0)$. Right : extension of $\Gamma(x_0)$ to $\Gamma_1^0$}}\label{PropagG2}
\end{figure}

Recall that $\delta'$ is in the direction of the propagation. If $L(\delta(s_i)) \leq 0$ for $ i=1$ (resp. $i=2$), then we define $s_1^{1,0}:=s_1$ (resp. $s_2^{1,0}:=s_2$). Otherwise, if $L(\delta(s_i)) > 0$ for $i=1$ (resp. $i=2$), denote $s_1^{2,0}$ (resp. $s_2^{2,0}$) the largest $s<s_1$ (resp. smallest $s>s_2$) such that $L(\gamma(s))=0$. With this definition of $s_i^{1,0},i=1,2$, we define 
\[
\Gamma_1^0:=\{ \delta(s) \, | \, s_1^{1,0}< s < s_2^{1,0} \}.
\]
Figure \ref{PropagG2} illustrate the extension of $\Gamma_1^0$. 

We prove the following for $\Gamma_1^0$. 
\begin{lem}\label{Lemtildegammax0}
Neighborhoods of points $\rho_1 \in T^*(\Gamma_1^0\times (0,T))$ do not belong to the support of the measure $\mu_1$.
\end{lem}

\beginpf

\begin{figure}
\begin{center}
\begin{tabular}{c}
 \mbox{\includegraphics[height=4cm]{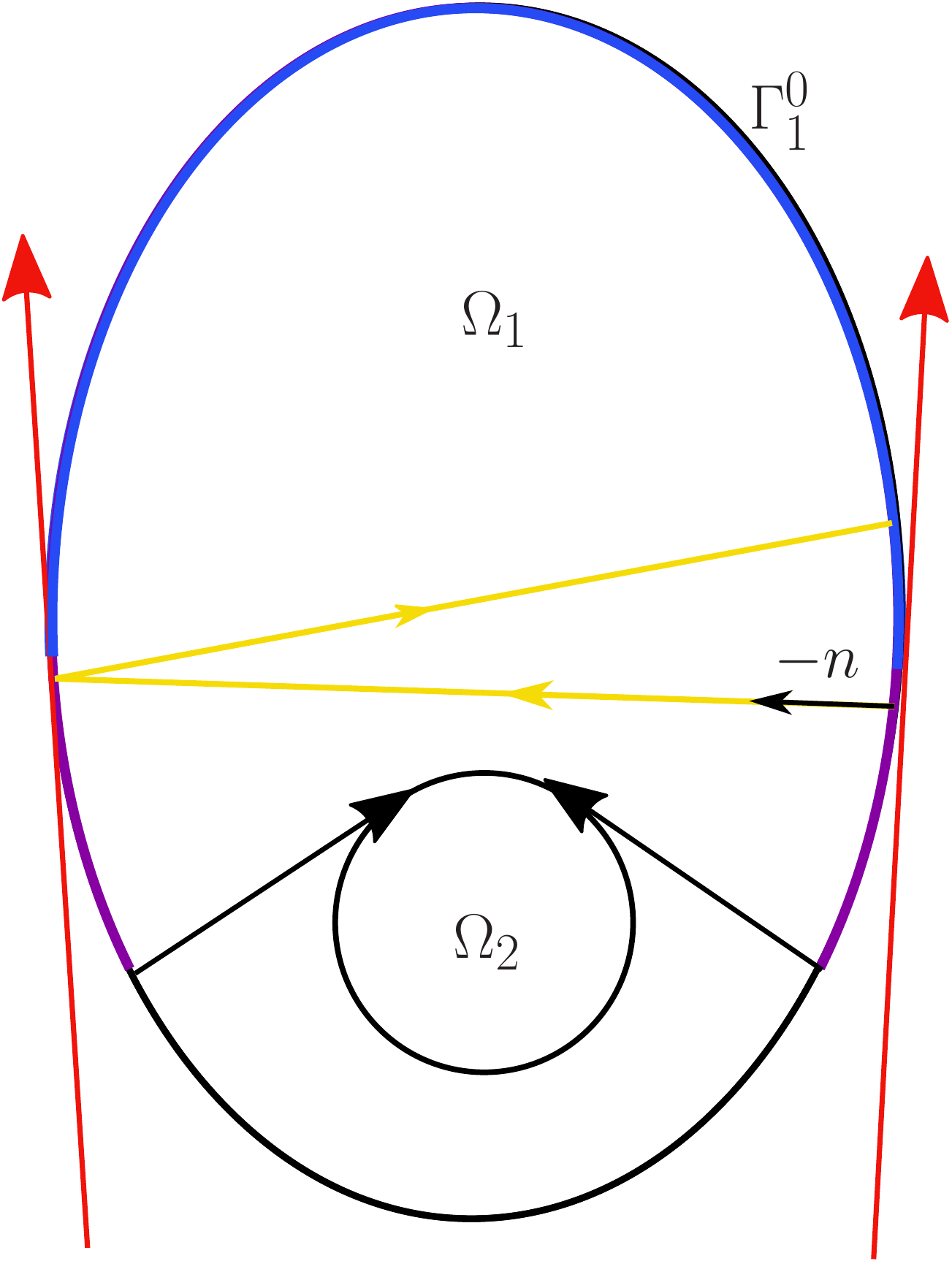}} \quad 
 \mbox{\includegraphics[height=4cm]{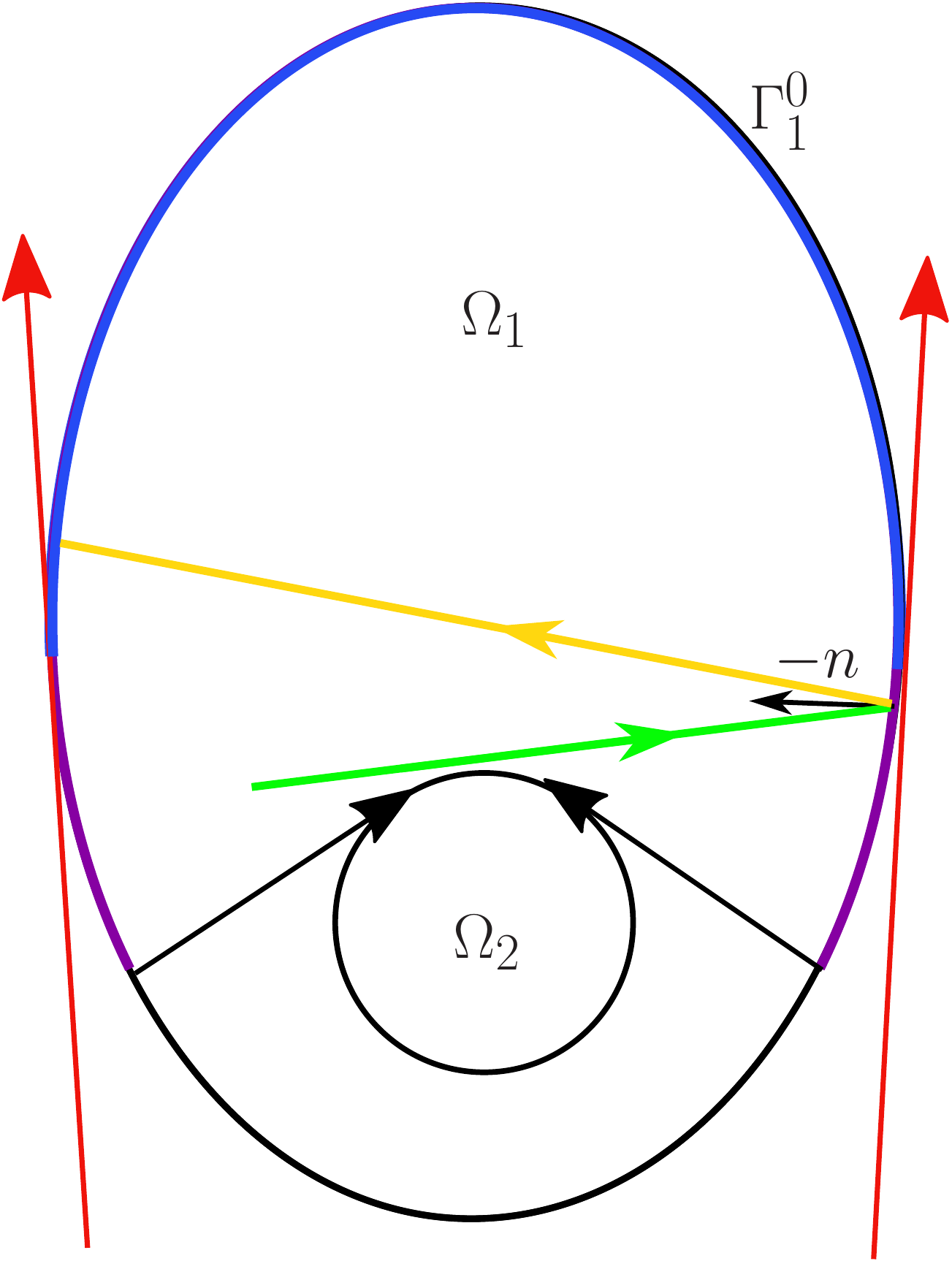}} 
\end{tabular}
\end{center}
\caption{Left : ray in the $-n$ direction. Right : propagation of the half-ray toward $\Gamma(x_0)$. In purple is the $\Gamma_1^0 \setminus \Gamma(x_0)$ region.}\label{S1}
\end{figure}

We consider $\rho_1=(t,x,\tau_1,\xi_1) \in T^*\left(\Gamma_1^0 \setminus \Gamma(x_0) \times (0,T)\right)$. The standard elliptic theory implies that $\mu_1=0$ near $\rho_1 \in \E$ and the gliding ray of $\d \Omega$ eventually reaches $\Gamma(x_0)$ so $\mu_1=0$ near $\rho_1\in \G^+$. Therefore, suppose $\rho_1 \in \H$. In this case, the conical assumption of $\Gamma(x_0)$ ensures that for every outgoing rays from $\rho_1=(t,x,\tau_1,-n(x)/c_1) \in \H$, there exists $n\in \NN^*$ such that $\Pi_x(\F^n(x,-n(x)/c_1))\in \Gamma(x_0)$ and the bicharacteristic intersects $\Gamma(x_0)$ in a non-diffractive way. Notice that by construction, no such rays intersects $\d \Omega_2$. Therefore, by Corollary \ref{ondesbordpropag}, one obtains that $\mu_1=0$ near $\rho_1=(t,x,\tau_1,-n(x)/c_1) \in \H$. In the general case $\rho_1= (t,x,\tau_1,\xi_1)\in \H$ it suffices to notice that there is always a half-ray ($\xi_1'>0$ if $x=\delta(s)$ is such that $s_1^{1,0}\leq s \leq s_1$ and $\xi_1'<0$ if $x=\delta(s)$ is such that $s_2^{1,0}\geq s \geq s_2$) confined between $\Gamma_1^0$ and the bicharacteristic starting from $x$ in the $-n(x)/c_1$ direction (see figure \ref{S1}). Therefore this half-ray intersects $\Gamma(x_0)$ in a non-diffractive way and one concludes with Corollary \ref{ondesbordpropag}.  

\endpf

\textbf{Step 2 : Definition of $\Gamma_2^0$}\newline

We define $0\leq s_1^{2,0}, s_2^{2,0} < 1$ to be such that $\delta_2(s_i^{2,0})=\textrm{argmax}(L(\delta(s_i^{1,0}))), i=1,2$. Since $\Omega$ and $\Omega_2$ are strictly convex, the existence and uniqueness of such points is ensured. Moreover, $0\leq s_1^{2,0}< s_2^{2,0} < 1$. We define 
\begin{equation}\label{defg20}
\Gamma_2^0:=\{ \delta_2(s) \, | \, s_1^{2,0}< s < s_2^{2,0} \}.
\end{equation}

We prove the following.

\begin{lem}\label{obsg20}
Neighborhoods of points $(\rho_1,\rho_2) \in T^*(\Gamma_2^0 \times (0,T))$ do not belong to the support of the measures $\mu_1,\mu_2$.
\end{lem}

First, we prove that 
\begin{lem}
Let $\Gamma_2^0$ be given by \eqref{defg20}. Then $\Gamma_2^0\neq \emptyset$. 
\end{lem}

\beginpf 

It follows easily from the strict convexity of $\Omega_2$ in the case where $L(\delta(s_i))<0$ for $i=1$ or $i=2$ so we focus on the the case where $L(\delta(s_i))=0, i=1,2$. But in this case, by definition \eqref{tangentOmega2} of $L$, $\Pi_x(\F(\delta(s_i),-n(\delta(s_i))/c_i))=\delta(s_j), i+j=3$ and $\Pi_x(\F^2(\delta(s_i),-n(\delta(s_i))/c_i))=\delta(s_i), i=1,2$. Therefore this ray is trapped and is not observed by $\Gamma(x_0)$, which is a contradiction with the definition of $\Gamma(x_0)$ and the fact that it implies GCC.

\endpf

We proceed with the proof of Lemma \ref{obsg20}.

\beginpf

\begin{figure}
\begin{center}
\begin{tabular}{c}
 \mbox{\includegraphics[height=5cm]{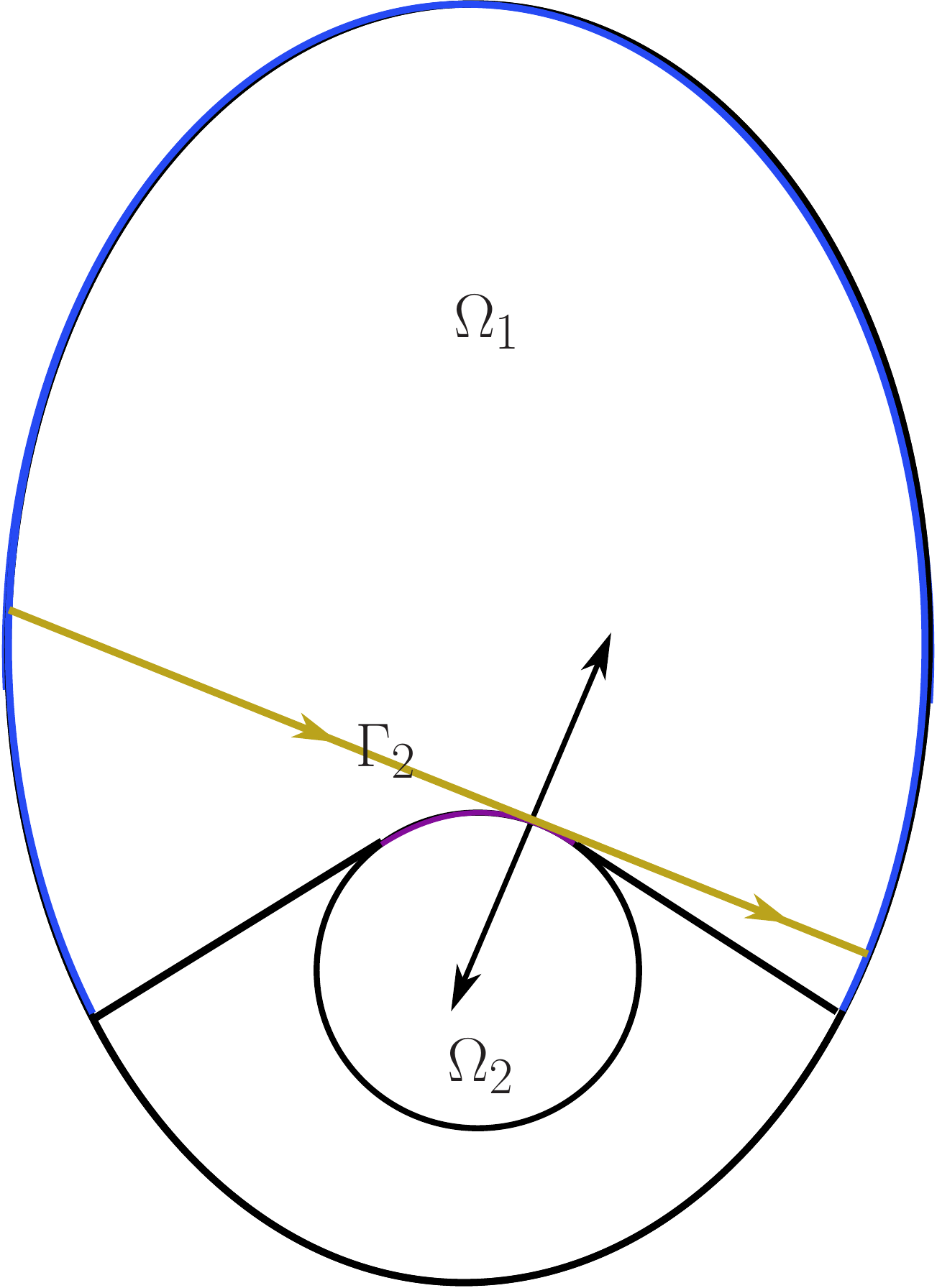}} 
\quad
 \mbox{\includegraphics[height=5cm]{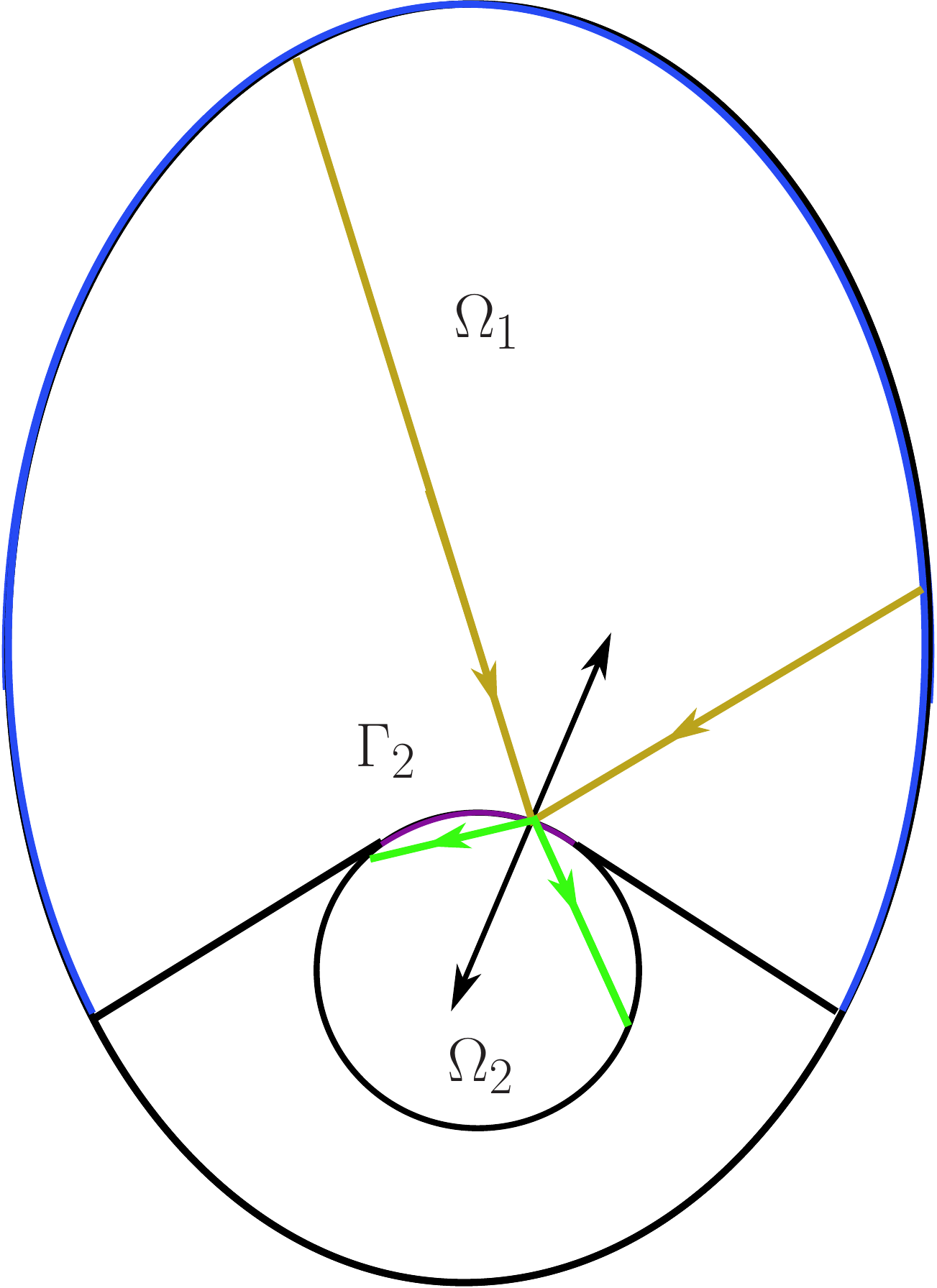}} 
\end{tabular}
\end{center}
\caption{Left : the two-half rays of $(\rho_1,\rho_2)\in \G^{1,-}\times \E^2$ intersect $\Gamma_1^0$ non-diffractively. Right : the two half-rays $\gamma_1^\pm$ of $(\rho_1,\rho_2)\in \H^1 \times \H^2$ intersect $\Gamma_1^0$ non-diffractively.}\label{S2}
\end{figure}

By construction, the line $l(\delta(s_i^{1,0}),\delta_2(s_i^{2,0}))$ is tangent at $\delta_2(s_i^{2,0}) \in \d \Omega_2$ for $i=1,2$. Moreover, the bicharacteristic starting from $\delta(s_i^{1,0})$ and passing through $\delta_2(s_i^{2,0})$ intersects $\Gamma_1^0$ non-diffractively. Therefore, the two half-bicharacteristics $\gamma_1^-(\rho_1,\rho_2)$ and $\gamma_1^+(\rho_1,\rho_2)$ of $(\rho_1,\rho_2)\in \G^{1,-}\times \E^2 \subset T^*(\Gamma_2^0 \times (0,T))$ intersect $\Gamma_1^0$ in a non-diffractive way (see figure \ref{S2}). Proposition \ref{propag} implies that $\mu_1, \mu_2=0$ near $(\rho_1,\rho_2)\in \G^{1,-}\times \E^2 \subset T^*(\Gamma_2^0 \times (0,T))$. It also follows that $\mu_1, \mu_2=0$ near $(\rho_1,\rho_2)\in \H^1\times (\H^2 \cup \G^{2,+} \cup \E^2) \subset T^*(\Gamma_2^0 \times (0,T))$ since both half-bicharacteristics $\gamma_1^+(\rho_1)$ and $\gamma_1^-(\rho_1)$ intersects $\Gamma_1^0$ non-diffractively. Finally, $\mu_1, \mu_2=0$ near $(\rho_1,\rho_2)\in \E^1 \times  \E^2 \subset T^*(\Gamma_2^0 \times (0,T))$ follows from the classical elliptic theory.  

\endpf

\textbf{Step 3. Extension of $\Gamma_2^{n-1}$ to $\Gamma_2^{n}$}\newline

\begin{figure}
\begin{center}
\begin{tabular}{c}
 \mbox{\includegraphics[height=5cm]{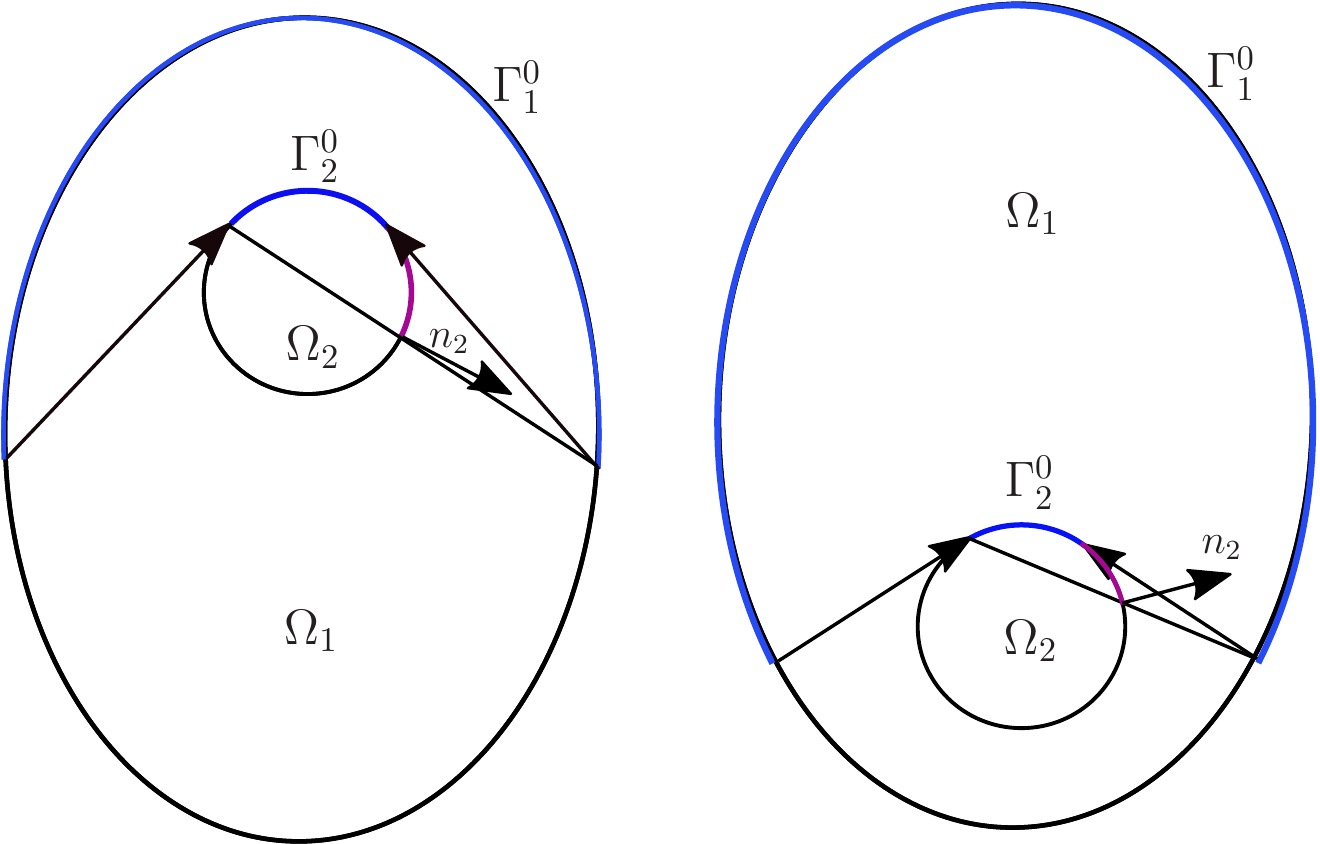}} 
\end{tabular}
\end{center}
\caption{Extension of one of the boundary of $\Gamma_2^{n-1}$ to $\Gamma_2^{n}$ in two different cases.}\label{S3}
\end{figure}

For $n\in \NN$, consider the line $l(\delta(s_1^{1,n-1}),\delta_2(s_2^{2,n-1}))$ (resp. $l(\delta(s_2^{1,n-1}),\delta_2(s_1^{2,n-1}))$). Then there is at most one intersection of this line with $\d \Omega_2$, aside from $\delta_2(s_2^{2,n-1})$ (resp. $\delta_2(s_2^{1,n-1})$). When it exists, we denote this intersection $\delta_2(\check{s}_1^{2,n})$ with $\check{s}_1^{2,n}\in [0,1]$ (resp. $\delta_2(\check{s}_2^{2,n})$ with $\check{s}_2^{2,n}\in [0,1]$). The uniqueness of the intersection is guaranteed by the strict convexity of $\Omega_2$. When such an intersection exists, we extend $\Gamma_2^{n-1}$ to $\Gamma_2^n$ such that 
\begin{equation}\label{defntildeg2}
\Gamma_2^n:=\{ \delta_2(s) \, | \, \check{s}_1^{2,n} < s < \check{s}_2^{2,n} \textrm{ and } \Pi_x(\F(\delta_2(s), n_2(\delta_2(s))) \in \Gamma_1^{n-1}\}, 
\end{equation}
and such that $\Gamma_2^n \supseteq \Gamma_2^{n-1}$ is open and connected. The boundary of $\Gamma_2^n$ is denoted $\delta_2(s_i^{2,n}), 0\leq s_1^{2,n}<s_2^{2,n}<1$. See figure \ref{S3} for a representation of the extension on one side of $\Gamma_2^{n-1}$.

We proceed to prove.
\begin{lem}\label{Lemtildegamma2}
Neighborhoods of points $(\rho_1,\rho_2)\in T^*(\Gamma_2^n \times (0,T)), n\in \NN^*$ do not belong to the support of the measures $\mu_1,\mu_2$
\end{lem}
 
\beginpf
 
Consider $\rho_1=(t,\delta_2(s),\tau_1,\xi_1)\in T^*(\Gamma_2^n \times (0,T))$ with $s_1^{2,n} \leq s < s_1^{2,n-1}$. The case $s_2^{2,n-1} \leq s < s_2^{2,n}$ is covered in the same fashion and the case $s_1^{2,n-1}< s < s_2^{2,n-1}$ was covered in the previous iteration. 

\begin{figure}[!ht]
\begin{center}
	\includegraphics[height=5cm]{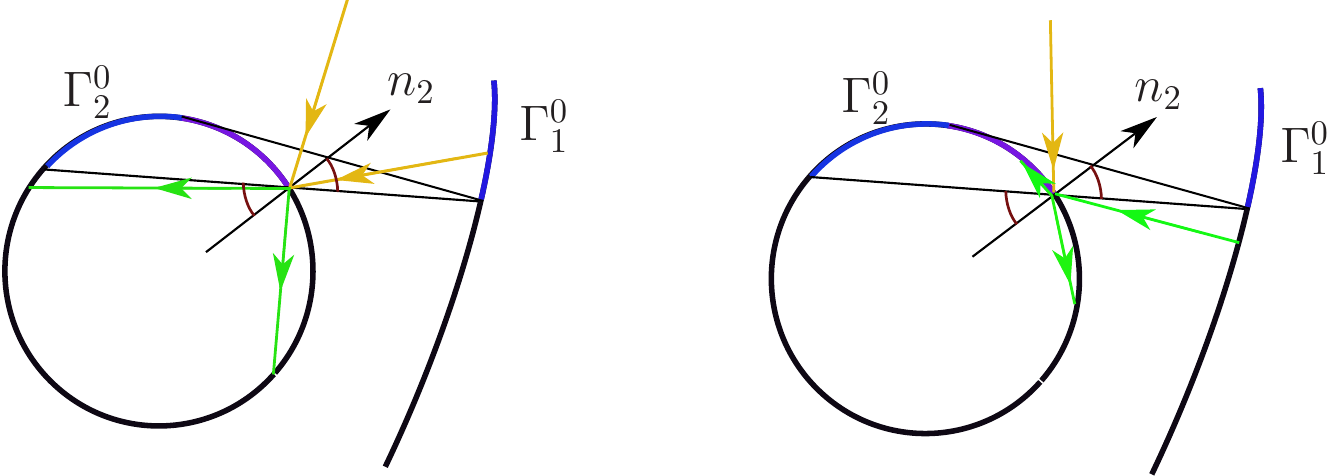}
\end{center}
\caption{\footnotesize{Left : both half-rays (in yellon) in $\Omega_1$ intersect $\Gamma_1$. Right, an incoming half-ray from $\d \Omega_1 \setminus \Gamma_1$. In this case, one of the transmitted half-ray is confined the region bounded by the line $l(\delta_2(s_2^{2,n-1}),\delta(s_1^{1,n-1}))$ and $\hat{\Gamma}_2^{n}$.}}\label{n2Gamma1}
\end{figure}

By definition of $\Gamma_2^n$, $\Pi_x(\F(\delta_2(s), n_2(\delta_2(s))/c_1))\in \Gamma_1^{n-1}$ and the intersection of the bicharacteristic is non-diffractive. Therefore, the half-bicharacteristic propagating (for small $\sigma$) in the positive tangential direction, say $\gamma_1^+$, intersects non-diffractively $\Gamma_1^{n-1}$. Notice that this is sufficient to conclude the cases $(\rho_1,\rho_2)\in (\H^1 \cup \G^{1,-}) \times (G^{2,+} \cup \E^2)$ by Corollary \ref{equivsuppmes} since the gliding ray will eventually reach $\Gamma_2^{n-1}$. The case $(\rho_1,\rho_2) \in \E^1 \times \E^2$ follows by the standard elliptic theory. The analysis of the remaining case $(\rho_1,\rho_2)\in \H^1 \times \H^2$ is divided in two cases : either $\gamma_1^-$ (recall \eqref{bicharreflfrontiere2}) intersects  non-diffractively $\Gamma_1^{n-1}$ or $\d \Omega \setminus \Gamma_1^{n-1}$. In the first case, we are able to conclude by Corollary \ref{equivsuppmes} that $\mu_1,\mu_2 =0$ near these points $(\rho_1,\rho_2)\in \H^1 \times \H^2$ (see figure \ref{n2Gamma1}, on the left).  

Otherwise, consider the line $l(\delta_2(s_2^{2,n-1}),\delta(s_1^{1,n-1}))$. Since $\gamma_1^-$ does not intersect $\Gamma_1^{n-1}$, the angle $\theta_1$ between $\gamma_1^-$ and the normal $n_2(\delta_2(s))$ is strictly larger than the angle $\theta$ between the line $l(\delta_2(s_2^{2,n-1}),\delta(s_1^{1,n-1}))$ and $n_2(\delta_2(s))$. The Snell-Descartes law implies that $\theta_2>\theta_1$. Therefore $\gamma_2^+$, the half-bicharacteristic propagating in the positive tangential direction for $\sigma$ small, associated to $\rho_2=(t,x,\tau_2,\xi_2)$ is confined to the region bounded between the line $l(\delta_2(s_2^{2,n-1}),\delta(s_1^{1,n-1}))$ and $\hat{\Gamma}_2^{n}:=\{\delta_2(s) \, | \, s_2^{1,n} < s < s_2^{2,n-1}\} \subset \Gamma_2^{n}$ (see figure \ref{n2Gamma1}, on the right). The half-ray therefore intersects non-diffractively $\hat{\Gamma}_2^{n}$ at a point $\Pi_x(\F_2(x,\xi_2))$. If $\Pi_x(\F_2(x,\xi_2))\in \Gamma_2^{n-1}$, then one concludes by Corollary \ref{equivsuppmes} that $\mu_1,\mu_2 =0$ near $(\rho_1,\rho_2)$. Otherwise, $\Pi_x(\F_2(x,\xi_2))\in \hat{\Gamma}_2^{n} \setminus \Gamma_2^{n-1}$ in which case we have to carry on with the analysis (which is possible thanks to Corollary \ref{equivsuppmesgraph} and $\gamma_1^+$ intersecting $\Gamma_1^{n-1}$). Let us denote $(\rho_{\F_2,1},\rho_{\F_2,2})\in \H^1 \times \H^2$ the coordinates associated to the intersection of $\gamma_2^+$ with $\Gamma_2^{n} \setminus \Gamma_2^{n-1}$. In this case, $\gamma_2^+(\rho_2)=\gamma_2^-(\rho_{\F_2,2})$. 

\begin{figure}[!ht]
\begin{center}
	\includegraphics[height=4cm]{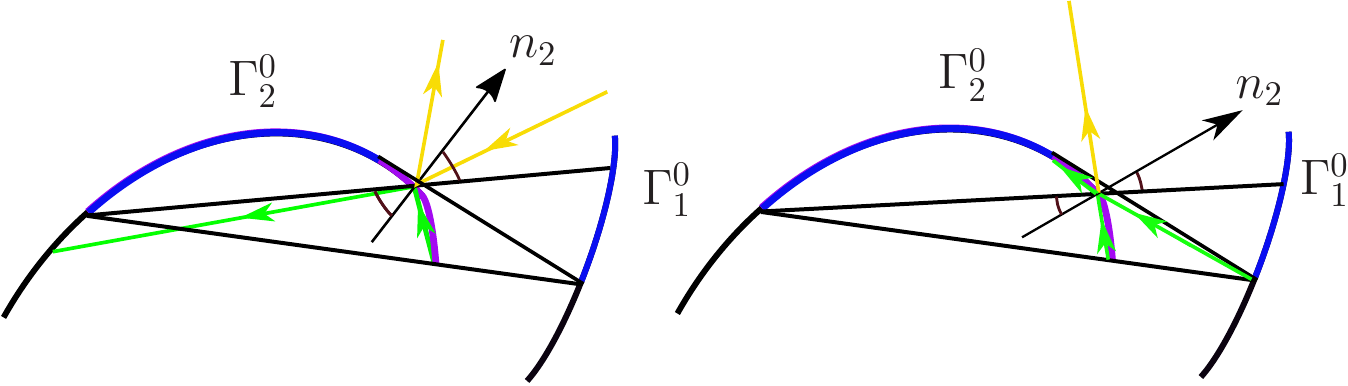}
\end{center}
\caption{\footnotesize{Left : the incoming half-ray $\gamma_2^+(\rho_{\F_2,2})$ as an angle of incidence smaller or equal than $\theta'$. The two half-rays in $\Omega_1$ therefore intersect$\Gamma_1^{n-1}$. Right : the incoming half-ray $\gamma_2^+(\rho_{\F_2,2})$ as an angle of incidence smaller or equal than $\theta'$. In this case, only one half-ray of $\Omega_1$ intersects $\Gamma_1^{n-1}$ and the outgoing half-ray of $\Omega_2$ is restricted to the region bounded by $\hat{\Gamma}_2^{n}$ and the line $l(\delta_2(s_2^{2,n-1}),\Pi_x(\F_2(x,\xi_2)))$.}}\label{n2Gamma1fds}
\end{figure}

Consider the line $l(\delta_2(s_2^{2,n-1}),\Pi_x(\F_2(x,\xi_2)))$. By construction, this line intersects $\Gamma_1^{n-1}$. Consider $\theta_{\F_2,2}$ the angle of incidence (and therefore the angle of reflection) at $\Pi_x(\F_2(x,\xi_2))$. Again, the analysis relies on a dichotomic argument. If $\theta_{\F_2,2}$ is less or equal than the angle $\theta'$ made between $-n_2(\Pi_x(\F_2(x,\xi_2)))$ and the line $l(\delta_2(s_2^{2,n-1}),\Pi_x(\F_2(x,\xi_2)))$, the half-ray $\gamma_2^+(\rho_{\F_2,2})$ intersects $\d \Omega_2 \setminus \Gamma_2^{n}$. But in such case, since $\theta_{\F_2,1}<\theta_{\F_2,2}$ where $\theta_{\F_2,1}$ is the angle between $\gamma_1^{\pm}(\rho_{\F_2,1})$ and $n_2(\Pi_x(\F_2(x,\xi_2)))$, both half-rays $\gamma_1^{\pm}(\rho_{\F_2,1})$ intersects $\Gamma_1^{n-1}$ in a non-diffractive way (see figure \ref{n2Gamma1fds}, on the left). Therefore we conclude by Corollary \ref{equivsuppmesgraph} that $\mu_1,\mu_2=0$ near $(\rho_1,\rho_2)$. Otherwise, only $\gamma_1^+(\rho_{\F_2,1})$ intersects $\Gamma_1^{n-1}$ and $\gamma_2^+(\rho_{\F_2,2})$ intersects $\hat{\Gamma}_2^{n}$ at a point $\Pi_x(\F^2_2(x,\xi_2))$ (see Figure \ref{n2Gamma1fds}, on the right). If $\Pi_x(\F^2_2(x,\xi_2)) \in \Gamma_2^{n-1}$, we conclude by Corollary \ref{equivsuppmesgraph} that $\mu_1,\mu_2=0$ near $(\rho_1,\rho_2)$. Otherwise, $\Pi_x(\F^2_2(x,\xi_2)) \in \hat{\Gamma}_2^{n} \setminus \Gamma_2^{n-1}$, in which case we iterate the present argument (thanks to Corollary \ref{equivsuppmesgraph} and to $\gamma_1^+(\rho_{\F_2,1})$ intersecting $\Gamma_1^{n-1}$).

We are able to conclude after a finite number of steps that $\mu_1,\mu_2=0$ near $(\rho_1,\rho_2)$ since, by hypothesis, $\Omega_2$ has no contact of infinite order with his boundary. Therefore, there exists $n\in \NN^*$ such that $\Pi_x(\F^n_2(x,\xi_2))\in \Gamma_2^{n-1}$ or $\Pi_x(\F^n_2(x,\xi_2)) \in \d \Omega_2 \setminus \hat{\Gamma}_2^n$, in which cases $\gamma_1^{\pm}(\rho_{\F^n_2,1})$ intersect $\Gamma_1^{n-1}$. 

\endpf

\textbf{Step 4 : Extension of $\Gamma_1^{n-1}$ to $\Gamma_1^{n}$}\newline

We now extend $\Gamma_1^{n-1}$ to $\Gamma_1^{n}$. Denote $s_1^{1,n}$ the largest $s \leq s_1^{1,n-1}$ such that
\[ 
\Pi_x(\F(\delta(s),-n(\delta(s))))=\delta_2(s_1^{2,n-1}).
\]
Likewise, denote $s_2^{1,n}$ the smallest $s\geq s_2^{1,n-1}$ such that  
\[
\Pi_x(\F(\delta(s),-n(\delta(s))))=\delta_2(s_2^{2,n-1}).
\] 
Then $\Gamma_1^{n}$ is the open and connected part of $\d \Omega$ such that $\Gamma_1^{n-1} \subset \Gamma_1^n$ and $\{\delta(s_1^{1,n}),\delta(s_2^{1,n}) \}= \d \Gamma_1^{n}$.

The proof of Lemma \ref{Lemtildegamma2}, with the slight modification that an half-ray can intersect $\Gamma_2^n$ gives Lemma \ref{Lemtildegammax0}. 

\begin{lem}
Neighborhoods of points $\rho_1 \in T^*(\Gamma_1\times (0,T))$ do not belong to the support of the measure $\mu_1$.
\end{lem}

\textbf{Step 5 : Iteration until $\Gamma_2^{n-1}=\Gamma_2^n$}\newline

We iterate step 3 and 4 until $\Gamma_2^{n-1}=\Gamma_2^n$ for $n\in \NN^*$. We begin by proving the following

\begin{lem}\label{extensiongamma12}
For any $\Omega$ and $\Omega_2$ as defined in Section \ref{SecState} and any $\Gamma(x_0)$ be defined as \eqref{defgammax0}, there exists $n\in \NN$ such that $\Gamma_2^{n-1}=\Gamma_2^n$.
\end{lem}

\beginpf

To prove Lemma \ref{extensiongamma12}, it suffices to prove that there are no accumulations points for the sequence $\{(\delta_2(s_1^{2,n}),\delta_2(s_2^{2,n}))\}_{n\in \NN}$, except when 
\begin{equation}\label{critstopfront}
(\delta_2(s_1^{2,n-1}),\delta_2(s_2^{2,n-1}))=(\delta_2(s_1^{2,n}),\delta_2(s_2^{2,n})) \textrm{ for some } n\in \NN^*. 
\end{equation}
Denote $\theta_{n,i}, i=1,2$ the angle between $-n(\delta(s_i^{1,n}))$ and the line $l(\delta(s_i^{1,n}),\delta_2(s_j^{2,n})), j=3-i$. Then it is easy to see that there exists $\theta_{min}>0$ such that $\theta_{n,1}+\theta_{n,2} \geq \alpha_{min}, i=1,2$ and for all $n\in \NN$ until \eqref{critstopfront} thanks to the fact that $\Gamma(x_0)$ defines a conical region of $\d \Omega\setminus \Gamma(x_0)$, that $\Omega$ and $\Omega_2$ are strictly convex and that there are uniform lower and upper bound for the distance between $\d \Omega_2$ and $\d \Omega\setminus \Gamma(x_0)$. 

\endpf

Therefore the iteration process ends for a finite $n\in \NN^*$ either because there are no longer intersections between the line $l(\delta(s_1^{1,n-1}),\delta_2(s_2^{2,n-1}))$ and $\d \Omega_2$ and the line $l(\delta(s_2^{1,n-1}),\delta_2(s_1^{2,n-1}))$ and $\d \Omega_2$ or because 
\begin{equation}\label{normalpara}
\Pi_{\xi}(\F^2(\delta(s_i),-n(\delta(s_i))))=-n(\delta(s_i)).
\end{equation}
At the end of the iteration process, one defines $\Gamma_1:=\Gamma_1^n$ and $\Gamma_2:=\Gamma_2^n$ and $\delta(s_i^1) \in \d \Gamma_1, \delta_2(s_i^2) \in \d \Gamma_2, i=1,2$ and $s_1^i <  s_2^i, i=1,2$. We highlight at this point that $\Gamma_1 \subsetneq \d \Omega$ and $\Gamma_2 \subsetneq \d \Omega_2$. However, we are able to prove Lemma \ref{Omega2gammax0}.

\begin{lem}\label{Omega2gammax0}
There exists $x_0^2\in \RR^2 \setminus \overline{\Omega}_2$ such that $\Gamma_2(x_0^2) = \Gamma_2$.  
\end{lem}

\beginpf

If $\delta(s_i^1) \in \d \Gamma_1, \delta_2(s_i^2) \in \d \Gamma_2, i=1,2$ are such that \eqref{normalpara}, $\delta'(s_i^1)=\delta_2'(s_i^2)$. Since $\Gamma(x_0) \subset \Gamma_1$, then the tangents at $\delta_2(s_i^2), i=1,2$ must intersect at a point $x_0^2 \in \RR^2 \setminus \Omega_2$ such that  $\Gamma_2(x_0^2) = \Gamma_2$. If, say $\delta_2(s_1^2)$, is such that there are no intersections between $\d \Omega_2$ and $l(\delta(s_2^2),\delta_2(s_1^2))$, then this line has an intersection with the tangent line at $\delta_2(s_2^2)$ if $\delta_2(s_2^2)$ satisfies \eqref{normalpara} or the line $l(\delta(s_1^2),\delta_2(s_2^2))$, defining in either cases $x_0^2 \in \RR^2 \setminus \overline{\Omega}_2$. 

\endpf

We then define $\Omega_1^f\subset \Omega_1$ the remaining part of the geometry bounded by $\d \Omega \setminus \Gamma_1$, $\d \Omega_2 \setminus \Gamma_2$ and the lines $l(\delta(s_i^1),\delta_2(s_i^2)), \delta(s_i^1)\in \d \Gamma_1,  \delta_2(s_i^2) \in \d \Gamma_2, i=1,2$ (see figure \ref{Defof}).\newline 

\textbf{Step 6 : Uniformly escaping geometry $\Omega_1^f$.}\newline

We prove 

\begin{lem}\label{lemunif}
Suppose $\Omega_1^f$ is a uniformly escaping geometry. Then $\mu_1,\mu_2=0$ near $(\rho_1,\rho_2)\in T^*(\Omega_1^f \times (0,T))$.
\end{lem}

When $\Omega_1^f$ is a uniformly escaping geometry, it suffices to let the bicharacteristics propagates until they exit $\Omega_1^f$ to obtain Lemma \ref{lemunif}. We recall the definition \eqref{directionmap} of $\M$. 

\beginpf

\begin{figure}[!ht]
\begin{center}
	\includegraphics[height=2cm]{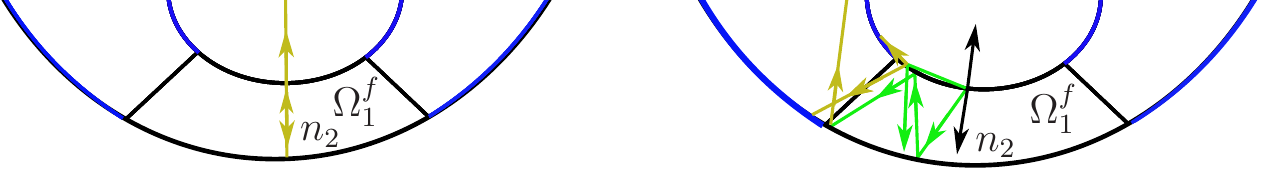}
\end{center}
\caption{\footnotesize{Left : the ray starting from $x$ in the $n_2$ direction such that $\M(x,n_2(x))=0$ is directly observed when transmitted to $\Omega_2$. Right : the propagation of the half-rays in the negative tangential direction from a point $x$ such that $\M(x,n_2(x))<0$.}}\label{figpropueg}
\end{figure}

Consider $\rho_1=(t,x,\tau_1,\xi_1) \in T^*((\d \Omega_2 \setminus \Gamma_2) \times (0,T))$ such that $\M(x,n_2(x))=0$. By definition of $\M$, $\F^2(x,n_2(x))=(x,n_2(x))$. Moreover, using Lemma \ref{Omega2gammax0}, we obtain that $\Pi_x(\F_2(x,-n_2(x))) \in \Gamma_2$ which allows to conclude that $\mu_1,\mu_2=0$ near $\rho_1=(t,x,\tau_1,n_2(x)/c_1), \rho_2=(t,x,\tau_2,-n_2(x)/c_2)$ (see figure \ref{figpropueg}, on the left). For $(\rho_1,\rho_2) \in \H^1 \times \H^2$ such that $\M(x,\xi_1)=0$ and $\xi_1\neq n_2$, each half-ray propagates toward $\M(x,\xi_1)<0$ or $\M(x,\xi_1)>0$ according to their direction of propagation. This case is therefore included in the case $(\rho_1,\rho_2) \in \H^1 \times \H^2$ such that $\M(x,\xi_1) \neq 0$.

Let $(\rho_1,\rho_2) \in \H^1 \times \H^2$ and, without loss of generality, suppose $\M(x,\xi_1) > 0$. Consider $\rho_2=(t,x,\tau_2,\xi_2) \in \H^2$ and denote the half-ray $\gamma_2^+$ outgoing from $\rho_2=(t,x,\tau_2,\xi_2)$ and moving, locally in $\sigma$, in the direction of the parametrisation. If $\Pi_x(\F(x,\xi_1))\in \Gamma_1$ and $\Pi_x(\F_2(x,\xi_2))\in \Gamma_2$, then we conclude that $\mu_1,\mu_2=0$ near $(\rho_1,\rho_2)$ by Corollary \ref{equivsuppmes}. Otherwise, $\Pi_x(\F(x,\xi_1))\in \d \Omega \setminus \Gamma_1$ or $\Pi_x(\F_2(x,\xi_2))\in \d  \Omega_2 \setminus \Gamma_2$ and one needs to carry on with the analysis (see figure \ref{figpropueg}, on the right).

We describe first the situation for $\Pi_x(\F_2(x,\xi_2))\in \d  \Omega_2 \setminus \Gamma_2$. The intersection is non-diffractive at this point and we denote $(\rho_{\F_2,1},\rho_{\F_2,2})\in \H^1 \times \H^2$. At this point we follow the outgoing rays $\gamma_1^+(\rho_{\F_2,2})$ and $\gamma_2^+(\rho_{\F_2,2})$ in the direction of propagation of $\gamma_2^+(\rho_2)$ (or equivalently $\gamma_2^-(\rho_{\F_2,2})$). It is important to notice that the direction of propagation of $\gamma_2^+(\rho_2)$  happens to move locally in the direction of the parametrisation but globally in the opposite direction, we follow the outgoing ray $\gamma_2^+(\rho_{\F_2,2})$ and the ray $\gamma_1^+(\rho_{\F_2,2})$, moving, locally, in the same direction than $\gamma_2^+(\rho_{\F_2,2})$. From Lemma \ref{Omega2gammax0}, we conclude that there exists $n^*\in \NN$ such that $\Pi_x(\F_2^{n^*}(x,\xi_2))\in \Gamma_2$.  

Let us now describe the situation for $\Pi_{x}(\F(x,\xi_1))\in \d  \Omega \setminus \Gamma_1$. Since we are considering the case $\M(x,\xi_1) > 0$, $\gamma_1^+$ move, not only locally but globally, in the direction of the parametrisation. If $\Pi_{x}(\F^2(x,\xi_1))\in \d  \Omega_2 \setminus \Gamma_2$, then the (possibly) transmitted ray of $\Omega_2$ falls in the description of the previous paragraph and the reflected ray $\gamma_1^+(\rho_{\F^2(x,\xi_1),1})$ propagate in the direction of the parametrisation. Otherwise, if $\Pi_{x}(\F^2(x,\xi_1))\in \d  \Omega \setminus \Gamma_1$, then the monotonicity of $\M$ ensures that the outgoing ray $\gamma_1^+(\rho_{\F^2(x,\xi_1),1})$ propagates in the direction of the parametrization. We conclude that there exists $n'\in \NN$ such that $\Pi_x(\F^{n'}(x,\xi_1)) \in \Gamma_1$. 

We are then able to conclude thanks to Corollary \ref{equivsuppmesgraph}. Indeed, for $(\rho_1,\rho_2)\in \H^1 \times \H^2$, then one has to follow the ray given by the previous description according to the sign of $\M$. When the ray intersects $\d \Omega_2 \setminus \Gamma_2$, one has to follow the new rays $\gamma_1^+$ or $\gamma_2^+$ which propagates according to the sign of $\M$. There exist a uniform time of observability for these rays from Lemma \ref{Omega2gammax0}, the monotonicity of $\M$ and the assumption that $\d \Omega$ and $\d \Omega_2$ has no contact of order $k-1$ with its tangents. We underline that $\rho_2 \in \G^{2,+} \cup \E^2$ does not affect the validity of the propagation of the rays. Therefore, $\mu_1,\mu_2=0$ near $(\rho_1,\rho_2)\in T^*(\d \Omega \setminus \Gamma_1) \times T^*(\d \Omega_2 \setminus \Gamma_2)$.  

\endpf

We prove Theorem \ref{thmtrap}

\beginpf

We follow the proof of Theorem \ref{main} until Step 6.  \newline

\textbf{Step 6 : non-uniformly escaping geometry $\Omega_1^f$} \newline 

The remaining case is the event where $\Omega_1^f$ is a trapping region. In such case, $c_1, c_2$ and $\Omega_1^f$ are assumed to satisfy the following : for every $\rho=(t,x,\tau,\xi) \in T^*((\d \Omega \setminus \d \Omega_1) \times (0,T))$, $\Pi_x(\F(x,\xi))\in \d \Omega_2 \setminus \Gamma_2$ then $\F_2^{\pm}(\F(x,\xi))\in \Gamma_2$. Therefore, $\mu_1=0$ near $\rho$ and $\mu_1,\mu_2=0$ near $(\rho_1,\rho_2) \in T^*(\d \Omega \setminus \Gamma_1) \times T^*(\d \Omega_2 \setminus \Gamma_2)$ as $\F^{-1}(x,\xi_1) \in \Gamma_1$ for $\rho_1$ not described by the previous analysis.    

\endpf

\section{GCC is more general than $\Gamma(x_0)$, even if $\Omega$ is strictly convex and $\Gamma$ connected}\label{SecGCCmoregen}

Consider $\Omega$ to be an ellipse of focii $F_1$ and $F_2$. We assume that the center of the ellipse is on $(0,0)$ and that the major axis lies on the $x$ axis. We denote the focii with the natural notations $F_1=(c,0)$ and $F_2=(-c,0), c>0$, the endpoints of the major axis $(-a,0)$ and $(a,0)$ and the endpoints of the minor axis $(0,-b)$ and $(0,b)$, with $a,b>0$ and $c=\sqrt{a^2-b^2}$. Let $x_1 \in \d \Omega$ (without loss of generality we assume $x_1=(x_1^1,x_1^2)$ such that $x_1^1,x_2^2<0$). Then, consider $x_2\in \d \Omega, x_1\neq x_2$ to be the only point intersected by the line $l(x_1,F_1)$. Then, 

\begin{lem}\label{lemel}
Every $\Gamma \subset \d \Omega$ open and connected including $x_1$ and $x_2$ satisfy GCC. 
\end{lem}

It is easy to see that if $\textrm{Lenght}(\Gamma)$ is stricly larger than half the perimeter of $\Omega$, then $\Gamma$ may be written of the form $\Gamma(x_0)$ for $x_0 \in \RR^2 \setminus \overline{\Omega}$. The case where $\textrm{Lenght}(\Gamma)$ is less or equal than the perimeter of $\Omega$ is the new part of the proof. The proof relies on the dynamic of the billiards. Proof of these dynamics may be found in \cite{Billiards}. \newline

\beginpf

\textbf{Major and minor axis}\newline

The bicharacteristics travelling along the major and minor axis are periodic ones, as they go back and forth on these axis. Since one $(a,0), (-b,0) \in \Gamma$, these rays are observed.\newline

\textbf{Bicharacteristics going through one foci}\newline

\begin{figure}[!ht]
\begin{center}
	\includegraphics[height=4cm]{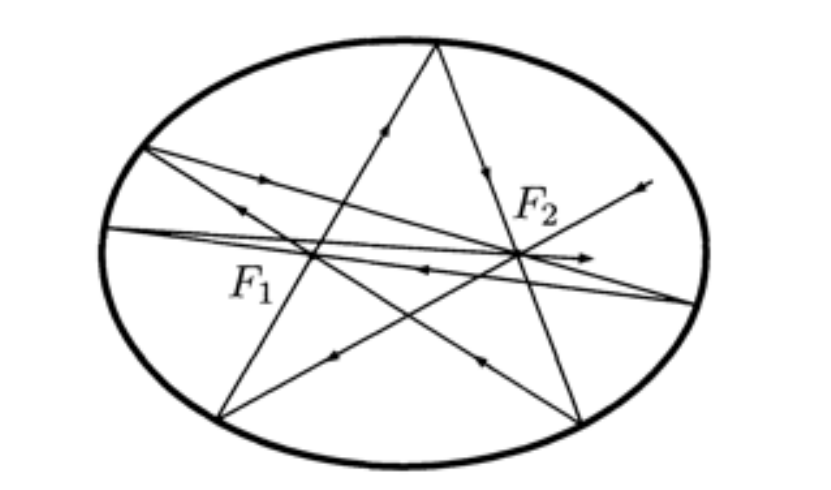}
\end{center}
\caption{\footnotesize{Example of a ray converging to the main axis (\cite{Billiards}).}}\label{convfoci}
\end{figure}

A bicharacteristic going through one foci bounces off $\d \Omega$ and go through the other foci. Every such bicharacteristics converge to the major axis and are therefore observed by $\Gamma$ (see figure \ref{convfoci}). \newline

\textbf{Hyperbolic caustics}\newline

\begin{figure}[!ht]
\begin{center}
	\includegraphics[height=4cm]{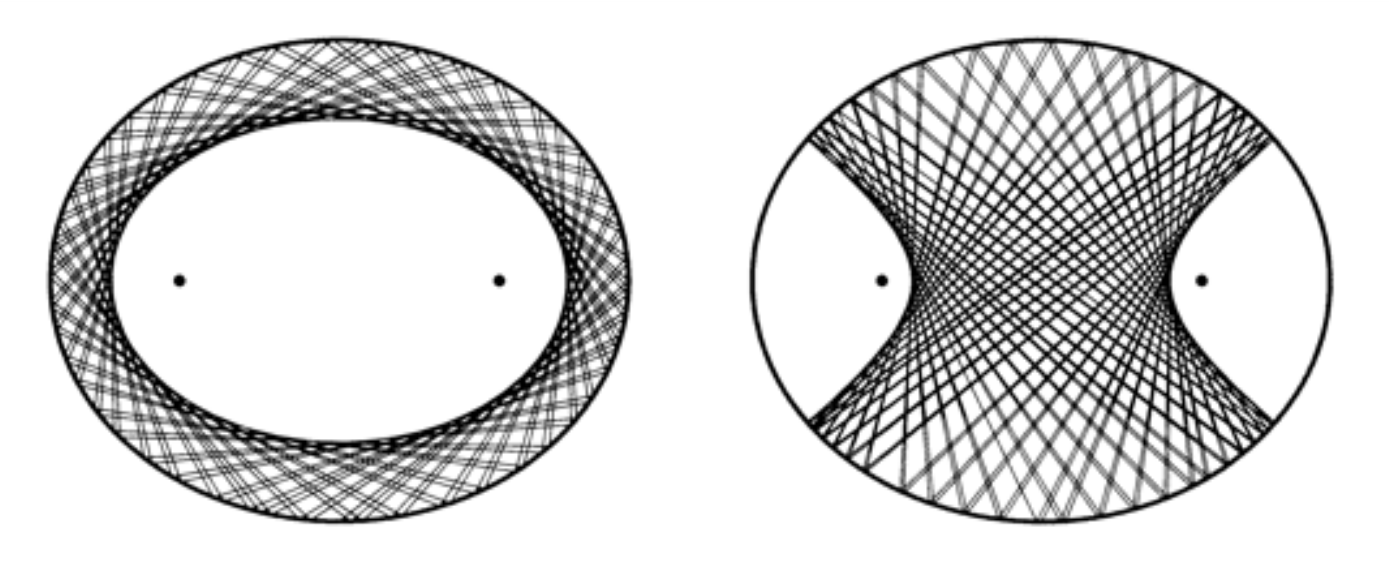}
\end{center}
\caption{\footnotesize{Left : examples of elliptic caustics. Right : examples of hyperbolic caustics (\cite{Billiards}).}}\label{caustic}
\end{figure}

Every bicharacteristics that cross the line between $F_1$ and $F_2$ will cross this line again after bouncing off $\d \Omega$. Such a bicharacteristic has all his rays tangential to an hyperbol (see figure \ref{caustic} on the right). Since $\{(x,y) \in \d \Omega \, | \, x\geq 0, y \leq 0, \} \subset \Gamma$, it is easy to see that every such bicharacteristics have to intersect $\Gamma$ in some time $T>0$. \newline

\textbf{Elliptic caustics}\newline

Every bicharacteristics that do not cross the line between $F_1$ and $F_2$ and the focii won't cross the line again after bouncing off $\d \Omega$. Moreover, each such rays is tangential to an ellipse of focii $F_1$ and $F_2$ with smaller minor and major axis (see figure \ref{caustic} on the left). Notice first that every elliptic caustic starting of $\d \Omega \cap \{y\leq 0\}$ and crossing the line $l(F_1,(a,0))$ either starts from $\Gamma$ or intersects $\Gamma$. Then, it is easy to see that every bicharacteristics starting from $\d \Omega \cap \{y\geq 0\}$ and crossing the line $l(F_2,(-a,0))$ has to intersect $\Gamma$. \newline

\textbf{Gliding ray}\newline

The gliding ray intersects $\Gamma$ in some time $T>0$.

\endpf

\begin{remark}
Notice that when $F_1 \rightarrow F_2$, one recovers that GCC is equivalent to $\Gamma(x_0)$ for the circle if $\Gamma$ is connected.
\end{remark}

\section{Conclusion}

The use of domains in $\RR^2$ was crucial in our analysis. In higher dimension, the proof does not follow immediately. Indeed, one uses in step 3 of the proof of Theorem \ref{main} the fact that the transmitted rays in $\Omega_2$ are confined in a certain region. This ensures that the rays propagating in this region either intersects $\Gamma_2^1$ or exit this region, in which case we obtain observability. This is no longer the case in dimension 3. One can only prove that the rays move toward $\Gamma_2^1$ but one can construct examples where the rays intersect regions of $\d \Omega_2 \setminus \d \Gamma_2^1$ before reaching $\Gamma_2^1$. One can circumvent this issue if one considers domains in $\RR^3$ obtained by the revolution around an axis of strictly convex curves and by considering $\Gamma(x_0)$ such that $x_0$ lies along the axis of revolution. The symmetry allows one to use the proof of Theorem \ref{main} straightforwardly. In the general case however one needs to analyse the rays transmitted back to $\Omega_1$ through $\d \Omega_2 \setminus \d \Gamma_2^1$. 

\bibliographystyle{plain}
\bibliography{biblio}

\end{document}